\documentclass[11pt,reqno]{NumPDEsArticle}

\usepackage[utf8]{inputenc}
\usepackage[T1]{fontenc}
\usepackage[english]{babel}

\usepackage{NumPDEsMacros}

\usepackage{csquotes}
\usepackage{enumitem}
\usepackage{amsmath,amssymb}
\usepackage{mathrsfs}
\usepackage{empheq}
\usepackage{graphicx}
\usepackage{tikz}
\usepackage{pgfplots}
\usepackage{caption}
\usepackage{subcaption}
\usepackage[capitalize,nameinlink]{cleveref}
\usepackage{float}
\usepackage{booktabs}
\usepackage[width=\textwidth]{caption}
\usepackage{bm}
\usepackage[dvipsnames]{xcolor}

\newcommand{\si}{\boldsymbol{s}}
\newcommand{\dtau}{\d_\tau}
\newcommand{\Dtau}{\mathrm{D}_\tau}
\newcommand{\aaa}{\boldsymbol{a}}
\newcommand{\bb}{\boldsymbol{b}}
\newcommand{\ee}{\boldsymbol{e}}
\newcommand{\mm}{\boldsymbol{m}}

\newcommand{\ff}{\boldsymbol{f}}
\newcommand{\jj}{\boldsymbol{j}}
\newcommand{\uu}{\boldsymbol{u}}
\newcommand{\vv}{\boldsymbol{v}}
\newcommand{\ww}{\boldsymbol{w}}
\newcommand{\qq}{\boldsymbol{q}}
\newcommand{\xx}{\vec{x}}
\newcommand{\zz}{\vec{z}}
\newcommand{\js}{\boldsymbol{J}_s}
\newcommand{\pp}{\boldsymbol{\pi}}
\newcommand{\fra}{\mathfrak{a}}
\newcommand{\heff}{\boldsymbol{h}_{\rm eff}}
\newcommand{\phii}{\boldsymbol{\varphi}}
\newcommand{\pphi}{\boldsymbol{\phi}}
\newcommand{\psii}{\boldsymbol{\psi}}
\newcommand{\zetaa}{\boldsymbol{\zeta}}
\newcommand{\etaa}{\boldsymbol{\eta}}
\newcommand{\bL}{\mathbf{L}}
\newcommand{\bH}{\mathbf{H}}
\newcommand{\bW}{\mathbf{W}}
\newcommand{\III}{\boldsymbol{\mathcal{I}}}

\definecolor{TUPlotBlue}{HTML}{006699}
\definecolor{TUPlotRed}{HTML}{A92E25}
\definecolor{TUPlotGreen}{HTML}{35948B}
\definecolor{TUPlotOrange}{HTML}{E18922}
\definecolor{TUPlotPurple}{HTML}{BA4682}
\definecolor{TUPlotYellow}{HTML}{D8C35A}
\definecolor{TUPlotGrey}{HTML}{9D9D9C}
\definecolor{TUPlotCyan}{HTML}{72ADD5}
\definecolor{TUPlotBrown}{HTML}{9B7140}
\definecolor{TUPlotPink}{HTML}{E4A9C7}

\DeclareMathOperator{\sech}{sech}

\title{A decoupled BDF2-type integrator for the \\spin-diffusion Landau--Lifshitz--Gilbert system \\in ferromagnetic multilayers}

\author{Michele Ald\'e}
\email{Michele.Alde@asc.tuwien.ac.at}

\author{Dirk Praetorius}
\email{Dirk.Praetorius@asc.tuwien.ac.at}

\author{Agus L. Soenjaya}
\email{Agus.Soenjaya@asc.tuwien.ac.at \quad\rm (corresponding author)}

\address{TU Wien, Institute of Analysis and Scientific Computing, Wiedner Hauptstra\ss{}e 8--10, 1040 Wien, Austria}


\thanks{This research was funded by the Austrian Science Fund (FWF) projects \href{https://www.doi.org/10.55776/F65}{10.55776/F65} (SFB F65 ``Taming complexity in PDE systems''), and \href{https://www.doi.org/10.55776/I6802}{10.55776/I6802} (international project I6802 ``Functional error estimates for PDEs on unbounded domains''). Additionally, the \href{https://www.vsmath.at}{Vienna School of Mathematics} supports Michele Ald\'e, and the Australian Mathematical Society Lift-off Fellowship supports Agus L. Soenjaya}

\begin{document}
	
\maketitle

\begin{abstract}
	We propose and analyze a fully discrete, linear, decoupled BDF2-type finite element integrator for the spin-diffusion Landau--Lifshitz--Gilbert (SDLLG) system arising in spintronics. The SDLLG equations couple the magnetization dynamics with a diffusion equation for spin accumulation. The proposed numerical method combines a projection-free BDF2-type tangent-plane discretization of the magnetization dynamics with BDF2 time-stepping for the spin accumulation. By treating the coupling terms explicitly through extrapolation, each time step requires only two independent linear solves that can be performed in parallel. The analysis accommodates the usual exchange field and a general, possibly nonlinear, lower-order effective field. Unconditional weak convergence towards a global weak solution is established on shape-regular triangulations, without imposing
	global quasi-uniformity or mesh-angle conditions. We also establish bounds on the deviation from the unit-length constraint. Under additional assumptions on the associated energy functional, an energy inequality for the limiting solution is obtained. Numerical experiments support the theoretical findings.
\end{abstract}


\section{Introduction}\label{sec:introduction}

Spin-dependent electron transport in ferromagnetic-nonmagnetic multilayers forms the basis of modern \emph{spintronics}. A landmark development was the discovery of giant magnetoresistance by Fert and Gr\"unberg, for which they were jointly awarded the 2007 Nobel Prize in Physics. This phenomenon describes the strong dependence of the electrical resistance of magnetic multilayers on
the relative orientation of their magnetizations. The diffusion-based theory in~\cite{ValFer93} subsequently provided a description of spin accumulation and transport in multilayers. A complementary manifestation of spin-dependent transport is spin-transfer torque, through which a spin-polarized electric current can induce magnetization dynamics and switching. This mechanism is central to spin-transfer-torque magnetoresistive random-access memory (STT-MRAM), a promising non-volatile magnetic memory technology~\cite{BendraEtAl2026,FiorentiniEtAl2022,WorledgeHu2024}, in which information is stored via the relative orientations of the magnetizations in magnetic layers.

At the continuum level, magnetization dynamics are commonly described by the \emph{Landau--Lifshitz--Gilbert} (LLG) equation, in which spin-polarized currents may be incorporated through prescribed spin-transfer torque terms~\cite{Slo96,ZhaLevFer02,ZhaLi04}. However, such descriptions do not adequately resolve the transport and relaxation of electron spins across multilayers. A more comprehensive description is provided by the coupling of the LLG equation to a diffusion equation for the non-equilibrium spin accumulation~\cite{ShpLevZha03,GarWan07,AbeRug_etal15,AbeRug_etal16}. The resulting \emph{spin-diffusion Landau--Lifshitz--Gilbert} (SDLLG) system captures the mutual interaction between magnetization dynamics and spin transport: the magnetization influences the spin accumulation, which in turn exerts a torque on the magnetization. This coupling is particularly relevant to multilayer structures, where spin accumulation generated at magnetic interfaces can mediate spin-transfer effects between separated ferromagnetic layers.

In this work, we consider the time-dependent SDLLG system in a bounded conducting domain $\Omega\subset\R^d$, $d=1,2,3$, containing a possibly disconnected ferromagnetic subdomain $\omega\subset\Omega$. For a final time $T>0$, we write $\omega_T\coloneqq \omega\times (0,T)$ and $\Omega_T\coloneqq \Omega\times (0,T)$. We denote by
\[
	\Gamma_J\coloneqq\partial\Omega\cap\partial\omega
\]
the part of the boundary of the ferromagnetic region that coincides with the boundary of the conducting domain. The unknowns are the magnetization $\mm:\omega_T\to\R^3$, subject to the pointwise constraint
$|\mm|=1$, and the spin accumulation $\si:\Omega_T\to\R^3$. In non-dimensional variables~\cite{GarWan07}, the coupled SDLLG system reads~\cite{AbeHrk_etal14,AbeRug_etal15,AbeRug_etal16}:
\begin{subequations}\label{equ:SDLLG}
	\begin{align}
		\partial_t\mm
		&=
		-\mm\times
		\left(
		\heff(\mm)
		+
		c\si
		\right)
		+
		\alpha \mm\times\partial_t\mm
		&&\text{in }\omega_T,
		\label{equ:SDLLG m}
		\\
		\partial_t\si
		&=
		-\nabla\cdot
		\js
		-
		D_0\si
		-
		D_0(\si \times\mm)
		&&\text{in }\Omega_T,
		\label{equ:SDLLG s}
		\\
		\js
		&=
		\beta \mm\otimes \jj 
		- 
		D_0 \left[ \nabla \si- \beta\beta' \mm \otimes \left(\nabla \si\cdot \mm\right) \right]
		&&\text{in }\Omega_T,
		\label{equ:SDLLG Js}
		\\
		\partial_{\boldsymbol{\nu}}\mm
		&= \mathbf{0}
		&&\text{on }\partial\omega\times(0,T),
		\label{equ:SDLLG m boundary}
		\\
		\partial_{\boldsymbol{\nu}} \si
		&=\mathbf{0}
		&&\text{on }\partial\Omega\times(0,T),
		\label{equ:SDLLG s boundary}
		\\
		\mm(0)
		&=\mm_0
		&&\text{in }\omega,
		\label{equ:SDLLG m initial}
		\\
		\si(0)
		&=\si_0
		&&\text{in }\Omega,
		\label{equ:SDLLG s initial}
	\end{align}
\end{subequations}
Here, $\heff(\mm)=\const{C}{ex}\Delta\mm+\pp(\mm)+\ff$ denotes the effective magnetic field, $\js$ is the spin-current tensor, and $\jj$ is the prescribed electric current density. The constants $\alpha>0$ and $c>0$ represent the Gilbert damping and spin-coupling strengths, respectively, while $D_0>0$ is the spin-diffusion coefficient and $\beta,\beta'$ are spin-polarization parameters. It is implicitly assumed that $\mm=0$ in the nonmagnetic but conducting material $\Omega\setminus \overline{\omega}$.

\subsection{Literature review}

The mathematical analysis of spin-polarized transport coupled to magnetization dynamics was initiated in~\cite{GarWan07,GarCerWan_dcds_07}, where the authors proposed the model and established the existence of global weak solutions. The existence of a global smooth solution to the SDLLG system was shown for $d=1$ in~\cite{PuGuo10} and, for small data and periodic boundary conditions, for $d=2$ in~\cite{GuoPu08}. More recently, \cite{CheWan23} established the local existence of arbitrarily regular solutions, while \cite{PuWan20,Ron23D,Ron23} proved the existence of partially regular solutions except at finitely many singular points. Optimal control of the SDLLG system, with the electric current density as the control variable, was studied in~\cite{AnMajProTra22}.

A significant development in the numerical analysis of SDLLG was achieved in~\cite{AbeHrk_etal14}, where the authors proposed a decoupled, linear finite element method based on the tangent-plane scheme. Their method employs a projection-free magnetization update and requires only two successive linear solves per time step. They established unconditional convergence towards global weak solutions of SDLLG, with an energy inequality for the limiting solutions; see also~\cite{AbeRug_etal15} for numerical simulations of the spin-diffusion model. To reduce the computational cost of the coupling, a stationary spin-diffusion approximation was considered in~\cite{RugAbeHrk_etal16} and a finite element integrator was proposed to solve this reduced problem. The validity of this stationary approximation was later investigated analytically in~\cite{FraJunPraSla23}.

For the LLG equation without spin diffusion ($\si=\boldsymbol{0}$), the development of structure-preserving finite element methods has a longer history. Fundamental approaches include the implicit midpoint method of Bartels and Prohl~\cite{BarPro06} and the linear tangent-plane scheme of Alouges and Jaisson~\cite{AloJai06,Alo08}. Higher-order methods employing BDF-type integrators and the tangent-plane scheme were subsequently developed in~\cite{AloKriSteTou14,AkrFeiKovLub21,AldFeiPra26,AldFeiPraSoen26_part3}. These developments provide the numerical foundation for the BDF2 approach employed in the present work, but do not address its coupling with the time-dependent spin accumulation.

\subsection{Main contributions}

In this work, we propose and analyze a fully discrete, linear, projection-free BDF2 finite element scheme for the time-dependent SDLLG system~\eqref{equ:SDLLG}. The method combines a linearly implicit BDF2-type integrator with a tangent-plane discretization of the LLG equation. Inspired by the fully decoupled scheme for the spin-diffusion Landau--Lifshitz--Bloch equation~\cite{Soe26}, we treat the coupling terms explicitly using extrapolated quantities from previous time levels. Consequently, each time step requires only the solution of two independent linear systems, which can be performed in parallel. This differs from~\cite{AbeHrk_etal14}, where the magnetization is updated first and subsequently used to compute spin accumulation. 

Moreover, our method updates the magnetization without nodal renormalization, thus avoiding mesh-angle restrictions associated with projection-based methods~\cite{Bar05}.
The discrete magnetization does not satisfy the unit-length constraint exactly, but its violation can be controlled. Nodal normalization is used only in the spin-diffusion part of the system to ensure the coercivity of its anisotropic diffusion operator.

Our main theoretical result establishes the unconditional subsequential convergence of the proposed scheme towards a global weak solution of the SDLLG system.
{This result requires only shape-regularity of the triangulations, thereby relaxing the global quasi-uniformity assumption in~\cite{AbeHrk_etal14}.}
To this end, we derive uniform stability estimates for both the discrete magnetization and spin accumulation, together with quantitative bounds on the violation of the nodal unit-length constraint. The convergence analysis presents several challenges. In particular, the BDF2 multistep structure complicates the derivation of discrete energy estimates, while the explicit treatment of the coupling requires careful
control of extrapolated quantities. Moreover, the magnetization enters the principal part of the
spin-diffusion operator, while its discrete approximation is not exactly of unit length. We overcome these difficulties by exploiting the BDF2 energy structure and the $G$-stability identity~\cite{Dah78,NevOde81}, as well as the coercivity of the spin-diffusion operator. Suitable compactness arguments are employed to identify the limits of the nonlinear coupling terms.

Under additional assumptions, we establish a micromagnetic energy inequality
for the limiting weak solution. In particular, our analysis allows nonlinear lower-order effective fields admitting an associated energy functional satisfying a suitable chain-rule property, rather than requiring linearity and self-adjointness as in the earlier energy analysis~\cite{AbeHrk_etal14};
see also~\cite{AldFeiPraSoen26_part3}.

The corresponding \emph{a-priori} analysis will be addressed in a subsequent work, creating a complete picture of a temporally second-order unconditional and linear fully discrete scheme for the SDLLG system.

\subsection{Outline of the paper}

The remainder of the paper is organized as follows.
Section~\ref{sec:notation and prelim} introduces the relevant notation and preliminary results. Section~\ref{sec:model} describes the SDLLG model, the relevant assumptions, the notion of weak solutions, the scheme, and the main result.
Section~\ref{sec:stability conv} is devoted to the stability and convergence analysis. Section~\ref{sec:energy} proves that the constructed weak solution satisfies a weak counterpart of the energy law.
Finally, Section~\ref{sec:num-exp} presents numerical experiments illustrating
the accuracy and performance of the proposed method.


\section{Notation and preliminaries}\label{sec:notation and prelim}

Throughout this work we use standard notation for Lebesgue, Sobolev, and Bochner spaces. 
For a measurable set $\mathcal{O}\subset\mathbb{R}^d, p \in [1,\infty]$, and $s\geq 0$, we denote by $L^p(\mathcal O)$ and $W^{s,p}(\mathcal O)$ the usual Lebesgue and Sobolev spaces, respectively, and write $H^s(\mathcal O)\coloneq W^{s,2}(\mathcal O)$. Bold letters 
denote spaces of vector- or matrix-valued functions, namely
\[
\bL^p(\mathcal{O}) \coloneq L^p(\mathcal{O};\R^3),
\qquad 
\bW^{s,p}(\mathcal{O}) \coloneq W^{s,p}(\mathcal{O};\R^3),
\qquad
\bH^s(\mathcal O) \coloneq H^s(\mathcal O; \R^3).
\]
We write $\norm[\big]{\cdot}_X$ and $|\cdot|_X$ for the norm and seminorm on a normed space $X$, respectively. In particular, for integer $k\ge 1$, the Sobolev seminorm involves derivatives
of order $k$, e.g., $|\uu|_{\mathbf W^{k,p}(\mathcal O)}= \norm[]{D^k\uu}_{L^p(\mathcal O)}$,
where the norm on the right-hand side is understood componentwise with the Euclidean matrix norm.
We write $\widetilde{\bH}^{-1}(\mathcal{O}) := \bH^1(\mathcal{O})^\ast$ for the dual space of $\bH^1(\mathcal{O})$
and denote the corresponding duality pairing by $\dual{\cdot}{\cdot}_{\widetilde{\bH}^{-1}(\mathcal{O}),\bH^1(\mathcal{O})}$, to be understood 
in the sense of the Gelfand triple
\[
\bH^1(\mathcal{O})\subset\bL^2(\mathcal{O})\subset\widetilde{\bH}^{-1}(\mathcal{O}).
\]
For readability, the duality pairing between $\widetilde{\bH}^{-1}(\mathcal{O})$ and $\bH^1(\mathcal{O})$ is often written without subscripts when the context is clear.
In particular, we shall use 
\begin{align}\label{equ:duality L2}
	\dual{\uu}{\zetaa}_{\widetilde{\bH}^{-1}(\mathcal{O}),\bH^1(\mathcal{O})}
	=
	\scalarproduct{\uu}{\zetaa}_{\mathcal{O}}\qquad \text{ for all }\uu\in \bL^2(\mathcal{O}), \;\zetaa\in \bH^1(\mathcal{O}),
\end{align}
as well as the interpolation inequality
\begin{align}\label{equ:Hminus1 interpolation}
	\norm[\big]{\ww}_{\bL^2(\mathcal{O})}^2
	\leq
	\norm[\big]{\ww}_{\widetilde{\bH}^{-1}(\mathcal{O})}
	\norm[\big]{\ww}_{\bH^1(\mathcal{O})}
	\qquad\text{for all }\ww\in\bH^1(\mathcal{O}).
\end{align}
For a Banach space $X$ and $T>0$, we use the standard Bochner spaces $L^p(0,T;X)$, $W^{s,p}(0,T;X)$, and $C([0,T];X)$ equipped with their natural norms.

For vector- or matrix-valued functions $\uu$ and $\vv$ of compatible dimensions, we denote by
\[
    (\uu,\vv)_{\mathcal{O}} \coloneq \int_{\mathcal{O}} \uu\cdot\vv\,\d x
\]
the corresponding $\bL^2(\mathcal{O})$ inner product, where the product inside the integral denotes the
Euclidean scalar product for vectors and the Frobenius product for matrices.
If $\Gamma$ is a measurable subset of the boundary of a bounded Lipschitz domain, we analogously write
\[ 
	(\uu,\vv)_{\Gamma}:=\int_{\Gamma}\uu\cdot\vv\,\d s
\]
for the corresponding surface integral. 

We shall use the continuity of the trace operator from $\bH^1(\mathcal{O})$ to $\bL^2(\partial\mathcal{O})$, namely there exists a constant $C_{\rm tr}>0$ depending only on $\mathcal{O}$, such that
\begin{align}\label{equ:trace}
	\norm[\big]{\zetaa}_{\bL^2(\partial\mathcal{O})}
	\leq
	C_{\rm tr}
	\norm[\big]{\zetaa}_{\bH^1(\mathcal{O})}
	\qquad
	\text{for all } \; \zetaa\in \bH^1(\mathcal{O}).
\end{align}

For $\aaa,\bb\in\R^3$, we write $\aaa\otimes\bb\in\R^{3\times3}$ for the tensor product
defined by $(\aaa\otimes\bb)_{jk}:=\aaa_j\bb_k$ for all $1\leq j,k\leq3$.
For vector fields $\uu,\vv:\mathcal O\subset \R^d \to\R^3$, we define $\nabla\uu\cdot\vv\in\R^d$ componentwise by
$(\nabla\uu\cdot\vv)_j \coloneq \sum_{i=1}^3 (\partial_j u_i)v_i$ for $1\leq j\leq d$.


\section{Model problem and main result}\label{sec:model}

\subsection{Mathematical model}\label{subsec:math model}

In this section, we describe the SDLLG system~\eqref{equ:SDLLG} in more detail.
Let $\omega\subset \Omega\subset\R^d$, $d=1,2,3$, be bounded polyhedral Lipschitz domains. The domain $\Omega$ represents the entire conducting multilayer, while $\omega$ denotes its ferromagnetic part. All results in this paper hold for $d\leq 3$, but we fix $d=3$ in the analysis for ease of presentation. We denote by
$\boldsymbol{\nu}:\partial\Omega\to\mathbb{S}^{d-1}\coloneqq\{\boldsymbol{x}\in\R^d\colon|\boldsymbol{x}|=1\}$ the outward unit normal to $\Omega$.
The unknowns are the normalized magnetization $\mm:\omega_T \to \mathbb{S}^2$ describing the orientation of the magnetic moments in the ferromagnetic region, and the spin accumulation $\si :\Omega_T\to \R^3$ representing the non-equilibrium polarization of the conduction electrons throughout the multilayer. 

Let $\jj:\Omega_T \to \R^d$ denote a prescribed electric current density. The associated spin-current tensor $\js:\Omega_T \to \R^{3\times d}$ is defined by
\begin{align}\label{equ:Js}
	\js
	\coloneqq 
	\beta \mm\otimes \jj 
	- 
	D_0 \left[ \nabla \si- \beta\beta' \mm \otimes \left(\nabla \si\cdot \mm\right) \right],
\end{align}
where $D_0\in L^\infty(\Omega)$ is the renormalized diffusion coefficient and $\beta,\beta'\in (0,1)$ are the nondimensional spin-polarization parameters. This spin-current tensor enters into the equation for spin diffusion in \eqref{equ:SDLLG s}.

Besides the coupling with the spin accumulation, the magnetization dynamics is driven by
the LLG equation~\eqref{equ:SDLLG m} with effective field $\heff$, which collects the contributions of the micromagnetic energy. We assume that the effective field takes the form
\begin{align}\label{equ:heff}
	\heff(\mm):= \const{C}{ex} \Delta \mm + \pp(\mm) + \ff,
\end{align}
where $\const{C}{ex}>0$ is the exchange constant, $\ff\colon\omega_T\to \R^3$ is an applied magnetic field, and $\pp\colon \bL^2(\omega)\to \bL^2(\omega)$ collects the lower-order, possibly nonlinear, contributions to the effective field. Typical examples include the magnetostatic stray field and crystalline anisotropy, see~\cite{AldFeiPraSoen26_part3} for details.

Whenever $\pp$ is linear and self-adjoint, the effective field~\eqref{equ:heff} is the
negative variational derivative of the micromagnetic energy
\begin{align}\label{equ:energy}
	\mathcal{E}(\mm,\ff)
	\coloneqq
	\frac{\const{C}{ex}}{2}\norm[\big]{\nabla\mm}_{\bL^2(\omega)}^2
	-
	\frac12\scalarproduct{\pp(\mm)}{\mm}_{\omega}
	-
	\scalarproduct{\ff}{\mm}_{\omega} ,
\end{align}
	i.e., $\heff(\mm)=-\delta\mathcal{E}(\mm,\ff)/\delta\mm$. We stress that $\mathcal{E}$
	is the micromagnetic energy of the ferromagnet alone: the spin accumulation does not
	contribute to it, and in~\eqref{equ:SDLLG m} the total driving field is
	$\heff(\mm)+c\si$, where the additional torque $c\si$ is supplied
	by~\eqref{equ:SDLLG s} rather than derived from an energy. 
Performing algebraic computations on~\eqref{equ:SDLLG m}, noting the modulus constraint $\abs{\mm}=1$ and the resulting orthogonality property $\mm\cdot\partial_t\mm=0$, we infer that
\begin{align}\label{equ:heff property}
	\heff\cdot\partial_t\mm = \alpha |\partial_t\mm|^2 - c \si \cdot \partial_t\mm.
\end{align}
A formal application of the chain rule to \eqref{equ:energy}
combined with \eqref{equ:heff property} yields the energy law
\begin{align}\label{equ:energy law}
	\frac{\d}{\d t}\mathcal{E}(\mm(t),\ff(t))
	=
	-\alpha\norm[\big]{\partial_t\mm(t)}_{\bL^2(\omega)}^2
	-
	\scalarproduct{\partial_t\ff(t)}{\mm(t)}_{\omega}
	+
	c\scalarproduct{\si(t)}{\partial_t\mm(t)}_{\omega} ,
\end{align}
which reveals the dissipative behavior of~\eqref{equ:SDLLG m} in the absence of the spin
accumulation and for $\ff$ constant in time.

We assume that there exist constants $D^\ast,D_\ast>0$ such that
\begin{equation*}
	D_\ast \leq D_0(\boldsymbol{x})\leq D^\ast \quad \text{for a.e. } \vec{x} \in \Omega.
\end{equation*}
This, together with the fact that $1-\beta\beta'>0$, ensures coercivity of the principal part of the spin-diffusion operator in \eqref{equ:SDLLG s}; see~\cite{AbeHrk_etal14} for further details. We assume the following regularity on the initial data: $\mm_0\in \bH^1(\omega)$ with $|\mm_0|=1$ and $\si_0\in \bH^1(\Omega)$. Furthermore, we assume that $\ff\in \bL^2(\omega_T)$ and $\jj\in L^2(0,T;\bH^1(\Omega))$ are given.

We now recall the definition of a weak solution to the SDLLG system~\eqref{equ:SDLLG} from~\cite{AbeHrk_etal14}. For ease of presentation, given $\mm(t)\in \bL^\infty(\omega)$ {with $|\mm(t)| \leq 1$ a.e. in $\omega$}, we define a time-dependent bilinear form $\fra(\mm(t);\cdot,\cdot):\bH^1(\Omega)\times \bH^1(\Omega)\to \R$ by
\begin{align}\label{equ:a}
	\fra(\mm(t);\zetaa_1,\zetaa_2)
	&\coloneqq
	\scalarproduct{D_0\nabla\zetaa_1}{\nabla\zetaa_2}_\Omega 
	- \beta\beta' \scalarproduct{D_0\mm(t) \otimes \bigl(\nabla\zetaa_1 \cdot \mm(t)\bigr)}{\nabla\zetaa_2}_\omega
	\nonumber\\
	&\qquad
	+ \scalarproduct{D_0\zetaa_1}{\zetaa_2}_\Omega 
	+ \scalarproduct{D_0 \bigl(\zetaa_1\times \mm(t) \bigr)}{\zetaa_2}_\omega.
\end{align}
It is known that $\fra$ is continuous and coercive~\cite[Lemma~5]{AbeHrk_etal14}: there
exists a constant $C_{\mathfrak a}>0$, depending only on $D^*$, such that
\begin{equation*}
  |\mathfrak a(\mm(t);\zetaa_1,\zetaa_2)|
  \le C_{\mathfrak a}\lVert\zetaa_1\rVert_{ \bH^1(\Omega)}\lVert\zetaa_2\rVert_{\bH^1(\Omega)},
\end{equation*}
and, since $\lvert\mm(t)\rvert\le 1$ a.e.\ in $\omega$ and $D_0\ge 0$, it holds that
\begin{align}\label{equ:a coercive}
	\fra(\mm(t);\zetaa,\zetaa) \geq (1-\beta\beta') D_\ast \norm[\big]{\zetaa}_{\bH^1(\Omega)}^2,
\end{align}
for almost all $t\in [0,T]$.

\begin{definition}\label{def:weak solution}
	Let $\mm_0\in \bH^1(\omega)$ with $|\mm_0|=1$ a.e. in $\omega$, and $\si_0\in \bH^1(\Omega)$. The pair $(\mm,\si)$ is a weak solution of the SDLLG system~\eqref{equ:SDLLG} if:
	\begin{enumerate}[label=\rm (\roman*), ref=\rm \roman*]
		\item $\mm\in \bH^1(\omega_T)$ with $|\mm|=1$ a.e. in $\omega_T$ and $\mm(0)=\mm_0$ in the sense of traces; \label{def: weak solution i}
		\item $\si\in L^2(0,T;\bH^1(\Omega))\cap L^\infty(0,T;\bL^2(\Omega)) \cap H^{1}(0,T;\widetilde{\bH}^{-1}(\Omega))$ and $\si(0)=\si_0$ in the sense of traces; \label{def: weak solution ii}
		\item for all $\phii\in \bH^1(\omega_T)$,
		\begin{align*}
			&\int_0^T\scalarproduct{\partial_t \mm(t)}{\phii(t)}_{\omega} \d t
			+ \alpha\int_0^T\scalarproduct{\partial_t \mm(t) \times \mm(t)}{\phii(t)}_{\omega} \d t
			\\
			&=
			- \const{C}{ex}\int_0^T \scalarproduct{\nabla \mm(t) \times\mm(t)}{\nabla\phii(t)}_{\omega} \d t
			+ \int_0^T\scalarproduct{\pp(\mm(t))\times\mm(t)}{\phii(t)}_{\omega} \d t
			\\
			&\quad
			+ \int_0^T\scalarproduct{\ff(t)\times \mm(t)}{\phii(t)}_{\omega} \d t
			+ c\int_0^T \scalarproduct{\si(t)\times\mm(t)}{\phii(t)}_{\omega} \d t;
		\end{align*} 
		\label{def: weak solution iii}
		\item for almost all $t\in [0,T]$ and all $\zetaa\in \bH^1(\Omega)$,
		\begin{align*}
			& \dual{\partial_t \si(t)}{\zetaa}_{\widetilde{\bH}^{-1}(\Omega),\bH^1(\Omega)}
			+ \fra(\mm(t);\si(t),\zetaa)
			=
			\beta \scalarproduct{\mm(t)\otimes \jj(t)}{\nabla \zetaa}_{\omega}
			- \beta \scalarproduct{\jj(t) \cdot \boldsymbol{\nu}}{\mm(t)\cdot \zetaa}_{\Gamma_J}.
		\end{align*}
		\label{def: weak solution iv}
	\end{enumerate}
\end{definition}


\subsection{Discretization}\label{subsec:discretization}

We now introduce the temporal discretization used throughout the paper. Let $N\in \N$ and set $\tau:= \frac{T}{N}$ and $t_n:=n\tau$ for $n=0,1,\ldots,N$. For a sequence $\{\ww^n\}_{n=0}^N$, we define the backward difference
\[
	\dtau \ww^{n+1}\coloneqq \frac{\ww^{n+1}-\ww^n}{\tau}, \qquad n=0,\ldots,N-1.
\]
For $n=1,\ldots,N-1$, set $\dtau^2 \ww^{n+1}\coloneqq \dtau(\dtau \ww^{n+1})$ and the BDF2 difference quotient 
\begin{align}\label{equ: Dtau}
	\Dtau \ww^{n+1} \coloneqq \frac{3\ww^{n+1}-4 \ww^n + \ww^{n-1}}{2\tau}
	=
	\frac32 \dtau \ww^{n+1}
	- \frac12 \dtau \ww^n.
\end{align}
We shall also use the second-order extrapolation given by
\begin{equation}\label{equ:hat w}
	\widehat{\ww}^{n+1}\coloneqq 2 \ww^n- \ww^{n-1}, \qquad n=1,\ldots, N.
\end{equation}

Next, we introduce the necessary setup for the spatial discretization. Let $\{\mathcal{T}_h^\Omega\}_{h>0}$ be a $\gamma$-{shape-regular} family of conforming simplicial triangulations of $\Omega$ with mesh-size $h:= \max_{K\in \mathcal{T}_h^\Omega} h_K$, where $h_K:=\operatorname{diam}(K)$. {In particular, $h_K/\rho_K\leq\gamma$ for all $K\in\mathcal T_h^\Omega$, where $\rho_K$ is the inradius of $K$.} We denote by $\mathcal{T}_h^\omega$ the restriction of $\mathcal{T}_h^\Omega$ to $\omega$ and assume that the ferromagnetic region $\omega$ is resolved by $\mathcal{T}_h^\Omega$, namely
\begin{align*}
	\mathcal{T}_h^\omega 
	=
	\left. \mathcal{T}_h^\Omega \right|_{\omega} 
	=
	\left\{K\in \mathcal{T}_h^\Omega: K\cap \omega\neq \emptyset\right\}
	\quad \text{and}\quad
	\overline{\omega}= \bigcup_{K\in \mathcal{T}_h^{\omega}} K.
\end{align*}
The corresponding set of nodes of $\mathcal{T}_h^\Omega$ (resp. $\mathcal{T}_h^\omega $) is denoted by $\mathcal{N}_h^\Omega$ (resp. $\mathcal{N}_h^\omega$).

For $\mathcal{O} \in \{\omega,\Omega\}$, let
\[
	\V_h^\mathcal{O} \coloneqq \left\{ \vv_h \in C(\overline{\mathcal{O}}; \R^3): \vv_h|_K \in \P_1(K;\R^3), \; \forall \,K\in \mathcal{T}_h^\mathcal{O}\right\}
\]
define the space of continuous and piecewise affine vector fields from $\overline{\mathcal{O}}$ to $\R^3$. We write~$\III_h^\mathcal{O}: C(\overline{\mathcal{O}}; \R^3) \to \V_h^\mathcal{O}$ for the corresponding nodal interpolation operator. Since the mesh resolves $\omega$, the interpolation operators are compatible with the restriction to $\omega$, thus we simply write $\III_h$.

The pointwise constraint $\abs{\mm}=1$ and the orthogonality relation $\mm\cdot \partial_t \mm=0$ motivate the following sets. For any $\phii_h\in \V_h^\omega$, define the nodal tangent space
\begin{align}\label{equ:discrete tangent sp}
	\mathcal{K}_h(\phii_h)
	\coloneqq
	\left\{ \psii_h \in \V_h^\omega: \phii_h(\vec{z}) \cdot \psii_h(\vec{z})=0, \;\forall\; \vec{z}\in \mathcal{N}_h^\omega \right\}.
\end{align}
We further introduce the sets
\begin{align*}
	\mathcal{M}_h
	\coloneqq
	\left\{
	\phii_h\in \V_h^\omega:
	\abs{\phii_h(\vec{z})}=1,
	\; \forall\; \vec{z}\in \mathcal{N}_h^\omega
	\right\},
\end{align*}
and
\begin{align*}
	\mathcal{U}_h
	\coloneqq
	\left\{
	\phii_h\in \V_h^\omega:
	\abs{\phii_h(\vec{z})}\geq 1,
	\; \forall\; \vec{z}\in \mathcal{N}_h^\omega
	\right\}.
\end{align*}
For $\phii_h\in \mathcal{U}_h$, we define the nodal projection
$\Pi_h\phii_h\in \mathcal{M}_h$ by
\begin{align*}
	\Pi_h\phii_h(\vec{z})
	\coloneqq
	\frac{\phii_h(\vec{z})}
	{\abs{\phii_h(\vec{z})}},
	\qquad
	\text{for all } \; \vec{z}\in \mathcal{N}_h^\omega.
\end{align*}
For this nodal projection~\cite{AbeHrk_etal14}, it holds that
\begin{equation}\label{equ:Pi Linfty}
	\norm[\big]{\Pi_h\phii_h}_{\bL^\infty(\omega)}
	\leq 1,
	\qquad
	\text{for all } \; \phii_h\in \mathcal{U}_h.
\end{equation}
Moreover, there exists a constant $C_\Pi>0$, depending only on $\gamma$, such that
\begin{align}\label{equ:Pi H1 stability}
	\norm[\big]{\nabla \Pi_h\phii_h}_{\bL^2(\omega)}
	\leq
	C_\Pi
	\norm[\big]{\nabla\phii_h}_{\bL^2(\omega)},
	\qquad
	\text{for all } \; \phii_h\in \mathcal{U}_h.
\end{align}
In particular, no quasi-uniformity or angle condition on the triangulations is required for~\eqref{equ:Pi H1 stability}. We shall prove later that all extrapolated magnetizations to which $\Pi_h$ is applied belong to $\mathcal{U}_h$, so that the projection is well defined throughout the discrete evolution.

We also denote by $P_h^\mathcal{O}: \bL^2(\mathcal{O}) \to\V_h^\mathcal{O}$
the $\bL^2(\mathcal{O})$-orthogonal projection, given by
\begin{align}\label{equ:L2 projection}
	\left(
	P_h^\mathcal{O}\zetaa-\zetaa,
	\zetaa_h
	\right)_\mathcal{O}
	=0,
	\qquad
	\text{for all } \; \zetaa_h\in\V_h^\mathcal{O}.
\end{align}
Clearly, $P_h^\mathcal{O}$ is stable in $\bL^2(\mathcal{O})$.
{By the best approximation property, $P_h^\mathcal O\zetaa\to\zetaa$ in $\bL^2(\mathcal O)$ for every $\zetaa\in\bL^2(\mathcal O)$ as $h\to0$.
For the analysis, let $R_h:\bH^1(\Omega)\to\V_h^\Omega$ be the Ritz projection with respect to the standard $\bH^1(\Omega)$ inner product. By the best approximation property, the density of smooth functions in $\bH^1(\Omega)$, and the standard interpolation estimates, we have
\begin{align}\label{equ:H1 projection}
    \norm[\big]{R_h\zetaa}_{\bH^1(\Omega)}\leq\norm[\big]{\zetaa}_{\bH^1(\Omega)},
\end{align}
and $R_h\zetaa\to\zetaa$ in $\bH^1(\Omega)$ for all $\zetaa\in\bH^1(\Omega)$ as $h\to 0$. Regarding $R_h$ as a bounded linear operator on $\bH^1(\Omega)$, we define its dual adjoint $R_h^*:\widetilde{\bH}^{-1}(\Omega)\to\widetilde{\bH}^{-1}(\Omega)$ by
\begin{align}\label{equ:dual reconstruction}
    \dual{R_h^*\ww}{\zetaa}\coloneq\dual{\ww}{R_h\zetaa}
    \quad\text{for all }\ww\in\widetilde{\bH}^{-1}(\Omega),\ \zetaa\in\bH^1(\Omega).
\end{align}
By~\eqref{equ:H1 projection}, we have
\[
    \bigl|\dual{\ww}{R_h\zetaa}\bigr|
    \leq\norm[\big]{\ww}_{\widetilde{\bH}^{-1}(\Omega)}\norm[\big]{\zetaa}_{\bH^1(\Omega)}.
\]
Thus, $R_h^*$ is well-defined, linear, and satisfies
$\norm[\big]{R_h^*\ww}_{\widetilde{\bH}^{-1}(\Omega)}\leq\norm[\big]{\ww}_{\widetilde{\bH}^{-1}(\Omega)}$ for all $\ww\in\widetilde{\bH}^{-1}(\Omega)$.
In particular, using~\eqref{equ:duality L2} and $R_h\ww_h=\ww_h$ for $\ww_h\in\V_h^\Omega$, we have
\begin{align}\label{equ:discrete Hminus1 interpolation}
    \norm[\big]{\ww_h}_{\bL^2(\Omega)}^2
    =\dual{R_h^*\ww_h}{\ww_h}
    \leq\norm[\big]{R_h^*\ww_h}_{\widetilde{\bH}^{-1}(\Omega)}\norm[\big]{\ww_h}_{\bH^1(\Omega)}.
\end{align}
Here and below, $\bL^2(\Omega)$ is identified with its image in $\widetilde{\bH}^{-1}(\Omega)$ as in~\eqref{equ:duality L2}.}

We will repeatedly use the {weighted} discrete norm
equivalence (cf.~\cite[Lemma~3.4]{Bar15}).
{For $\mathcal O\in\{\omega,\Omega\}$ and
	$\vec z\in\mathcal N_h^\mathcal O$, let
	\begin{align}\label{equ:nodal weights}
		\mu_{\vec z}^{\mathcal O}
		\coloneq \frac14
		\sum_{\substack{K\in\mathcal T_h^\mathcal O\\\vec z\in K}}|K|.
	\end{align}
	Standard scaling arguments show that, for every $r\in[1,\infty)$,
	there exists a constant $C_r>0$, depending only on $r$ and the
	spatial dimension, such that}
\begin{align}\label{equ:discrete nodal norm}
	C_r^{-1}
	\norm[\big]{\vv_h}_{L^r(\mathcal O)}^r
	\leq
	\sum_{\vec z\in\mathcal N_h^\mathcal O}
	{\mu_{\vec z}^{\mathcal O}}
	\abs{\vv_h(\vec z)}^r
	\leq
	C_r\norm[\big]{\vv_h}_{L^r(\mathcal O)}^r,
	\quad\text{for all }\vv_h\in\V_h^\mathcal O.
\end{align}

Finally, we set discrete initial data $(\mm_h^0, \si_h^0) \in\mathcal{M}_h \times \V_h^\Omega$, and let $\pp_h: \bL^2(\omega) \to \bL^2(\omega)$ denote a discrete realization of the lower-order field contribution $\pp$. For $n=1,\ldots,N$, let $\ff_h^n$ and $\jj_h^n$ denote suitable
discrete approximations of the prescribed fields $\ff$ and $\jj$ at time $t_n$. Our analysis requires only square integrability in time of $\ff$ and $\jj$, without pointwise-in-time evaluations. The stability and consistency assumptions on $\mm_h^0$, $\si_h^0$, $\pp_h$, $\ff_h^n$, and $\jj_h^n$ required for the
convergence analysis will be stated together with the main convergence
result in the next subsection.


\subsection{Numerical algorithm and main result}\label{sec:num alg}

For $\phii_h\in\mathcal{M}_h$ and $\si_h,\zetaa_h\in\V_h^\Omega$, we define
the discrete bilinear form
\begin{align}
	\label{equ:discrete spin bilinear form}
	\fra_h(\phii_h;\si_h,\zetaa_h)
	&\coloneqq
	\scalarproduct{D_0\nabla\si_h}{\nabla\zetaa_h}_{\Omega}
	-
	\beta\beta'
	\scalarproduct{D_0\phii_h\otimes \left(\nabla\si_h\cdot\phii_h\right)}{\nabla\zetaa_h}_{\omega}
	\nonumber
	\\
	&\qquad
	+
	\scalarproduct{D_0\si_h}{\zetaa_h}_{\Omega}
	+
	\scalarproduct{D_0\left(\si_h\times\phii_h\right)}{\zetaa_h}_{\omega}.
\end{align}
Using $\|\phii_h\|_{\bL^\infty(\omega)}\le 1$ for all
$\phii_h\in\mathcal M_h$, the discrete bilinear form is uniformly continuous:
there exists a constant $C>0$, depending only on $D^*$, $\beta$, and $\beta'$, such that
\[
    |\fra_h(\phii_h;\si_h,\zetaa_h)|
    \le
    C
    \|\si_h\|_{\bH^1(\Omega)}
    \|\zetaa_h\|_{\bH^1(\Omega)}
\]
for all $\phii_h\in\mathcal M_h$ and all
$\si_h,\zetaa_h\in \mathbb{V}_h^\Omega$.

As in~\eqref{equ:a coercive}, the discrete bilinear form $\fra_h$ satisfies the coercivity estimate
\begin{align}\label{equ:ah coercive}
	\fra_h(\phii_h;\zetaa_h,\zetaa_h) \geq (1-\beta\beta') D_\ast \norm[\big]{\zetaa_h}_{\bH^1(\Omega)}^2.
\end{align}

We propose the following BDF2-type integrator for solving the SDLLG system~\eqref{equ:SDLLG}.

\begin{algorithm}\label{alg: full discr}
	
	\textbf{Input:}
	Discrete initial data
	$\mm_h^0\in\mathcal{M}_h$ and $\si_h^0\in\V_h^\Omega$,
	discrete approximations
	$\{\ff_h^n\}_{n=1}^N$ and $\{\jj_h^n\}_{n=1}^N$
	of the prescribed fields,
	$\tau>0$, and $N\in\N$.
	
	\begin{enumerate}[label=\rm (\roman*), ref=\rm \roman*]
		
		\item\label{item:first step full discr}
		
		Compute $\vv_h^0$ and $\si_h^1$ in parallel as follows:
		
		Compute $\vv_h^0\in\mathcal{K}_h(\mm_h^0)$ such that, for all
		$\phii_h\in\mathcal{K}_h(\mm_h^0)$,
		\begin{align}
			\label{equ:first step full discr m}
			&
			\alpha
			\scalarproduct{\vv_h^0}{\phii_h}_{\omega}
			+
			\scalarproduct{\mm_h^0\times\vv_h^0}{\phii_h}_{\omega}
			+
			\const{C}{ex}\tau
			\scalarproduct{\nabla\vv_h^0}{\nabla\phii_h}_{\omega}
			\nonumber
			\\
			&=
			-
			\const{C}{ex}
			\scalarproduct{\nabla\mm_h^0}{\nabla\phii_h}_{\omega}
			+
			\scalarproduct{
				\pp_h(\mm_h^0)+\ff_h^1+c\si_h^0
			}{
				\phii_h
			}_{\omega}.
		\end{align}
		
		Compute $\si_h^1\in\V_h^\Omega$ such that, for all
		$\zetaa_h\in\V_h^\Omega$,
		\begin{align}
			\label{equ:first step full discr s}
			&
			\scalarproduct{\dtau\si_h^1}{\zetaa_h}_{\Omega}
			+
			\fra_h
			\left(
			\Pi_h(\mm_h^0);
			\si_h^1,
			\zetaa_h
			\right)
			\nonumber
			\\
			&=
			\beta
			\scalarproduct{
				\Pi_h(\mm_h^0)\otimes\jj_h^1
			}{
				\nabla\zetaa_h
			}_{\omega}
			-
			\beta
			\scalarproduct{
				\jj_h^1\cdot\boldsymbol{\nu}
			}{
				\Pi_h(\mm_h^0)\cdot\zetaa_h
			}_{\Gamma_J}.
		\end{align}
		
		\item\label{item:def mm1}
		
		Set
		\begin{equation}\label{equ:first step m1}
		\mm_h^1\coloneqq\mm_h^0+ \tau\vv_h^0.
		\end{equation}
		
		\item\label{item:next steps}
		
		For $n=1,\ldots,N-1$, repeat the following
		steps~\eqref{item:def predictors}--\eqref{item:def mmhj}:
		
		\begin{enumerate}[label=\rm(\alph*), ref=\rm \alph*]
			\item\label{item:def predictors}
			Set
			\begin{align}
				\label{equ:def predictors}
				\widehat{\mm}_h^{n+1}
				\coloneqq
				2\mm_h^n-\mm_h^{n-1},
				\qquad
				\widehat{\si}_h^{n+1}
				\coloneqq
				2\si_h^n-\si_h^{n-1}.
			\end{align}
			
			\item\label{item:BDF2 m step}
			Compute
			$\vv_h^n\in\mathcal{K}_h(\widehat{\mm}_h^{n+1})$
			such that, for all
			$\phii_h\in\mathcal{K}_h(\widehat{\mm}_h^{n+1})$,
			\begin{align}
				\label{equ:full discr alg m}
				&
				\alpha
				\scalarproduct{\vv_h^n}{\phii_h}_{\omega}
				+
				\scalarproduct{
					\widehat{\mm}_h^{n+1}\times\vv_h^n
				}{
					\phii_h
				}_{\omega}
				+
				\frac23\const{C}{ex}\tau
				\scalarproduct{\nabla\vv_h^n}{\nabla\phii_h}_{\omega}
				\nonumber
				\\
				&=
				-
				\frac13\const{C}{ex}
				\scalarproduct{\nabla\left(4\mm_h^n-\mm_h^{n-1}\right)}{\nabla\phii_h}_{\omega}
				+
				\scalarproduct{\pp_h\bigl(\widehat{\mm}_h^{n+1}\bigr) +\ff_h^{n+1} +c\widehat{\si}_h^{n+1}}{\phii_h}_{\omega}.
			\end{align}
			
			\item\label{item:BDF2 s step}
			In parallel with step~\eqref{item:BDF2 m step}, compute
			$\si_h^{n+1}\in\V_h^\Omega$ such that, for all $\zetaa_h\in\V_h^\Omega$,
			\begin{align}
				\label{equ:full discr alg s}
				&
				\scalarproduct{\Dtau\si_h^{n+1}}{\zetaa_h}_{\Omega}
				+
				\fra_h
				\bigl(
				\Pi_h\bigl(\widehat{\mm}_h^{n+1}\bigr);
				\si_h^{n+1},
				\zetaa_h
				\bigr)
				\nonumber
				\\
				&=
				\beta
				\scalarproduct{
					\Pi_h\bigl(\widehat{\mm}_h^{n+1}\bigr)
					\otimes\jj_h^{n+1}
				}{
					\nabla\zetaa_h
				}_{\omega}
				-
				\beta
				\scalarproduct{
					\jj_h^{n+1}\cdot\boldsymbol{\nu}
				}{
					\Pi_h\bigl(\widehat{\mm}_h^{n+1}\bigr)
					\cdot\zetaa_h
				}_{\Gamma_J}.
			\end{align}
			
			\item\label{item:def mmhj}
			
			Set
			\begin{equation}\label{equ:gen step mn1}
			\mm_h^{n+1}
			\coloneqq
			\frac43\mm_h^n
			-
			\frac13\mm_h^{n-1}
			+
			\frac23\tau\vv_h^n.
			\end{equation}
			
		\end{enumerate}
	\end{enumerate}
	
	\textbf{Output:}
	Sequences
	$\{\mm_h^n\}_{n=0}^N$,
	$\{\si_h^n\}_{n=0}^N$, and
	$\{\vv_h^n\}_{n=0}^{N-1}$.
\end{algorithm}

Thanks to the Lax--Milgram lemma as well as
the continuity and coercivity of $\fra_h$, Algorithm~\ref{alg: full discr} is well-posed.
Notice that, for the steps in which the nodal projections $\Pi_h(\mm^0_h)$ and
$\Pi_h(\widehat{\mm}^{\,n+1}_h)$ are involved, we need to guarantee that
$\mm^0_h,\widehat{\mm}^{\,n+1}_h\in\mathcal U_h$. Since $\mm^0_h\in\mathcal M_h$, 
the first time step is well defined. Moreover, at the
$n$-th step the extrapolant $\widehat{\mm}^{\,n+1}_h=2\mm^n_h-\mm^{n-1}_h$ is already
determined by the previous iterates, and an induction on $n$ shows that
$\widehat{\mm}^{\,n+1}_h\in\mathcal U_h$ for every $n=1,\dots,N-1$; see
Lemma~\ref{lem:admissibility predictor}, whose proof uses only the iterates computed in the previous steps.
In particular, $\widehat{\mm}^{\,n+1}_h(\zz)\neq\bm 0$ for every
$\zz\in\mathcal N^\omega_h$, so that the nodal constraint in \eqref{equ:discrete tangent sp} is non-degenerate
at every node and $\mathcal K_h(\widehat{\mm}^{\,n+1}_h)$ is a non-degenerate discrete tangent
space. The stronger uniform bound
$|\widehat{\mm}^{\,n+1}_h(\zz)|\ge 1$ is what guarantees the stability properties
\eqref{equ:Pi Linfty}--\eqref{equ:Pi H1 stability} of the nodal projection.

Using the discrete solutions generated by Algorithm~\ref{alg: full discr}, we introduce the following time reconstructions. For
$n=0,\ldots,N-1$ and $t\in[t_n,t_{n+1})$, let
\begin{align}\label{equ:time reconstructions m}
	\mm_{h\tau}(t)
	&:=
	\frac{t-t_n}{\tau}\mm_h^{n+1}
	+ \frac{t_{n+1}-t}{\tau}\mm_h^n,
	&
	\mm_{h\tau}^-(t)
	&:=\mm_h^n,
	&
	\mm_{h\tau}^+(t)
	&:=\mm_h^{n+1},
\end{align}
and
\begin{align}\label{equ:time reconstructions s}
	\si_{h\tau}(t)
	&:=
	\frac{t-t_n}{\tau}\si_h^{n+1}
	+ \frac{t_{n+1}-t}{\tau}\si_h^n,
	&
	\si_{h\tau}^-(t)
	&:=\si_h^n,
	&
	\si_{h\tau}^+(t)
	&:=\si_h^{n+1}.
\end{align}
Under the convention
\[
\widehat{\mm}_h^1:=\mm_h^0,
\qquad
\widehat{\si}_h^1:=\si_h^0,
\]
we also define the extrapolated piecewise constant functions
\begin{align}
	\widehat{\mm}_{h\tau}^+(t)
	&:=\widehat{\mm}_h^{n+1},
		\qquad t\in[t_n,t_{n+1}),\quad n=0,\ldots,N-1, \label{equ:time reconstructions predictor m}
	\\
	\widehat{\si}_{h\tau}^+(t)
	&:=\widehat{\si}_h^{n+1},
		\qquad t\in[t_n,t_{n+1}),\quad n=0,\ldots,N-1. \label{equ:time reconstructions predictor s}
\end{align}
Moreover, we define
\begin{align*}
	\vv_{h\tau}^-(t):=\vv_h^n,
	\qquad
	t\in[t_n,t_{n+1}).
\end{align*}
We also let $\Dtau\si_h^1:= \dtau \si_h^1$ and further set for $n=0,\ldots,N-1$,
\begin{align}\label{equ:time reconstruction BDF2 s}
	\Dtau\si_{h\tau}^+(t)
	:=
	\Dtau\si_h^{n+1},
	\qquad t\in[t_n,t_{n+1}).
\end{align}
Note that we have the following identities:
\begin{align}\label{equ:dtm dts cont disc}
	\partial_t\mm_{h\tau}(t)
	=
	\dtau\mm_h^{n+1},
	\qquad
	\partial_t\si_{h\tau}(t)
	=
	\dtau\si_h^{n+1},
	\qquad
	\text{for }t\in [t_n,t_{n+1}).
\end{align}
We impose the following assumptions on the discrete data and the lower-order effective field operator.

\begin{enumerate}[label=\rm (H\arabic*), ref=\rm H\arabic*]
	\item\label{ass:initial data}
	The discrete initial data satisfy
	\begin{alignat}{2}
		\mm_h^0
		&\rightharpoonup \mm_0
		&&\qquad \text{in } \bH^1(\omega),
		\label{equ:initial convergence m}
		\\
		\si_h^0
		&\rightharpoonup \si_0
		&&\qquad \text{in } \bH^1(\Omega),
		\label{equ:initial convergence s}
	\end{alignat}
	as $h\to0$. We will denote
	\[
		C_{\mm,0}:=\sup_{h>0}\,\norm[\big]{\mm_h^0}_{\bH^1(\omega)}\quad \text{ and } \quad C_{\si,0}:=\sup_{h>0}\,\norm[\big]{\si_h^0}_{\bH^1(\Omega)}
	\]
	and notice that both constants are finite by the weak convergences in \eqref{equ:initial convergence m}--\eqref{equ:initial convergence s}.
	We further assume that the discrete approximations of the prescribed data
	satisfy
	\begin{alignat}{2}
		\ff_{h\tau}^+
		&\rightharpoonup \ff
		&&\qquad \text{in } \bL^2(\omega_T),
		\label{equ:f discrete convergence}
		\\
		\jj_{h\tau}^+
		&\to \jj
		&&\qquad \text{in }
		\bL^2\bigl(0,T;\bH^1(\Omega)\bigr),
		\label{equ:j discrete convergence}
	\end{alignat}
	where
	\[
	\ff_{h\tau}^+(t):=\ff_h^{n+1},
	\qquad
	\jj_{h\tau}^+(t):=\jj_h^{n+1}
	\qquad
	\text{for }t\in[t_n,t_{n+1}).
	\]
	Since every weakly convergent sequence is bounded, there exist constants $\const{C}{$f$}>0$ and $\const{C}{$j$}>0$ independent of $h$, $\tau$, and $N$, such that
	\begin{equation}\label{equ:Cf Cj}
		\tau \sum_{n=0}^{N-1} \norm[\big]{\ff_h^{n+1}}_{\bL^2(\omega)}^2 \leq \const{C}{$f$}, 
		\qquad
		\tau \sum_{n=0}^{N-1} \norm[\big]{\jj_h^{n+1}}_{\bH^1(\Omega)}^2 \leq \const{C}{$j$}.
	\end{equation}
	
	\item\label{ass:pi bounded}
	The discrete operators $\pp_h:\bL^2(\omega)\to\bL^2(\omega)$
	are uniformly bounded, i.e., there exists $\const{C}{$\pi$}>0$ independent of $h$
	such that
	\begin{align}
		\norm[\big]{\pp_h(\ww)}_{\bL^2(\omega)}
		\leq
		\const{C}{$\pi$} \norm[\big]{\ww}_{\bL^2(\omega)}
		\qquad
		\text{for all }\; \ww\in\bL^2(\omega).\nonumber
	\end{align}
	
	\item\label{ass:pi consistency}
	For any sequence
	$\ww_{h\tau}\to\ww$ strongly in $\bL^2(\omega_T)$, it holds
	\begin{align}
		\pp_h(\ww_{h\tau})
		\rightharpoonup
		\pp(\ww)
		\qquad
		\text{in }\bL^2(\omega_T)\nonumber
	\end{align}
	as $(h,\tau)\to(0,0)$.
\end{enumerate}

We now state the main theorem of this paper.

\begin{theorem}\label{the:main convergence}
	Let $\{\mathcal{T}_h^\Omega\}_{h>0}$ be a 
	$\gamma$-{shape-regular} family of conforming simplicial triangulations of $\Omega$
	which resolves $\omega$, and suppose that
	assumptions~\eqref{ass:initial data}--\eqref{ass:pi consistency} are satisfied.
	Then, the sequences of discrete functions
	\[
		\mm_{h\tau}^\star
		\in
		\bigl\{
			\mm_{h\tau},\;
			\mm_{h\tau}^{\pm},\;
			\widehat{\mm}_{h\tau}^{+},\;
			\Pi_h\bigl(\widehat{\mm}_{h\tau}^{+}\bigr)
		\bigr\}
		\qquad\text{and}\qquad
		\si_{h\tau}^\star
		\in
		\bigl\{
			\si_{h\tau},\;
			\si_{h\tau}^{\pm},\;
			\widehat{\si}_{h\tau}^{+}
		\bigr\}
	\]
	generated by Algorithm~\ref{alg: full discr} admit subsequences (not relabeled)
	that unconditionally converge towards a pair $(\mm,\si)$ which is a weak solution of
	the SDLLG system~\eqref{equ:SDLLG} in the sense of
	Definition~\ref{def:weak solution}. More precisely, it holds that
	\begin{align}\label{equ:main convergence}
		\mm_{h\tau} \rightharpoonup \mm
		\quad\text{in } \bH^1(\omega_T)
		\qquad\text{and}\qquad
		\si_{h\tau} \rightharpoonup \si
		\quad\text{in } L^2\bigl(0,T;\bH^1(\Omega)\bigr)
	\end{align}
	as $(h,\tau)\to(0,0)$, and hence $\mm_{h\tau}\to\mm$ strongly in $\bL^2(\omega_T)$.
	Moreover, $\si_{h\tau}\to\si$ strongly in $\bL^2(\Omega_T)$,
	$\vv_{h\tau}^{-}\rightharpoonup\partial_t\mm$ in $\bL^2(\omega_T)$, as well as
	$\mm_{h\tau}^\star\overset{*}{\rightharpoonup}\mm$ in
	$L^\infty\bigl(0,T;\bH^1(\omega)\bigr)$ and
	$\si_{h\tau}^\star\overset{*}{\rightharpoonup}\si$ in
	$L^\infty\bigl(0,T;\bL^2(\Omega)\bigr)$, where all convergences hold for the same
	subsequence.
\end{theorem}


\begin{corollary}\label{cor:full sequence}
	Under the assumptions of Theorem~\ref{the:main convergence} and provided that the
	weak solution $(\mm,\si)$ of the SDLLG system~\eqref{equ:SDLLG} is unique, the
	limits of Theorem~\ref{the:main convergence} do not only hold for subsequences, but
	also for the full sequences $\mm_{h\tau}^\star$ and $\si_{h\tau}^\star$ as
	$(h,\tau)\to(0,0)$.
\end{corollary}


\section{Stability and convergence analysis}\label{sec:stability conv}

We aim to prove Theorem~\ref{the:main convergence}. To this end, we will derive the necessary stability estimates and use them to perform a compactness argument.


\subsection{Stability estimates}
\label{subsec:stability}

We begin with the following lemma, which shows that $\Pi_h(\widehat{\mm}_h^{n+1})$ is well-defined.

\begin{lemma}\label{lem:admissibility predictor}
	Let $\{(\mm_h^n,\vv_h^n)\}$ be generated by Algorithm~\ref{alg: full discr}. Then for $n=0,\ldots,N$,
	\begin{align}\label{equ:nodal m geq one}
		\abs{\mm_h^n(\vec{z})}\geq 1,
		\qquad
		\text{ for all}\quad \vec{z}\in\mathcal{N}_h^\omega.
	\end{align}
	Moreover, for every $n=1,\ldots,N-1$,
	\begin{align}
		\abs{\widehat{\mm}_h^{n+1}(\vec{z})}\geq 1,
		\qquad
		\text{ for all}\quad  \vec{z}\in\mathcal{N}_h^\omega.
		\label{equ:nodal mhat geq one}
	\end{align}
	In particular, for $n=1,\ldots,N-1$, we have $\widehat{\mm}_h^{n+1}\in\mathcal{U}_h$, thus $\Pi_h(\widehat{\mm}_h^{n+1})$ is well-defined.
\end{lemma}

\begin{proof}
	Fix $\vec{z}\in\mathcal{N}_h^\omega$. Since
	$\mm_h^0\in\mathcal{M}_h$ and
	$\vv_h^0\in\mathcal{K}_h(\mm_h^0)$, the initial update gives
	\begin{align}
		\abs{\mm_h^1(\vec{z})}^2
		&=
		\abs{
			\mm_h^0(\vec{z})
			+\tau\vv_h^0(\vec{z})
		}^2=
		\abs{\mm_h^0(\vec{z})}^2
		+
		\tau^2\abs{\vv_h^0(\vec{z})}^2=
		1+\tau^2\abs{\vv_h^0(\vec{z})}^2
		\geq1.
		\label{equ:first nodal constraint}
	\end{align}
	
	For $n\geq1$, the BDF2 update and the nodal tangency $\widehat{\mm}_h^{n+1}(\vec{z}) \cdot \vv_h^n(\vec{z})=0$ yield the nodewise identity
	\begin{align}\label{equ:nodal BDF2 recurrence}
		\frac32
		\abs{\mm_h^{n+1}(\vec{z})}^2
		-
		2\abs{\mm_h^n(\vec{z})}^2
		+
		\frac12
		\abs{\mm_h^{n-1}(\vec{z})}^2
		=
		\frac32\tau^4
		\abs{\dtau^2\mm_h^{n+1}(\vec{z})}^2.
	\end{align}
	This identity is derived in~\cite[Proposition~3.8]{AldFeiPra26}, where its proof uses only the BDF2 identities and the nodal orthogonality condition, and therefore applies unchanged to the present scheme.
	For brevity, set $r_n:= \abs{\mm_h^n(\vec{z})}^2- \abs{\mm_h^{n-1}(\vec{z})}^2.$
	Multiplying~\eqref{equ:nodal BDF2 recurrence} by $2/3$ and
	rearranging gives
	\begin{align}
		r_{n+1}
		=
		\frac13 r_n
		+
		\tau^4
		\abs{\dtau^2\mm_h^{n+1}(\vec{z})}^2.
		\label{equ:increment recurrence}
	\end{align}
	By~\eqref{equ:first nodal constraint}, $r_1= \tau^2\abs{\vv_h^0(\vec{z})}^2
	\geq 0$.
	Hence, an induction argument based on~\eqref{equ:increment recurrence} shows that
	$r_n\geq0$ for every $n\geq1$, and therefore that the sequence $\{|\mm_h^n(\vec{z})|^2\}_{n=0}^N$ is non-decreasing. Consequently, $\abs{\mm_h^n(\vec{z})}^2\geq |\mm_h^0(\zz)|^2 = 1$, which proves~\eqref{equ:nodal m geq one}.
	
	It remains to consider the extrapolated magnetization. Using
	$\widehat{\mm}_h^{n+1}=2\mm_h^n-\mm_h^{n-1}$, we compute
	\begin{align}\label{equ:predictor nodal identity}
		\abs{\widehat{\mm}_h^{n+1}(\vec{z})}^2
		&=
		\abs{2\mm_h^n(\vec{z}) -\mm_h^{n-1}(\vec{z})}^2
		=
		\abs{\mm_h^n(\vec{z})}^2
		+r_n
		+
		2\tau^2
		\abs{\dtau\mm_h^n(\vec{z})}^2.
	\end{align}
	Since $\abs{\mm_h^n(\vec{z})}^2 \geq1$ and $r_n\geq0$, it follows that
	\[
	\abs{\widehat{\mm}_h^{n+1}(\vec{z})}^2
	\geq 1.
	\]
	This proves~\eqref{equ:nodal mhat geq one}, and hence
	$\widehat{\mm}_h^{n+1}\in\mathcal{U}_h$, as required.
\end{proof}


We now derive various stability estimates for the discrete solutions generated by Algorithm~\ref{alg: full discr}. To this end, we recall the following identity, and we refer to \cite{Dah78,NevOde81}, or \cite[Equation (31)]{AldFeiPra26} for more details on it. Let $X$ be a Hilbert space. For any sequence $\{\ww^n\}_{n=0}^N\subset X$ and $n=1,\ldots,N-1$, we have
\begin{align}\label{equ:BDF2 energy identity}
	\tau\scalarproduct{\Dtau\ww^{n+1}}{\ww^{n+1}}_{X}
	&=
	\frac14
	\left(
	\norm[\big]{\ww^{n+1}}_{X}^2
	-
	\norm[\big]{\ww^n}_{X}^2
	+
	\norm[\big]{\widehat{\ww}^{n+2}}_{X}^2
	-
	\norm[\big]{\widehat{\ww}^{n+1}}_{X}^2 \right) 
	+
	\frac{\tau^4}{4}  \norm[\big]{\dtau^2 \ww^{n+1}}_{X}^2.
\end{align}


\begin{lemma}\label{lem:stab shn}
	Let $1\leq n\leq N$. The discrete solutions generated by
	Algorithm~\ref{alg: full discr} satisfy
	\begin{align}\label{equ:stab shn}
		\norm[\big]{\si_h^{n}}_{\bL^2(\Omega)}^2
		+ \norm[\big]{\widehat{\si}_h^{n+1}}_{\bL^2(\Omega)}^2
		+ \tau \sum_{j=0}^{n-1} \norm[\big]{\si_h^{j+1}}_{\bH^1(\Omega)}^2
		+ \tau^4 \sum_{j=1}^{n-1} \norm[\big]{\dtau^2 \si_h^{j+1}}_{\bL^2(\Omega)}^2
		\lesssim 
		1,
	\end{align}
	where the constant depends only on $C_{\si,0}$, $\const{C}{$j$}$, $\beta$, $\beta'$,
	$D_*$, $\Gamma_J$, and $\Omega$, but is independent of $h$, $\tau$, $N$.
\end{lemma}

\begin{proof}
	We first consider the initial step. Choosing
	$\zetaa_h=\tau\si_h^1$ in \eqref{equ:first step full discr s} gives
	\begin{align*}
		&
		\scalarproduct{\si_h^1-\si_h^0}{\si_h^1}_{\Omega}
		+
		\tau \fra_h \left(\Pi_h(\mm_h^0); \si_h^1, \si_h^1 \right)
		\\
		&=
		\tau\beta
		\scalarproduct{\Pi_h(\mm_h^0)\otimes\jj_h^1}{\nabla\si_h^1}_{\omega}
		- \tau\beta \scalarproduct{\jj_h^1\cdot\boldsymbol{\nu}}{\Pi_h(\mm_h^0)\cdot\si_h^1}_{\Gamma_J}.
	\end{align*}
	Using \eqref{equ:ah coercive}, together with the identity $(\aaa-\bb,\aaa) = \frac{1}{2}\|\aaa\|^2 - \frac{1}{2}\|\bb\|^2 + \frac{1}{2}\|\aaa-\bb\|^2$, we obtain
	\begin{align*}
		&
		\frac12
		\norm[\big]{\si_h^1}_{\bL^2(\Omega)}^2
		-
		\frac12
		\norm[\big]{\si_h^0}_{\bL^2(\Omega)}^2
		+
		\frac12
		\norm[\big]{\si_h^1-\si_h^0}_{\bL^2(\Omega)}^2
		+
		D_*(1-\beta\beta')\tau
		\norm[\big]{\si_h^1}_{\bH^1(\Omega)}^2
		\\
		&\leq
		\tau\beta
		\scalarproduct{\Pi_h(\mm_h^0)\otimes\jj_h^1}{\nabla\si_h^1}_{\omega}
		- \tau\beta \scalarproduct{\jj_h^1\cdot\boldsymbol{\nu}}{\Pi_h(\mm_h^0)\cdot\si_h^1}_{\Gamma_J}.
	\end{align*}
	Notice that the trace inequality \eqref{equ:trace} gives
	\begin{align*}
		\scalarproduct{\jj_h^1\cdot\boldsymbol{\nu}}{\Pi_h(\mm_h^0)\cdot\si_h^1}_{\Gamma_J}
		&\leq
		\norm[\big]{\jj_h^1}_{\bL^2(\Gamma_J)}
		\norm[\big]{\Pi_h(\mm_h^0)\cdot\si_h^1}_{\bL^2(\Gamma_J)}
		\\
		&\eqreff*{equ:Pi Linfty}{\lesssim}
		\norm[\big]{\jj_h^1}_{\bH^1(\Omega)}
		\norm[\big]{\si_h^1}_{\bH^1(\Omega)},
	\end{align*}
	and therefore, using \eqref{equ:Pi Linfty} also on the volume scalar product on the right-hand side, we have
	\begin{align*}
		\frac12
		\norm[\big]{\si_h^1}_{\bL^2(\Omega)}^2
		-
		\frac12
		\norm[\big]{\si_h^0}_{\bL^2(\Omega)}^2
		+
		\frac12
		\norm[\big]{\si_h^1-\si_h^0}_{\bL^2(\Omega)}^2
		&+
		D_*(1-\beta\beta')\tau
		\norm[\big]{\si_h^1}_{\bH^1(\Omega)}^2
		\\
		&\lesssim
		\tau
		\norm[\big]{\jj_h^1}_{\bH^1(\Omega)}
		\norm[\big]{\si_h^1}_{\bH^1(\Omega)}.
	\end{align*}
	By applying Young's inequality with parameter $\varepsilon = D_*(1-\beta\beta')$ and absorbing the $\bH^1$-norm of $\si_h^1$ to the left-hand side, we obtain
	\begin{align}
		\label{equ:initial spin stability}
		&
		\norm[\big]{\si_h^1}_{\bL^2(\Omega)}^2
		+
		\norm[\big]{\si_h^1-\si_h^0}_{\bL^2(\Omega)}^2
		+
		D_*(1-\beta\beta')\tau
		\norm[\big]{\si_h^1}_{\bH^1(\Omega)}^2
		\lesssim
		\norm[\big]{\si_h^0}_{\bL^2(\Omega)}^2
		+
		\tau
		\norm[\big]{\jj_h^1}_{\bH^1(\Omega)}^2.
	\end{align}
	Next, for $j=1,\ldots,n-1$, we take $\zetaa_h=\tau\si_h^{j+1}$ in \eqref{equ:full discr alg s} to obtain
	\begin{align*}
		&
		\tau
		\scalarproduct{\Dtau\si_h^{j+1}}{\si_h^{j+1}}_{\Omega}
		+
		\tau
		\fra_h
		\big(
		\Pi_h\bigl(\widehat{\mm}_h^{j+1}\bigr); \si_h^{j+1},\si_h^{j+1}\big)
		\\
		&=
		\tau\beta
		\scalarproduct{\Pi_h\bigl(\widehat{\mm}_h^{j+1}\bigr)\otimes\jj_h^{j+1}}{\nabla\si_h^{j+1}}_{\omega}
		-
		\tau\beta
		\scalarproduct{\jj_h^{j+1}\cdot\boldsymbol{\nu}}{\Pi_h\bigl(\widehat{\mm}_h^{j+1}\bigr)\cdot\si_h^{j+1}}_{\Gamma_J}.
	\end{align*}
	Therefore, using \eqref{equ:Pi Linfty}, \eqref{equ:ah coercive}, \eqref{equ:BDF2 energy identity}, and the trace inequality \eqref{equ:trace} as above, we infer that
	\begin{align}
		\label{equ:one step spin stability}
		&
		\frac14
		\left(
		\norm[\big]{\si_h^{j+1}}_{\bL^2(\Omega)}^2 - \norm[\big]{\si_h^j}_{\bL^2(\Omega)}^2
		+
		\norm[\big]{\widehat{\si}_h^{j+2}}_{\bL^2(\Omega)}^2 - \norm[\big]{\widehat{\si}_h^{j+1}}_{\bL^2(\Omega)}^2
		\right)
		+
		\frac{\tau^4}{4}
		\norm[\big]{\dtau^2 \si_h^{j+1}}_{\bL^2(\Omega)}^2
		\nonumber
		\\
		&\quad
		+
		D_*(1-\beta\beta')\tau
		\norm[\big]{\si_h^{j+1}}_{\bH^1(\Omega)}^2
		\nonumber
		\\
		&\lesssim
		\tau \norm[\big]{\jj_h^{j+1}}_{\bH^1(\Omega)}
		\norm[\big]{\si_h^{j+1}}_{\bH^1(\Omega)}
		\nonumber \\
		&\leq
		\frac12 D_*(1-\beta\beta')\tau \norm[\big]{\si_h^{j+1}}_{\bH^1(\Omega)}^2
		+
		C_1\tau \norm[\big]{\jj_h^{j+1}}_{\bH^1(\Omega)}^2,
	\end{align}
	where in the last step we used Young's inequality.
	Rearranging the terms, summing~\eqref{equ:one step spin stability} over $j=1,\ldots,n-1$, and combining the result with~\eqref{equ:initial spin stability} gives
	\begin{align*}
		\frac14\big(\norm[\big]{\si_h^n}_{\bL^2(\Omega)}^2 + \norm[\big]{\widehat{\si}_h^{n+1}}_{\bL^2(\Omega)}^2\big) + \frac{\tau^4}{4}\sum_{j=1}^{n-1} \norm[\big]{\dtau^2 \si_h^{j+1}}_{\bL^2(\Omega)}^2 &+ \frac12 D_*(1-\beta\beta')\tau \sum_{j=0}^{n-1} \norm[\big]{\si_h^{j+1}}_{\bH^1(\Omega)}^2 \\
		&
		\lesssim
		\norm[\big]{\si_h^0}_{\bL^2(\Omega)}^2 + \tau \sum_{j=0}^{n-1} \norm[\big]{\jj_h^{j+1}}_{\bH^1(\Omega)}^2.
	\end{align*}
	Finally, \eqref{equ:Cf Cj} and assumption~\eqref{ass:initial data} allow us to prove \eqref{equ:stab shn} and conclude the proof.
\end{proof}


\begin{lemma}
	The discrete solutions generated by Algorithm~\ref{alg: full discr} satisfy
	\begin{align}\label{equ:spin dt Hminus1}
		\tau
		\sum_{n=0}^{N-1} \norm[\big]{{R_h^*(\dtau\si_h^{n+1})}}_{\widetilde{\bH}^{-1}(\Omega)}^2
		\lesssim 1,
	\end{align}
	where the implicit constant depends only on $C_{\si,0}$, $\const{C}{$j$}$, $\beta$, $\beta'$, $D_*$, $D^\ast$, $\Omega$, $\Gamma_J$, and on $\gamma$, but is independent of $h$, $\tau$, and $N$.
\end{lemma}

\begin{proof}
	We first consider the initial step. For any $\zetaa \in \bH^1(\Omega)$, by using \eqref{equ:dual reconstruction} and \eqref{equ:first step full discr s}, together with $\norm[\big]{\Pi_h(\mm_h^0)}_{\bL^\infty(\omega)}\leq1$ and the continuity of $\fra_h$, we have
	\[
	|\fra_h(\Pi_h(\mm_h^0);\si_h^1,{R_h\zetaa})| \lesssim \norm[\big]{\si_h^1}_{\bH^1(\Omega)} \norm[\big]{{R_h\zetaa}}_{\bH^1(\Omega)}.
	\]
	{The stability of $R_h$ in~\eqref{equ:H1 projection}} and the trace inequality~\eqref{equ:trace} imply that
	\begin{align*}
		\dual{{R_h^*(\dtau\si_h^1)}}{\zetaa}_{\widetilde{\bH}^{-1}(\Omega),\bH^1(\Omega)}
		&\eqreff*{equ:dual reconstruction}{=}
		\scalarproduct{\dtau\si_h^1}{{R_h\zetaa}}_{\Omega}
		\\
		&=
		\beta
		\scalarproduct{\Pi_h(\mm_h^0)\otimes\jj_h^1}{\nabla {R_h\zetaa}}_{\omega}
		-
		\beta
		\scalarproduct{\jj_h^1\cdot\boldsymbol{\nu}}{\Pi_h(\mm_h^0)\cdot {R_h\zetaa}}_{\Gamma_J}
		\\
		&\quad
		-
		\fra_h
		\left(
		\Pi_h(\mm_h^0);
		\si_h^1,
		{R_h\zetaa}
		\right)
		\\
		&\lesssim
		\left(
		\norm[\big]{\si_h^1}_{\bH^1(\Omega)}
		+
		\norm[\big]{\jj_h^1}_{\bH^1(\Omega)}
		\right)
		\norm[\big]{\zetaa}_{\bH^1(\Omega)}.
	\end{align*}
	Therefore,
	\begin{align}
		\label{equ:initial spin Hminus}
		\norm[\big]{{R_h^*(\dtau\si_h^1)}}_{\widetilde{\bH}^{-1}(\Omega)}
		\lesssim
		\norm[\big]{\si_h^1}_{\bH^1(\Omega)}
		+
		\norm[\big]{\jj_h^1}_{\bH^1(\Omega)}.
	\end{align}
	For $n=1,\ldots,N-1$, now taking $\zetaa_h = {R_h\zetaa}$ in \eqref{equ:full discr alg s}, the same argument yields
	\begin{align*}
		\dual{{R_h^*(\Dtau\si_h^{n+1})}}{\zetaa}_{\widetilde{\bH}^{-1}(\Omega),\bH^1(\Omega)}
		&\eqreff*{equ:dual reconstruction}{=}
		\scalarproduct{\Dtau\si_h^{n+1}}{{R_h\zetaa}}_{\Omega}
		\\
		&=
		\beta \scalarproduct{\Pi_h(\widehat{\mm}_h^{n+1})\otimes\jj_h^{n+1}}{\nabla {R_h\zetaa}}_{\omega}
		\\
		&\quad
		-
		\beta
		\scalarproduct{\jj_h^{n+1}\cdot\boldsymbol{\nu}}{\Pi_h(\widehat{\mm}_h^{n+1})\cdot {R_h\zetaa}}_{\Gamma_J}
		\\
		&\quad
		-
		\fra_h
		\left(
		\Pi_h(\widehat{\mm}_h^{n+1});
		\si_h^{n+1},
		{R_h\zetaa}
		\right)
		\\
		&\lesssim
		\left(
		\norm[\big]{\si_h^{n+1}}_{\bH^1(\Omega)}
		+
		\norm[\big]{\jj_h^{n+1}}_{\bH^1(\Omega)}
		\right)
		\norm[\big]{\zetaa}_{\bH^1(\Omega)}.
	\end{align*}
	This implies
	\begin{align}
		\label{equ:BDF2 spin Hminus}
		\norm[\big]{{R_h^*(\Dtau\si_h^{n+1})}}_{\widetilde{\bH}^{-1}(\Omega)}
		\lesssim
		\norm[\big]{\si_h^{n+1}}_{\bH^1(\Omega)}
		+
		\norm[\big]{\jj_h^{n+1}}_{\bH^1(\Omega)}.
	\end{align}	
	Squaring \eqref{equ:initial spin Hminus} and~\eqref{equ:BDF2 spin Hminus}, multiplying by $\tau$ and summing, we obtain by \eqref{equ:Cf Cj} and \eqref{equ:stab shn} 
	\begin{align}\label{equ:BDF2 spin Hminus sum}
		\tau \norm[\big]{{R_h^*(\dtau\si_h^1)}}_{\widetilde{\bH}^{-1}(\Omega)}^2
		+
		\tau \sum_{n=1}^{N-1} \norm[\big]{{R_h^*(\Dtau\si_h^{n+1})}}_{\widetilde{\bH}^{-1}(\Omega)}^2
		\lesssim 1.
	\end{align}
	Hence, noting that $R_h^*$ is linear, applying Lemma~\ref{lem:BDF2 norm equivalence} {to $\{R_h^*\si_h^n\}_{n=0}^N$ in $\widetilde{\bH}^{-1}(\Omega)$} and using~\eqref{equ:BDF2 spin Hminus sum}, we deduce~\eqref{equ:spin dt Hminus1}, as required.
\end{proof}


We now obtain some stability estimates for the discrete magnetization.


\begin{lemma}\label{lem:stability magnetization}
	The discrete solutions generated by Algorithm~\ref{alg: full discr} satisfy
	\begin{align}\label{equ:magnetization stability}
		&
		\max_{0\leq n\leq N}
		\norm[\big]{\mm_h^n}_{\bH^1(\omega)}^2
		+
		\tau
		\sum_{n=0}^{N-1}
		\norm[\big]{\vv_h^n}_{\bL^2(\omega)}^2
		+
		\tau^2
		\norm[\big]{\nabla\vv_h^0}_{\bL^2(\omega)}^2
		\nonumber\\
		&\qquad
		+
		\tau^3
		\sum_{n=1}^{N-1}
		\norm[\big]{\dtau^2 \mm_h^{n+1}}_{\bL^2(\omega)}^2
		+
		\tau^4
		\sum_{n=1}^{N-1}
		\norm[\big]{\dtau^2\mm_h^{n+1}}_{\bH^1(\omega)}^2
		\lesssim 1.
	\end{align}
	Consequently,
	\begin{align}\label{equ:predictor m stability}
		\max_{0\leq n\leq N-1}
		\norm[\big]{\widehat{\mm}_h^{n+1}}_{\bH^1(\omega)}
		+
		\max_{0\leq n\leq N-1}
		\norm[\big]{
			\Pi_h\bigl(\widehat{\mm}_h^{n+1}\bigr)
		}_{\bH^1(\omega)}
		\lesssim 1,
	\end{align}
	where the implicit constant depends only on $\alpha$, $\const{C}{ex}$, $\const{C}{$\pi$}$, $C_{\mm,0}$, $C_{\si,0}$, $\const{C}{$f$}$, $\const{C}{$j$}$, $c$, $\beta$, $\beta'$, $D_*$, $T$, $\omega$, $\Omega$, $\Gamma_J$, and on $\gamma$, but is independent of $h$, $\tau$, and $N$.
\end{lemma}

\begin{proof}
	We first consider the initial step. Set
	\[
	\boldsymbol{F}_h^1
	:=
	\pp_h(\mm_h^0)+\ff_h^1+c\si_h^0.
	\]
	Using $\mm_h^1=\mm_h^0+\tau\vv_h^0$, we can rewrite \eqref{equ:first step full discr m} as
	\[
	\alpha\scalarproduct{\vv_h^0}{\phii_h}_{\omega}
	+
	\scalarproduct{\mm_h^0\times\vv_h^0}{\phii_h}_{\omega}
	+
	\const{C}{ex}
	\scalarproduct{\nabla\mm_h^1}{\nabla\phii_h}_{\omega}
	=
	\scalarproduct{\boldsymbol{F}_h^1}{\phii_h}_{\omega}.
	\]
	Testing with $\phii_h=\tau\vv_h^0$, using
	$\scalarproduct{\mm_h^0\times\vv_h^0}{\vv_h^0}_{\omega}=0$,
	and applying Young's inequality give
	\begin{align}\label{equ:initial m stability}
		\const{C}{ex}
		\norm[\big]{\nabla\mm_h^1}_{\bL^2(\omega)}^2
		+
		\alpha\tau
		\norm[\big]{\vv_h^0}_{\bL^2(\omega)}^2
		+
		\const{C}{ex}\tau^2
		\norm[\big]{\nabla\vv_h^0}_{\bL^2(\omega)}^2
		&\leq
		\const{C}{ex}
		\norm[\big]{\nabla\mm_h^0}_{\bL^2(\omega)}^2
		+
		\frac{\tau}{\alpha}
		\norm[\big]{\boldsymbol{F}_h^1}_{\bL^2(\omega)}^2
		\nonumber \\
		&\lesssim 1,
	\end{align}
	where the final bound follows from
	assumptions~\eqref{ass:initial data}--\eqref{ass:pi bounded}.
	
	Next, for $n\geq1$, set
	\begin{align}\label{equ:discrete F}
		\boldsymbol{F}_h^{n+1}
		:=
		\pp_h\bigl(\widehat{\mm}_h^{n+1}\bigr)
		+ \ff_h^{n+1}
		+ c\widehat{\si}_h^{n+1}.
	\end{align}
	By the convention $\widehat{\mm}_h^1=\mm_h^0$ and $\widehat{\si}_h^1=\si_h^0$, definition~\eqref{equ:discrete F} with $n=0$ coincides with $\boldsymbol{F}_h^1$ defined above. Equation~\eqref{equ:full discr alg m} can be equivalently written as
	\begin{align}\label{equ:vhn equiv}
		\alpha\scalarproduct{\vv_h^n}{\phii_h}_{\omega}
		+
		\scalarproduct{
			\widehat{\mm}_h^{n+1}\times\vv_h^n
		}{
			\phii_h
		}_{\omega}
		+
		\const{C}{ex}
		\scalarproduct{\nabla\mm_h^{n+1}}{\nabla\phii_h}_{\omega}
		=
		\scalarproduct{\boldsymbol{F}_h^{n+1}}{\phii_h}_{\omega}.
	\end{align}
	Now, for ease of presentation, define for $n \ge 1$
	\begin{align}\label{equ:G definition}
		\mathcal{G}_h^n
		\coloneqq
		&\frac14
		\norm[\big]{\nabla\mm_h^{n-1}}_{\bL^2(\omega)}^2
		- \scalarproduct{\nabla\mm_h^{n-1}}{\nabla\mm_h^{n}}_{\omega}
		+ \frac54
		\norm[\big]{\nabla\mm_h^n}_{\bL^2(\omega)}^2
		\\
		=
		&\frac14
		\big(
		\norm[\big]{\nabla\mm_h^n}_{\bL^2(\omega)}^2
		+
		\norm[\big]{
			2\nabla\mm_h^n-\nabla\mm_h^{n-1}
		}_{\bL^2(\omega)}^2
		\big).\nonumber
	\end{align}
	Applying the identity~\eqref{equ:BDF2 energy identity} with $\ww_h^n = \nabla \mm_h^n$, we obtain
	\[
	\tau \scalarproduct{\nabla\mm_h^{n+1}}{\nabla\vv_h^n}_{\omega}
	=
	\mathcal{G}_h^{n+1}-\mathcal{G}_h^n+\frac{\tau^4}{4}\norm[\big]{\nabla\dtau^2\mm_h^{n+1}}_{\bL^2(\omega)}^2.
	\]
	Therefore, setting~$\phii_h=\tau\vv_h^n=\tau\Dtau\mm_h^{n+1}$ in~\eqref{equ:vhn equiv} gives
	\[
	\const{C}{ex}
	\left(
	\mathcal{G}_h^{n+1}-\mathcal{G}_h^n
	\right)
	+
	{\alpha\tau}
	\norm[\big]{\vv_h^n}_{\bL^2(\omega)}^2
	+
	\frac{\const{C}{ex}\tau^4}{4}
	\norm[\big]{\nabla\dtau^2\mm_h^{n+1}}_{\bL^2(\omega)}^2
	=
	\tau\scalarproduct{\boldsymbol{F}_h^{n+1}}{\vv_h^n}_\omega,
	\]
	and applying Young's inequality to the right-hand side yields
	\begin{align}
		\label{equ:m one step G stability}
		&
		\const{C}{ex}
		\left(
		\mathcal{G}_h^{n+1}-\mathcal{G}_h^n
		\right)
		+
		\frac{\alpha\tau}{2}
		\norm[\big]{\vv_h^n}_{\bL^2(\omega)}^2
		+
		\frac{\const{C}{ex}\tau^4}{4}
		\norm[\big]{\nabla\dtau^2\mm_h^{n+1}}_{\bL^2(\omega)}^2
		\leq
		\frac{\tau}{2\alpha}
		\norm[\big]{\boldsymbol{F}_h^{n+1}}_{\bL^2(\omega)}^2.
	\end{align}
	It remains to control the right-hand side of \eqref{equ:m one step G stability}.
	To this end, we use the explicit formula~\eqref{equ:explicit nodal BDF2 constraint} for the nodal values $\abs{\mm_h^\ell(\zz)}^2$ from Lemma~\ref{lem:explicit nodal BDF2 constraint}. {Multiplying it by $\mu_{\zz}^{\omega}$, summing over all nodes $\zz\in\mathcal N_h^\omega$,} and using the norm equivalence~\eqref{equ:discrete nodal norm}, we obtain, for all $1\leq\ell\leq N$,
	\begin{align}
		\label{equ:m L2 preliminary}
		\norm[\big]{\mm_h^\ell}_{\bL^2(\omega)}^2
		\lesssim
		\Big(1+ \tau^2 \norm[\big]{\vv_h^0}_{\bL^2(\omega)}^2
		+
		\tau^4 \sum_{j=2}^{\ell} \norm[\big]{\dtau^2\mm_h^j}_{\bL^2(\omega)}^2 \Big).
	\end{align}
	To estimate the last term, we use the definition of $\dtau^2$, together with \eqref{equ:BDF2 id 1}, \eqref{equ:BDF2 id 2}, and Lemma~\ref{lem:BDF2 norm equivalence} to obtain
	\begin{align}
		\tau
		\sum_{j=2}^{\ell}
		\norm[\big]{\dtau^2\mm_h^j}_{\bL^2(\omega)}^2
		&=
		\frac{1}{\tau}
		\sum_{j=2}^{\ell}
		\norm[\big]{
			\dtau\mm_h^j-\dtau\mm_h^{j-1}
		}_{\bL^2(\omega)}^2
		\nonumber
		\leq
		\frac{4}{\tau}
		\sum_{j=1}^{\ell}
		\norm[\big]{\dtau\mm_h^j}_{\bL^2(\omega)}^2
		\nonumber \\
		&=
		\frac{4}{\tau^2}
		\Big(
		\tau
		\sum_{j=0}^{\ell-1}
		\norm[\big]{\dtau\mm_h^{j+1}}_{\bL^2(\omega)}^2
		\Big)
		\nonumber
		\eqreff*{equ:BDF2 norm equivalence}{\lesssim}
		\frac{1}{\tau^2}
		\Big(
		\tau
		\sum_{j=0}^{\ell-1}
		\norm[\big]{\vv_h^j}_{\bL^2(\omega)}^2
		\Big).\nonumber
	\end{align}
	In particular,
	\begin{align}\label{equ:dtt m bound multiplied}
		\tau^4
		\sum_{j=2}^{\ell}
		\norm[\big]{\dtau^2\mm_h^j}_{\bL^2(\omega)}^2
		\lesssim
		\tau
		\Big(
		\tau \sum_{j=0}^{\ell-1} \norm[\big]{\vv_h^j}_{\bL^2(\omega)}^2
		\Big).
	\end{align}
	It then follows from~\eqref{equ:m L2 preliminary} and \eqref{equ:dtt m bound multiplied} that, for any $1 \leq \ell\leq N$,
	\begin{align}\label{equ:m L2 velocity bound}
		\norm[\big]{\mm_h^\ell}_{\bL^2(\omega)}^2
		\lesssim
		\bigg[
		1 + \tau \Big(\tau \sum_{j=0}^{\ell-1} \norm[\big]{\vv_h^j}_{\bL^2(\omega)}^2 \Big)
		\bigg].
	\end{align}
	Moreover, $\mm^0_h\in\mathcal M_h$ {implies $|\mm_h^0|\leq1$ on each element by convexity, thus $\norm[\big]{\mm^0_h}^2_{\bL^2(\omega)}\leq|\omega|$},
	and \eqref{equ:m L2 velocity bound} holds for every $0\le\ell\le N$.
	Consequently, using
	$\widehat{\mm}_h^{n+1}=2\mm_h^n-\mm_h^{n-1}$, we also have for $n=1,\ldots,N-1$,
	\begin{align}\label{equ:mhat L2 velocity bound}
		\norm[\big]{\widehat{\mm}_h^{n+1}}_{\bL^2(\omega)}^2
		\lesssim
		\bigg[
		1 + \tau \Big(\tau \sum_{j=0}^{n-1} \norm[\big]{\vv_h^j}_{\bL^2(\omega)}^2\Big)
		\bigg].
	\end{align}
	By assumption~\eqref{ass:pi bounded}, \eqref{equ:Cf Cj}, Lemma~\ref{lem:stab shn}, and
	\eqref{equ:mhat L2 velocity bound}, we infer that
	\begin{align}
		\label{equ:forcing m bound revised}
		&
		\tau
		\sum_{n=1}^{\ell-1}
		\norm[\big]{\boldsymbol{F}_h^{n+1}}_{\bL^2(\omega)}^2
		\lesssim
		1
		+
		\tau^2
		\sum_{n=1}^{\ell-1}
		\Big(
		\tau
		\sum_{j=0}^{n-1}
		\norm[\big]{\vv_h^j}_{\bL^2(\omega)}^2
		\Big).
	\end{align}
	Now, we continue from~\eqref{equ:m one step G stability}. Summing this expression over $n=1,\ldots,\ell-1$ yields
	\begin{align}\label{equ:summed G stability}
		\const{C}{ex}\mathcal{G}_h^\ell
		+
		\frac{\alpha\tau}{2}
		\sum_{n=1}^{\ell-1}
		\norm[\big]{\vv_h^n}_{\bL^2(\omega)}^2
		+
		\frac{\const{C}{ex}\tau^4}{4}
		\sum_{n=1}^{\ell-1}
		\norm[\big]{\nabla\dtau^2\mm_h^{n+1}}_{\bL^2(\omega)}^2
		\leq
		\const{C}{ex}\mathcal{G}_h^1
		+
		\frac{\tau}{2\alpha}
		\sum_{n=1}^{\ell-1}
		\norm[\big]{\boldsymbol{F}_h^{n+1}}_{\bL^2(\omega)}^2.
	\end{align}
	To conclude, we need some estimates for $\mathcal{G}_h^n$. To this end, let
	\[
	X:=\bL^2(\omega;\R^{3\times 3}),
	\qquad
	\boldsymbol{x}
	:=
	\begin{pmatrix}
		\nabla\mm_h^{n-1}\\
		\nabla\mm_h^n
	\end{pmatrix}
	\in X\times X.
	\]
	Observe that
	\[
	G
	:=
	\frac14
	\begin{pmatrix}
		1 & -2\\
		-2 & 5
	\end{pmatrix}
	\]
	is symmetric positive definite with eigenvalues $\mu^\pm = \frac14 \bigl( 3\pm2\sqrt{2}\bigr)$.
	Hence, by the spectral theorem,
	\[
	\mu^-
	\norm[\big]{\boldsymbol{x}}_{X\times X}^2
	\leq
	\scalarproduct{G\boldsymbol{x}}{\boldsymbol{x}}_{X\times X}
	\leq
	\mu^+
	\norm[\big]{\boldsymbol{x}}_{X\times X}^2.
	\]
	We then obtain
	\begin{align}\label{equ:G norm equivalence}
		\mu^-
		\Big(\norm[\big]{\nabla\mm_h^{n-1}}_{\bL^2(\omega)}^2
		+
		\norm[\big]{\nabla\mm_h^n}_{\bL^2(\omega)}^2
		\Big)
		\leq
		\mathcal{G}_h^n
		\leq
		\mu^+
		\Big(\norm[\big]{\nabla\mm_h^{n-1}}_{\bL^2(\omega)}^2
		+
		\norm[\big]{\nabla\mm_h^n}_{\bL^2(\omega)}^2
		\Big).
	\end{align}
	Moreover, the initial-step estimate~\eqref{equ:initial m stability} implies $\mathcal{G}_h^1\leq C$. Thus, from~\eqref{equ:forcing m bound revised} and~\eqref{equ:summed G stability}, after dropping the non-negative terms from the left-hand side and adding the initial contribution $\tau\norm[\big]{\vv_h^0}_{\bL^2(\omega)}^2\leq C$, we infer that
	\begin{equation*}
		\tau\sum_{j=0}^{\ell-1} \norm[\big]{\vv_h^j}_{\bL^2(\omega)}^2
		\lesssim
		1
		+
		\tau^2
		\sum_{n=1}^{\ell-1}
		\bigg(
		\tau \sum_{j=0}^{n-1} \norm[\big]{\vv_h^j}_{\bL^2(\omega)}^2
		\bigg).
	\end{equation*}
	An application of the discrete Gronwall inequality, see Lemma~\ref{lem: gronwall} in Appendix~\ref{sec:appendix}, gives
	\begin{align} \label{equ:velocity stability}
		\tau
		\sum_{j=0}^{\ell-1}
		\norm[\big]{\vv_h^j}_{\bL^2(\omega)}^2
		\lesssim \exp\big(\sum_{n=0}^{\ell-1}C_{\mathrm{Gr}}\tau^2\big)= \exp(C_{\mathrm{Gr}}\ell \tau^2)\leq \exp(C_{\mathrm{Gr}}T^2)\lesssim 1,
	\end{align}
	which implies 
	\begin{align*}
		\max_{1\leq\ell\leq N}
		\tau
		\sum_{j=0}^{\ell-1}
		\norm[\big]{\vv_h^j}_{\bL^2(\omega)}^2
		\lesssim 1.
	\end{align*}
	It now follows from~\eqref{equ:dtt m bound multiplied}, \eqref{equ:forcing m bound revised}, and~\eqref{equ:velocity stability} that
	\begin{align}\label{equ:bound dtt and F}
		\tau^3
		\sum_{n=1}^{N-1}
		\norm[\big]{\dtau^2\mm_h^{n+1}}_{\bL^2(\omega)}^2
		+
		\tau\sum_{n=1}^{N-1} \norm[\big]{\boldsymbol{F}_h^{n+1}}_{\bL^2(\omega)}^2
		\lesssim 1.
	\end{align}
	Returning to~\eqref{equ:summed G stability} and using the lower bound from~\eqref{equ:G norm equivalence}, we obtain, for every
	$\ell=2,\ldots,N$,
	\begin{align*}
		&
		\mu^- \const{C}{ex}
		\left(
		\norm[\big]{\nabla\mm_h^{\ell-1}}_{\bL^2(\omega)}^2
		+
		\norm[\big]{\nabla\mm_h^\ell}_{\bL^2(\omega)}^2
		\right)
		+
		\frac{\alpha\tau}{2}
		\sum_{n=1}^{\ell-1}
		\norm[\big]{\vv_h^n}_{\bL^2(\omega)}^2
		\\
		&\qquad
		+
		\tau^3
		\sum_{n=1}^{\ell-1}
		\norm[\big]{\dtau^2\mm_h^{n+1}}_{\bL^2(\omega)}^2
		+
		\frac{\const{C}{ex}\tau^4}{4}
		\sum_{n=1}^{\ell-1}
		\norm[\big]{\nabla\dtau^2\mm_h^{n+1}}_{\bL^2(\omega)}^2
		\lesssim 1.
	\end{align*}
	Together with the initial-step estimate \eqref{equ:initial m stability}, this implies
	\begin{align*}
		&\max_{0\leq n\leq N}
		\norm[\big]{\nabla\mm_h^n}_{\bL^2(\omega)}^2
		+
		\tau
		\sum_{n=0}^{N-1}
		\norm[\big]{\vv_h^n}_{\bL^2(\omega)}^2
		+
		\tau^2
		\norm[\big]{\nabla\vv_h^0}_{\bL^2(\omega)}^2
		\\
		&\qquad
		+
		\tau^3
		\sum_{n=1}^{N-1}
		\norm[\big]{\dtau^2\mm_h^{n+1}}_{\bL^2(\omega)}^2
		+
		\tau^4
		\sum_{n=1}^{N-1}
		\norm[\big]{\dtau^2\mm_h^{n+1}}_{\bH^1(\omega)}^2
		\lesssim 1.
	\end{align*}
	Moreover, combining~\eqref{equ:m L2 velocity bound} with~\eqref{equ:velocity stability} gives
	\[
	\max_{0\leq n\leq N}
	\norm[\big]{\mm_h^n}_{\bL^2(\omega)}^2
	\lesssim 1,
	\]
	from which we deduce the required inequality~\eqref{equ:magnetization stability}.
	
	Finally, since $\widehat{\mm}_h^{n+1}=2\mm_h^n-\mm_h^{n-1}$, from~\eqref{equ:magnetization stability}, \eqref{equ:Pi Linfty}, \eqref{equ:Pi H1 stability}, and Lemma~\ref{lem:admissibility predictor}, we infer the estimate \eqref{equ:predictor m stability} for $n\ge 1$.
	For $n=0$, the convention $\widehat{\mm}^{\,1}_h:=\mm^0_h\in\mathcal M_h$ gives
	$\Pi_h(\widehat{\mm}^{\,1}_h)=\mm^0_h$, so that both terms are bounded by $C_{\mm,0}$.
	This concludes the proof.
\end{proof}


The next lemma shows some estimates of the unit-length and projection defects.

\begin{lemma}[Unit-length and projection defects]\label{lem:unit length defect}
	The discrete solutions generated by Algorithm~\ref{alg: full discr} satisfy for $n=0,\ldots,N$,
	\begin{align}
		\label{equ:nodal interpolated modulus defect}
		\norm[\big]{\III_h\bigl(\abs{\mm_h^n}^2\bigr)-1}_{L^1(\omega)}
		&\lesssim
		\tau,
		\\
		\label{equ:modulus defect}
		\norm[\big]{\abs{\mm_h^n}^2 -1}_{L^1(\omega)}
		&\lesssim
		h^2+\tau.
	\end{align}
	Moreover, for $n=0,\ldots,N-1$,
	\begin{align}\label{equ:predictor projection defect}
		\norm[\big]{\widehat{\mm}_h^{n+1}- \Pi_h\bigl(\widehat{\mm}_h^{n+1}\bigr)}_{\bL^1(\omega)}
		\lesssim \tau.
	\end{align}
	Here, the implicit constant depends only on $\alpha$, $\const{C}{ex}$, $\const{C}{$\pi$}$, $C_{\mm,0}$, $C_{\si,0}$, $\const{C}{$f$}$, $\const{C}{$j$}$, $c$, $\beta$, $\beta'$, $D_*$, $T$, $\omega$, $\Omega$, $\Gamma_J$, and on $\gamma$, but is independent of $h$, $\tau$, and $N$.
\end{lemma}

\begin{proof}
	First, notice that if $n=0$, $\III_h(\abs{\mm_h^0}^2)-1 = 0,$ and the estimate~\eqref{equ:nodal interpolated modulus defect} follows immediately. For $n\ge 1$ we note that the right-hand side of \eqref{equ:explicit nodal BDF2 constraint} is nonnegative. We then integrate the nonnegative nodal interpolant using~\eqref{equ:nodal weights} and apply the discrete norm equivalence~\eqref{equ:discrete nodal norm} and estimate~\eqref{equ:magnetization stability} to obtain
	\begin{align*}
		\norm[\big]{\III_h\bigl(\abs{\mm_h^n}^2\bigr)-1}_{L^1(\omega)}
		&\lesssim
		\tau^2
		\norm[\big]{\vv_h^0}_{\bL^2(\omega)}^2
		+
		\tau^4
		\sum_{j=2}^{n}
		\norm[\big]{\dtau^2\mm_h^j}_{\bL^2(\omega)}^2
		\lesssim \tau.
	\end{align*}
	This proves~\eqref{equ:nodal interpolated modulus defect}.
	To show~\eqref{equ:modulus defect}, we first note that by the triangle inequality,
	\begin{align*}
		\norm[\big]{\abs{\mm_h^n}^2-1}_{L^1(\omega)}
		&\leq
		\norm[\big]{\III_h\bigl(\abs{\mm_h^n}^2\bigr)-1}_{L^1(\omega)}
		+
		\norm[\big]{\abs{\mm_h^n}^2 -\III_h\bigl(\abs{\mm_h^n}^2\bigr)}_{L^1(\omega)}.
	\end{align*}
	The first term is already $\mathcal{O}(\tau)$ by \eqref{equ:nodal interpolated modulus defect}. For the second term, since $\abs{\mm_h^n}^2$ is a polynomial of degree at most two
	on $K$, standard scaling arguments give
	\[
	\norm[\big]{\abs{\mm_h^n}^2- \III_h\bigl(\abs{\mm_h^n}^2\bigr)}_{L^1(K)}
	\lesssim
	h_K^2
	\big|\abs{\mm_h^n}^2\big|_{{W}^{2,1}(K)}.
	\]
	Since $\mm_h^n$ is affine on $K$, we have
	$\big|\abs{\mm_h^n}^2\big|_{{W}^{2,1}(K)}
	\lesssim \norm[\big]{\nabla\mm_h^n}_{\bL^2(K)}^2$.
	Therefore, summing over $K\in\mathcal T_h^\omega$ and using the stability estimate \eqref{equ:magnetization stability}, we obtain
	\[
	\norm[\big]{\abs{\mm_h^n}^2- \III_h\bigl(\abs{\mm_h^n}^2\bigr)}_{L^1(\omega)}
	\lesssim h^2,
	\]
	from which \eqref{equ:modulus defect} follows.
	Finally, notice that for $n=0$, the estimate~\eqref{equ:predictor projection defect} holds trivially and is exactly $0$, since $\widehat{\mm}_h^1 = \mm_h^0\in \mathcal{M}_h$. For $n\ge 1$, since
	$\widehat{\mm}_h^{n+1}\in\mathcal{U}_h$, we have for every
	$\vec{z}\in\mathcal{N}_h^\omega$,
	\begin{align*}
		\abs{\widehat{\mm}_h^{n+1}(\vec{z}) - \Pi_h\bigl(\widehat{\mm}_h^{n+1}\bigr)(\vec{z})}
		&=
		\abs{\widehat{\mm}_h^{n+1}(\vec{z})}-1
		\\
		&\leq
		\frac12 \left(\abs{\widehat{\mm}_h^{n+1}(\vec{z})}^2 -1 \right)
		\\
		&\eqreff*{equ:predictor nodal identity}{=}
		\frac12 \left( 2\abs{\mm_h^n(\vec{z})}^2- \abs{\mm_h^{n-1}(\vec{z})}^2 + 2\tau^2 \abs{\dtau\mm_h^n(\vec{z})}^2 -1 \right)
		\\
		&\eqreff*{equ:nodal m geq one}{\leq}
		\bigl(\abs{\mm_h^n(\vec{z})}^2 - 1\bigr) + \tau^2 \abs{\dtau\mm_h^n(\vec{z})}^2.
	\end{align*}
	Hence, using~\eqref{equ:discrete nodal norm} multiple times, it follows that
	\begin{align*}
		\norm[\big]{\widehat{\mm}_h^{n+1} -\Pi_h\bigl(\widehat{\mm}_h^{n+1}\bigr)}_{\bL^1(\omega)}
		&\eqreff*{equ:discrete nodal norm}{\lesssim}
		\sum_{\vec{z}\in\mathcal{N}_h^\omega}{\mu_{\vec z}^\omega}
		\bigl( \abs{\widehat{\mm}_h^{n+1}(\vec{z})}-1 \bigr)
		\\
		&\lesssim \sum_{\vec{z}\in\mathcal{N}_h^\omega}{\mu_{\vec z}^\omega} \bigl( \abs{\mm_h^n(\vec{z})}^2 - 1 + \tau^2 \abs{\dtau\mm_h^n(\vec{z})}^2 \bigr)
		\\
		& \eqreff*{equ:explicit nodal BDF2 constraint}{\lesssim} \sum_{\vec{z}\in\mathcal{N}_h^\omega}{\mu_{\vec z}^\omega}
		\bigl( \tau^2 |\vv_h^0(\vec{z})|^2 + \tau^4 \sum_{j=2}^{n} |\dtau^2\mm_h^j(\vec{z})|^2 + \tau^2 |\dtau\mm_h^n(\vec{z})|^2 \bigr)
		\\
		&\eqrefftwo*{equ:discrete nodal norm}{equ:magnetization stability}{\lesssim}
		\tau
		+
		\tau^2
		\norm[\big]{\dtau\mm_h^n}_{\bL^2(\omega)}^2
		\\		
		&\eqreff*{equ:BDF2 norm equivalence}{\lesssim}
		\tau
		+
		\tau \Big( \tau\norm[\big]{\dtau\mm_h^1}_{\bL^2(\omega)}^2
		+
		\tau\sum_{j=1}^{n-1}
		\norm[\big]{\Dtau\mm_h^{j+1}}_{\bL^2(\omega)}^2 \Big)
		\\
		&\eqreff*{equ:magnetization stability}{\lesssim}
		\tau.
	\end{align*}
	This completes the proof.
\end{proof}


\begin{remark}[Quadratic constraint violation]\label{rem:quadratic defect}
	The estimates of Lemma~\ref{lem:unit length defect} only rely on the stability estimate~\eqref{equ:magnetization stability}. If the discrete solutions additionally satisfy
	\begin{align*}
		\norm[\big]{\vv_h^0}_{\bL^2(\omega)}^2
		+
		\tau^2\sum_{j=1}^{N-1}\norm[\big]{\dtau^2\mm_h^{j+1}}_{\bL^2(\omega)}^2
		\leq C_{\rm reg}
	\end{align*}
	with $C_{\rm reg}>0$ independent of $h$ and $\tau$, then the first estimate in the proof of Lemma~\ref{lem:unit length defect} yields $\norm[\big]{\III_h\bigl(\abs{\mm_h^n}^2\bigr)-1}_{L^1(\omega)}\lesssim\tau^2$ and hence $\norm[\big]{\abs{\mm_h^n}^2-1}_{L^1(\omega)}\lesssim h^2+\tau^2$ for all $0\leq n\leq N$; cf.~\cite[Remark~3.9]{AldFeiPra26} and \cite[Remark~9]{AldFeiPraSoen26_part3}. This is consistent with the quadratic decay of the nodal unit-length defect observed in Section~\ref{sec:num-exp}.
\end{remark}

\subsection{Compactness of the discrete interpolants}
\label{subsec:compactness}

In the following section we collect the consequences of the discrete stability estimates for the space-time interpolants and extract weakly convergent subsequences. 

Recall the definition of various interpolants in Section~\ref{sec:num alg}.
The following lemma collects available bounds and consistency properties of interpolants.

\begin{lemma}[Bounds and consistency of the magnetization interpolants]
	\label{lem:bounds interpolants m}
	There exists a constant, independent of $h$, $\tau$, and $N$, such that
	\begin{align}\label{equ:m interpolant bounds}
		&
		\norm[\big]{\mm_{h\tau}}_{\bH^1(\omega_T)}
		+
		\norm[\big]{\vv_{h\tau}^-}_{\bL^2(\omega_T)}
		+
		\norm[\big]{\mm_{h\tau}}_{L^\infty(0,T;\bH^1(\omega))}
		\nonumber
		\\
		&\qquad
		+
		\norm[\big]{\mm_{h\tau}^{\pm}}_{L^\infty(0,T;\bH^1(\omega))}
		+
		\norm[\big]{\widehat{\mm}_{h\tau}^{+}}_{L^\infty(0,T;\bH^1(\omega))}
		+
		\norm[\big]{\Pi_h(\widehat{\mm}_{h\tau}^{+})}_{L^\infty(0,T;\bH^1(\omega))}
		\lesssim 1.
	\end{align}
	Moreover, we have
	\begin{align}
		\label{equ:m interpolant consistency}
		&
		\norm[\big]{
			\mm_{h\tau}-\mm_{h\tau}^{-}
		}_{\bL^2(\omega_T)}
		+
		\norm[\big]{
			\mm_{h\tau}-\mm_{h\tau}^{+}
		}_{\bL^2(\omega_T)}
		\nonumber
		\\
		&\qquad
		+
		\norm[\big]{
			\mm_{h\tau}^{+}-\mm_{h\tau}^{-}
		}_{\bL^2(\omega_T)}
		+
		\norm[\big]{
			\widehat{\mm}_{h\tau}^{+}-\mm_{h\tau}^{-}
		}_{\bL^2(\omega_T)}
		\lesssim \tau,
	\end{align}
	and
	\begin{align}\label{equ:m predictor H1 consistency}
		\norm[\big]{\widehat{\mm}_{h\tau}^{+}-\mm_{h\tau}^{+}}_{L^2(0,T;\bH^1(\omega))}
		&\lesssim \tau^{1/2},
		\\
		\label{equ:projected predictor consistency}
		\norm[\big]{\Pi_h(\widehat{\mm}_{h\tau}^{+})-\widehat{\mm}_{h\tau}^{+}}_{\bL^2(\omega_T)}
		&\lesssim
		\tau^{1/2}.
	\end{align}
	The implicit constant depends on $C_{\mm,0}$, $C_{\si,0}$, $\const{C}{$f$}$, $\const{C}{$j$}$, $\alpha$, $\const{C}{ex}$, $\const{C}{$\pi$}$, $c$, $\beta$, $\beta'$, $D_*$, $T$, $\omega$, $\Omega$, $\Gamma_J$, and on $\gamma$.
\end{lemma}

\begin{proof}
	We first show~\eqref{equ:m interpolant bounds}. Note that the bounds for $\mm_{h\tau}^{\pm}$ follow directly from Lemma~\ref{lem:stability magnetization}. Next, \eqref{equ:magnetization stability} from Lemma~\ref{lem:stability magnetization}, the identities~\eqref{equ:BDF2 id 1}--\eqref{equ:BDF2 id 2}, and the norm equivalence \eqref{equ:BDF2 norm equivalence} imply
	\begin{align}
		\label{equ:backward derivative m bound}
		\tau
		\sum_{n=0}^{N-1}
		\norm[\big]{\dtau\mm_h^{n+1}}_{\bL^2(\omega)}^2
		\lesssim 1.
	\end{align}
	The identity~\eqref{equ:dtm dts cont disc} and estimate~\eqref{equ:backward derivative m bound}, together with the uniform $\bH^1(\omega)$-bound for $\mm_h^n$ in \eqref{equ:magnetization stability}, yield
	\[
	\norm[\big]{\mm_{h\tau}}_{\bH^1(\omega_T)}
	+
	\norm[\big]{\mm_{h\tau}}_{L^\infty(0,T;\bH^1(\omega))}
	\lesssim 1.
	\]
	Moreover, by the definition of $\widehat{\mm}_h^{n+1}$,
	we obtain the asserted bound for $\widehat{\mm}_{h\tau}^{+}$. The bound for $\vv_{h\tau}^{-}$ follows immediately from~\eqref{equ:magnetization stability}, while the $\bH^1(\omega)$-stability of $\Pi_h$ gives the corresponding bound for $\Pi_h(\widehat{\mm}_{h\tau}^{+})$. This proves~\eqref{equ:m interpolant bounds}.
	
	Next, we show~\eqref{equ:m interpolant consistency}. For $t\in[t_n,t_{n+1})$, observe that
	\begin{align*}
		\mm_{h\tau}(t)-\mm_h^n
		&=
		(t-t_n)\dtau\mm_h^{n+1},
		\\
		\mm_{h\tau}(t)-\mm_h^{n+1}
		&=
		-(t_{n+1}-t)\dtau\mm_h^{n+1}.
	\end{align*}
	Hence, we obtain
	\begin{align*}
		\norm[\big]{
			\mm_{h\tau}-\mm_{h\tau}^{\pm}
		}_{\bL^2(\omega_T)}^2
		&=
		\sum_{n=0}^{N-1}
		\int_{t_n}^{t_{n+1}}
		(t-t_n)^2 \norm[\big]{\dtau\mm_h^{n+1}}_{\bL^2(\omega)}^2 \d t
		\lesssim
		\tau^3
		\sum_{n=0}^{N-1}
		\norm[\big]{\dtau\mm_h^{n+1}}_{\bL^2(\omega)}^2
		\eqreff*{equ:backward derivative m bound}{\lesssim}
		\tau^2.
	\end{align*}
	Similarly, since $\mm_h^{n+1}-\mm_h^n= \tau\dtau\mm_h^{n+1}$, we have
	\[
	\norm[\big]{
		\mm_{h\tau}^{+}-\mm_{h\tau}^{-}
	}_{\bL^2(\omega_T)}
	\lesssim \tau.
	\]
	Furthermore, by our convention,
	$\widehat{\mm}_{h\tau}^{+}-\mm_{h\tau}^{-}=0$ on $[t_0,t_1)$,
	while, for $n\geq1$,
	\[
	\widehat{\mm}_h^{n+1}-\mm_h^n
	=
	\mm_h^n-\mm_h^{n-1}
	=
	\tau\dtau\mm_h^n.
	\]
	Thus, another application of~\eqref{equ:backward derivative m bound} proves~\eqref{equ:m interpolant consistency}.
	
	To establish~\eqref{equ:m predictor H1 consistency}, we use the fact that on the first time interval, $\mm_h^1-\widehat{\mm}_h^1 =
	\tau\vv_h^0$, while for $n\geq1$,
	\[
	\mm_h^{n+1}-\widehat{\mm}_h^{n+1}
	=
	\tau^2\dtau^2\mm_h^{n+1}.
	\]
	It then follows that
	\begin{align*}
		\norm[\big]{\mm_{h\tau}^{+}-\widehat{\mm}_{h\tau}^{+}}_{L^2(0,T;\bH^1(\omega))}^2
		&=
		\int_{t_0}^{t_1} \norm[\big]{\tau\vv_h^0}_{\bH^1(\omega)}^2 \d t
		+
		\sum_{n=1}^{N-1} \int_{t_n}^{t_{n+1}}
		\norm[\big]{\tau^2\dtau^2\mm_h^{n+1}}_{\bH^1(\omega)}^2
		\d t
		\\
		&=
		\tau^3 \norm[\big]{\vv_h^0}_{\bH^1(\omega)}^2
		+
		\tau^5 \sum_{n=1}^{N-1} \norm[\big]{\dtau^2\mm_h^{n+1}}_{\bH^1(\omega)}^2
		\eqreff*{equ:magnetization stability}{\lesssim}
		\tau,
	\end{align*}
	which proves~\eqref{equ:m predictor H1 consistency}.
	
	
	Finally, we prove~\eqref{equ:projected predictor consistency}.
	On the first time interval, the convention
	$\widehat{\mm}_h^1=\mm_h^0\in\mathcal{M}_h$ gives
	$\Pi_h(\widehat{\mm}_h^1)=\widehat{\mm}_h^1$.
	For $n=1,\ldots,N-1$ and $\vec{z}\in\mathcal{N}_h^\omega$,
	we have $\abs{\widehat{\mm}_h^{n+1}(\vec{z})}\geq1$
	by~\eqref{equ:nodal mhat geq one}. Hence, using
	$(r-1)^2\leq r^2-1$ for $r\geq1$, we obtain
	\begin{align*}
		\abs{
			\Pi_h\bigl(\widehat{\mm}_h^{n+1}\bigr)(\vec{z})
			-\widehat{\mm}_h^{n+1}(\vec{z})
		}^2
		&=\bigl(\abs{\widehat{\mm}_h^{n+1}(\vec{z})}-1\bigr)^2
		\leq
		\abs{\widehat{\mm}_h^{n+1}(\vec{z})}^2-1
		\\
		&\eqreff*{equ:predictor nodal identity}{=}
		2\abs{\mm_h^n(\vec{z})}^2
		-\abs{\mm_h^{n-1}(\vec{z})}^2
		+2\tau^2\abs{\dtau\mm_h^n(\vec{z})}^2-1
		\\
		&\eqreff*{equ:nodal m geq one}{\leq}
		2\bigl(\abs{\mm_h^n(\vec{z})}^2-1\bigr)
		+2\tau^2\abs{\dtau\mm_h^n(\vec{z})}^2.
	\end{align*}
	Since $\III_h(\abs{\mm_h^n}^2)-1$ is nonnegative,
	{the weighted summation identity~\eqref{equ:nodal weights} and} the norm equivalence~\eqref{equ:discrete nodal norm}
	{yield}
	\begin{align*}
		\norm[\big]{
			\Pi_h\bigl(\widehat{\mm}_h^{n+1}\bigr)
			-\widehat{\mm}_h^{n+1}
		}_{\bL^2(\omega)}^2
		&\lesssim
		\sum_{\vec{z}\in\mathcal{N}_h^\omega}{\mu_{\vec z}^\omega}
		\abs{
			\Pi_h\bigl(\widehat{\mm}_h^{n+1}\bigr)(\vec{z})
			-\widehat{\mm}_h^{n+1}(\vec{z})
		}^2
		\\
		&\lesssim
		\sum_{\vec{z}\in\mathcal{N}_h^\omega}{\mu_{\vec z}^\omega}
		(\abs{\mm_h^n(\vec{z})}^2-1)
		+\tau^2\sum_{\vec{z}\in\mathcal{N}_h^\omega}{\mu_{\vec z}^\omega}
		\abs{\dtau\mm_h^n(\vec{z})}^2
		\\
		&\lesssim
		\norm[\big]{\III_h\bigl(\abs{\mm_h^n}^2\bigr)-1}_{L^1(\omega)}
		+\tau^2\norm[\big]{\dtau\mm_h^n}_{\bL^2(\omega)}^2.
	\end{align*}
	The first term is bounded by $\lesssim\tau$
	by~\eqref{equ:nodal interpolated modulus defect}.
	For the second term, we use~\eqref{equ:backward derivative m bound}
	to obtain, for $n=1,\ldots,N-1$,
	\[
	\tau^2\norm[\big]{\dtau\mm_h^n}_{\bL^2(\omega)}^2
	\leq
	\tau\Bigl(
	\tau\sum_{j=0}^{N-1}
	\norm[\big]{\dtau\mm_h^{j+1}}_{\bL^2(\omega)}^2
	\Bigr)
	\lesssim \tau.
	\]
	Consequently,
	\[
	\norm[\big]{
		\Pi_h\bigl(\widehat{\mm}_h^{n+1}\bigr)
		-\widehat{\mm}_h^{n+1}
	}_{\bL^2(\omega)}^2
	\lesssim \tau
	\qquad\text{for all }n=0,\ldots,N-1.
	\]
	Integrating in time and using $N\tau=T$, we conclude that
	\begin{align*}
		\norm[\big]{
			\Pi_h(\widehat{\mm}_{h\tau}^{+})
			-\widehat{\mm}_{h\tau}^{+}
		}_{\bL^2(\omega_T)}^2
		&=
		\tau\sum_{n=0}^{N-1}
		\norm[\big]{
			\Pi_h\bigl(\widehat{\mm}_h^{n+1}\bigr)
			-\widehat{\mm}_h^{n+1}
		}_{\bL^2(\omega)}^2
		\lesssim N\tau^2
		=T\tau.
	\end{align*}
	Taking square roots proves~\eqref{equ:projected predictor consistency} and completes the proof of the lemma.
\end{proof}

\begin{lemma}[Bounds and consistency of the spin interpolants]
	\label{lem:bounds interpolants s}
	There exists a constant, independent of $h$, $\tau$, and $N$, such that
	\begin{align}\label{equ:s interpolant bounds}
		&
		\norm[\big]{\si_{h\tau}}_{L^\infty(0,T;\bL^2(\Omega))}
		+
		\norm[\big]{\si_{h\tau}^{\pm}}_{L^\infty(0,T;\bL^2(\Omega))}
		+
		\norm[\big]{\widehat{\si}_{h\tau}^{+}}_{L^\infty(0,T;\bL^2(\Omega))}
		\nonumber
		\\
		&\qquad
		+
		\norm[\big]{\si_{h\tau}}_{L^2(0,T;\bH^1(\Omega))}
		+
		\norm[\big]{\si_{h\tau}^{\pm}}_{L^2(0,T;\bH^1(\Omega))}
		+
		\norm[\big]{\widehat{\si}_{h\tau}^{+}}_{L^2(0,T;\bH^1(\Omega))}
		\nonumber
		\\
		&\qquad
		+
		\norm[\big]{{R_h^*(\partial_t\si_{h\tau})}}_{L^2(0,T;\widetilde{\bH}^{-1}(\Omega))}
		+
		\norm[\big]{{R_h^*(\Dtau\si_{h\tau}^{+})}}_{L^2(0,T;\widetilde{\bH}^{-1}(\Omega))}
		\lesssim 1.
	\end{align}
	Moreover,
	\begin{align} \label{equ:s interpolant Hminus consistency}
		&
		\norm[\big]{{R_h^*(\si_{h\tau}-\si_{h\tau}^{-})}}_{L^2(0,T;\widetilde{\bH}^{-1}(\Omega))}
		+
		\norm[\big]{{R_h^*(\si_{h\tau}-\si_{h\tau}^{+})}}_{L^2(0,T;\widetilde{\bH}^{-1}(\Omega))}
		\nonumber\\
		&\quad+\norm[\big]{{R_h^*(\widehat{\si}_{h\tau}^{+}-\si_{h\tau}^{-})}}_{L^2(0,T;\widetilde{\bH}^{-1}(\Omega))}
		\lesssim \tau,
	\end{align}
	and consequently
	\begin{align}\label{equ:s interpolant L2 consistency}
		&
		\norm[\big]{\si_{h\tau}-\si_{h\tau}^{-}}_{\bL^2(\Omega_T)}
		+
		\norm[\big]{\si_{h\tau}-\si_{h\tau}^{+}}_{\bL^2(\Omega_T)}
		+
		\norm[\big]{\widehat{\si}_{h\tau}^{+}-\si_{h\tau}^{-}}_{\bL^2(\Omega_T)}
		\lesssim \tau^{1/2}.
	\end{align}
	The implicit constant depends on $C_{\si,0}$, $\const{C}{$j$}$, $\beta$, $\beta'$, $D_*$, $D^*$, $T$, $\Omega$, $\Gamma_J$, and on $\gamma$.
\end{lemma}

\begin{proof}
	First, we show~\eqref{equ:s interpolant bounds}. Immediately from $\tau\norm[\big]{\si_h^0}_{\bH^1(\Omega)}^2\le TC_{\si,0}^2$ by~\eqref{equ:initial convergence s}, from~\eqref{equ:stab shn} and the fact that $\si_{h\tau}$ is a convex combination of $\si_h^n$ and
	$\si_h^{n+1}$ on each interval $[t_n,t_{n+1})$, we have
	\[
	\norm[\big]{\si_{h\tau}^{\pm}}_{L^\infty(0,T;\bL^2(\Omega))}
	+
	\norm[\big]{\si_{h\tau}^{\pm}}_{L^2(0,T;\bH^1(\Omega))}
	+
	\norm[\big]{\si_{h\tau}}_{L^\infty(0,T;\bL^2(\Omega))}
	+
	\norm[\big]{\si_{h\tau}}_{L^2(0,T;\bH^1(\Omega))}
	\lesssim 1.
	\]
	By~\eqref{equ:time reconstructions predictor s}, we also have
	\[
	\norm[\big]{\widehat{\si}_{h\tau}^{+}}_{L^\infty(0,T;\bL^2(\Omega))}
	+
	\norm[\big]{\widehat{\si}_{h\tau}^{+}}_{L^2(0,T;\bH^1(\Omega))}
	\lesssim 1.
	\]
	Noting~\eqref{equ:dtm dts cont disc}, we obtain
	\[
	\norm[\big]{{R_h^*(\partial_t\si_{h\tau})}}_{L^2(0,T;\widetilde{\bH}^{-1}(\Omega))}
	\eqreff*{equ:spin dt Hminus1}{\lesssim} 1.
	\]
	Furthermore, by~\eqref{equ:time reconstruction BDF2 s} we have
	\begin{align*}
		\norm[\big]{{R_h^*(\Dtau\si_{h\tau}^{+})}}_{L^2(0,T;\widetilde{\bH}^{-1}(\Omega))}^2
		&=
		\tau
		\norm[\big]{{R_h^*(\dtau\si_h^1)}}_{\widetilde{\bH}^{-1}(\Omega)}^2
		+
		\tau
		\sum_{n=1}^{N-1}
		\norm[\big]{{R_h^*(\Dtau\si_h^{n+1})}}_{\widetilde{\bH}^{-1}(\Omega)}^2
		\eqrefftwo{equ:spin dt Hminus1}{equ:BDF2 norm equivalence}{\lesssim}
		1.
	\end{align*}
	This proves~\eqref{equ:s interpolant bounds}.	
	Next, we show~\eqref{equ:s interpolant Hminus consistency}. For
	$t\in[t_n,t_{n+1})$, we have
	\begin{align*}
		\si_{h\tau}(t)-\si_h^n
		&=
		(t-t_n)\dtau\si_h^{n+1},
		\\
		\si_{h\tau}(t)-\si_h^{n+1}
		&=
		-(t_{n+1}-t)\dtau\si_h^{n+1}.
	\end{align*}
	It follows that
	\begin{align*}
		\norm[\big]{{R_h^*(\si_{h\tau}-\si_{h\tau}^{\pm})}}_{L^2(0,T;\widetilde{\bH}^{-1}(\Omega))}^2
		&\lesssim
		\tau^3
		\sum_{n=0}^{N-1}
		\norm[\big]{{R_h^*(\dtau\si_h^{n+1})}}_{\widetilde{\bH}^{-1}(\Omega)}^2
		\eqreff*{equ:spin dt Hminus1}{\lesssim}
		\tau^2.
	\end{align*}
	Moreover, by convention \eqref{equ:time reconstructions predictor s},
	$\widehat{\si}_{h\tau}^{+}-\si_{h\tau}^{-}=0$ on $[t_0,t_1)$,
	while, for $n\geq1$,
	\[
	\widehat{\si}_h^{n+1}-\si_h^n
	=
	\si_h^n-\si_h^{n-1}
	=
	\tau\dtau\si_h^n.
	\]
	Thus, similarly we obtain
	\begin{align*}
		\norm[\big]{{R_h^*(\widehat{\si}_{h\tau}^{+}-\si_{h\tau}^{-})}}_{L^2(0,T;\widetilde{\bH}^{-1}(\Omega))}^2
		&=
		\tau^3
		\sum_{n=1}^{N-1}
		\norm[\big]{{R_h^*(\dtau\si_h^n)}}_{\widetilde{\bH}^{-1}(\Omega)}^2
		\eqreff*{equ:spin dt Hminus1}{\lesssim}
		\tau^2.
	\end{align*}
	This proves~\eqref{equ:s interpolant Hminus consistency}.
	Finally, to obtain~\eqref{equ:s interpolant L2 consistency}, {we use the interpolation inequality~\eqref{equ:discrete Hminus1 interpolation}.}
	Since the unreconstructed differences appearing in
	\eqref{equ:s interpolant Hminus consistency} are uniformly bounded
	in $L^2(0,T;\bH^1(\Omega))$ by~\eqref{equ:s interpolant bounds}, after an integration in time and an application of the Cauchy--Schwarz inequality,~\eqref{equ:s interpolant L2 consistency} follows. This completes the proof.
\end{proof}


The preceding bounds and consistency estimates provide the compactness needed to extract convergent subsequences and identify common limits of the various time interpolants. We detail this in the following proposition.

\begin{proposition}[Extraction of convergent subsequences]
	\label{prop:compactness}
	Under the assumptions of Theorem~\ref{the:main convergence}, there exist
	a subsequence, not relabeled, and functions
	\[
	\mm
	\in
	\bH^1(\omega_T)
	\cap
	L^\infty(0,T;\bH^1(\omega))
	\]
	and
	\[
	\si
	\in
	L^\infty(0,T;\bL^2(\Omega))
	\cap
	L^2(0,T;\bH^1(\Omega))
	\cap
	H^1(0,T;\widetilde{\bH}^{-1}(\Omega))
	\]
	such that $\abs{\mm}=1$ almost everywhere in $\omega_T$ and, along
	the same subsequence, all convergences
	\begin{alignat}{2}
		\mm_{h\tau}
		&\rightharpoonup
		\mm
		&&\qquad
		\text{in }\bH^1(\omega_T),
		\label{equ:compactness m H1}
		\\
		\mm_{h\tau}
		&\to
		\mm
		&&\qquad
		\text{in }\bH^s(\omega_T)\quad \text{for all }\, 0<s<1,
		\label{equ:compactness m Hs}
		\\
		\mm_{h\tau},
		\mm_{h\tau}^{\pm},
		\widehat{\mm}_{h\tau}^{+},
		\Pi_h(\widehat{\mm}_{h\tau}^{+})
		&\to
		\mm
		&&\qquad
		\text{in }L^2(0,T;\bH^s(\omega))\quad \text{for all }\, 0<s<1,
		\label{equ:compactness m L2 Hs}
		\\
		\mm_{h\tau},
		\mm_{h\tau}^{\pm},
		\widehat{\mm}_{h\tau}^{+},
		\Pi_h(\widehat{\mm}_{h\tau}^{+})
		&\overset{*}{\rightharpoonup}
		\mm
		&&\qquad
		\text{in }L^\infty(0,T;\bH^1(\omega)),
		\label{equ:compactness m weak star}
		\\
		\mm_{h\tau},
		\mm_{h\tau}^{\pm},
		\widehat{\mm}_{h\tau}^{+},
		\Pi_h(\widehat{\mm}_{h\tau}^{+})
		&\rightharpoonup
		\mm
		&&\qquad
		\text{in }L^2(0,T;\bH^1(\omega)),
		\label{equ:compactness m weak H1}
		\\
		\mm_{h\tau},
		\mm_{h\tau}^{\pm},
		\widehat{\mm}_{h\tau}^{+},
		\Pi_h(\widehat{\mm}_{h\tau}^{+})
		&\to
		\mm
		&&\qquad
		\text{in }\bL^2(\omega_T),
		\label{equ:compactness m strong}
		\\
		\mm_{h\tau},
		\mm_{h\tau}^{\pm},
		\widehat{\mm}_{h\tau}^{+},
		\Pi_h(\widehat{\mm}_{h\tau}^{+})
		&\to
		\mm
		&&\qquad
		\text{pointwise a.e. in }\omega_T,
		\label{equ:compactness m pointwise}
		\\
		\partial_t\mm_{h\tau}
		&\rightharpoonup
		\partial_t\mm
		&&\qquad
		\text{in }\bL^2(\omega_T),
		\label{equ:compactness dt m}
		\\
		\vv_{h\tau}^{-}
		&\rightharpoonup
		\partial_t\mm
		&&\qquad
		\text{in }\bL^2(\omega_T)
		\label{equ:compactness v}
	\end{alignat}
	hold simultaneously as $(h,\tau)\to(0,0)$.
	Moreover, along the same subsequence, it holds that
	\begin{alignat}{2}
		\si_{h\tau},
		\si_{h\tau}^{\pm},
		\widehat{\si}_{h\tau}^{+}
		&\overset{*}{\rightharpoonup}
		\si
		&&\qquad
		\text{in }L^\infty(0,T;\bL^2(\Omega)),
		\label{equ:compactness s weak star}
		\\
		\si_{h\tau},
		\si_{h\tau}^{\pm},
		\widehat{\si}_{h\tau}^{+}
		&\rightharpoonup
		\si
		&&\qquad
		\text{in }L^2(0,T;\bH^1(\Omega)),
		\label{equ:compactness s weak H1}
		\\
		\si_{h\tau},
		\si_{h\tau}^{\pm},
		\widehat{\si}_{h\tau}^{+}
		&\to
		\si
		&&\qquad
		\text{in }\bL^2(\Omega_T),
		\label{equ:compactness s strong}
		\\
		{R_h^*(\partial_t\si_{h\tau})}
		&\rightharpoonup
		\partial_t\si
		&&\qquad
		\text{in }
		L^2(0,T;\widetilde{\bH}^{-1}(\Omega)),
		\label{equ:compactness dt s}
		\\
		{R_h^*(\Dtau\si_{h\tau}^{+})}
		&\rightharpoonup
		\partial_t\si
		&&\qquad
		\text{in }
		L^2(0,T;\widetilde{\bH}^{-1}(\Omega)).
		\label{equ:compactness BDF s}
	\end{alignat}
\end{proposition}

\begin{proof}
	Since the magnetization-related convergences follow from the same arguments as in \cite[Proposition~3.12]{AldFeiPra26}, we only sketch them.
	
	\smallskip
	{\bf Step~1 (Proof of limits~(\ref{equ:compactness m H1})--(\ref{equ:compactness dt m}) and the unit-length constraint).}
	
	Thanks to uniform boundedness~\eqref{equ:m interpolant bounds}, 
	the Eberlein--Šmulian theorem provides $\mm\in\bH^1(\omega_T)$ and a subsequence such that
	$\mm_{h\tau}\rightharpoonup\mm$ in $\bH^1(\omega_T)$.  This proves \eqref{equ:compactness m H1} and, as a consequence, \eqref{equ:compactness dt m}. The compact inclusions 
	$\bH^1(\omega_T)\Subset \bH^s(\omega_T)\Subset\bL^2(\omega_T)$ for all $s\in (0,1)$,
	imply \eqref{equ:compactness m Hs} and the strong convergence 
	\eqref{equ:compactness m strong} for $\mm_{h\tau}$.
	The continuous inclusions $\bH^1(\omega_T)\subset L^2(0,T;\bH^1(\omega))$ 
	and $\bH^s(\omega_T)\subset L^2(0,T;\bH^s(\omega))$ yield 
	$\mm_{h\tau}\rightharpoonup\mm$ in $L^2(0,T;\bH^1(\omega))$ and
	$\mm_{h\tau}\to\mm$ in $L^2(0,T;\bH^s(\omega))$.

	Let now $\mm_{h\tau}^\star\in\{\mm_{h\tau}, \mm_{h\tau}^\pm,\widehat{\mm}_{h\tau}^+,\Pi_h(\widehat{\mm}_{h\tau}^+)\}$. 
	The consistency estimates~\eqref{equ:m interpolant consistency}
	and~\eqref{equ:projected predictor consistency} show that $\mm_{h\tau}^{\pm}$,
	$\widehat{\mm}_{h\tau}^+$, and $\Pi_h(\widehat{\mm}_{h\tau}^+)$ converge to $\mm$
	in $\bL^2(\omega_T)$ as well, which completes~\eqref{equ:compactness m strong}.
	Together with the
	uniform bounds~\eqref{equ:m interpolant bounds}, the Banach--Alaoglu and the
	Eberlein--Šmulian theorems provide further subsequences converging weakly-$*$ in
	$L^\infty(0,T;\bH^1(\omega))$ and weakly in $L^2(0,T;\bH^1(\omega))$, with in principle different limits. Since
	$L^\infty(0,T;\bH^1(\omega))\subset L^2(0,T;\bH^1(\omega))\subset\bL^2(\omega_T)$
	with continuous inclusions, uniqueness of the limit identifies all of them
	with~$\mm$. This proves~\eqref{equ:compactness m weak star}--\eqref{equ:compactness m weak H1} and shows,
	in particular, that $\mm\in L^\infty(0,T;\bH^1(\omega))$. Upon extraction of a
	further subsequence, the pointwise convergence~\eqref{equ:compactness m pointwise} is a direct consequence
	of~\eqref{equ:compactness m strong}. 
	
	Finally, the interpolation inequality
	$\norm[\big]{\ww}_{\bH^s(\omega)}\lesssim \norm[\big]{\ww}_{\bL^2(\omega)}^{1-s}
	\norm[\big]{\ww}_{\bH^1(\omega)}^{s}$, valid for all $s\in (0,1)$, and Hölder's inequality
	in time with the exponents $1/(1-s)$ and $1/s$ give
	\[
		\norm[\big]{\ww}_{L^2(0,T;\bH^s(\omega))}^{2}
		\lesssim
		\norm[\big]{\ww}_{\bL^2(\omega_T)}^{2(1-s)}
		\norm[\big]{\ww}_{L^2(0,T;\bH^1(\omega))}^{2s}.
	\]
	Applying this to $\ww=\mm_{h\tau}^\star-\mm$, together with the uniform
	bounds~\eqref{equ:m interpolant bounds} and the strong
	convergences~\eqref{equ:compactness m strong} in $\bL^2(\omega_T)$,
	implies~\eqref{equ:compactness m L2 Hs} for all $\mm_{h\tau}^\star$.

	Moreover, Lemma~\ref{lem:unit length defect} gives
	\begin{align*}
		\norm[\big]{\abs{\mm}^2-1}_{L^1(\omega_T)}
		&\leq
		\norm[\big]{\abs{\mm}^2-\abs{\mm_{h\tau}^{-}}^2}_{L^1(\omega_T)}
		+
		\norm[\big]{\abs{\mm_{h\tau}^{-}}^2-1}_{L^1(\omega_T)}
		\\
		&\leq
		\bigl(
			\norm[\big]{\mm}_{\bL^2(\omega_T)}
			+\norm[\big]{\mm_{h\tau}^{-}}_{\bL^2(\omega_T)}
		\bigr)
		\norm[\big]{\mm-\mm_{h\tau}^{-}}_{\bL^2(\omega_T)}
		+
		\norm[\big]{\abs{\mm_{h\tau}^{-}}^2-1}_{L^1(\omega_T)}
		\\
		&\eqreff*{equ:modulus defect}{\lesssim}
		\norm[\big]{\mm-\mm_{h\tau}^{-}}_{\bL^2(\omega_T)}
		+
		T\bigl(h^2+\tau\bigr)
		\xrightarrow{(h,\tau)\to(0,0)}
		0 .
	\end{align*}
	Hence, $\abs{\mm}=1$ a.e. in $\omega_T$.
	
	\smallskip
	{\bf Step~2 (Proof of limits~(\ref{equ:compactness s weak star})--(\ref{equ:compactness dt s})).}
	{By~\eqref{equ:s interpolant bounds} and the definition of $R_h^*$, the sequences $\{\si_{h\tau}\}$ and $\{R_h^*\si_{h\tau}\}$ are bounded in $L^2(0,T;\bH^1(\Omega))$ and $H^1(0,T;\widetilde{\bH}^{-1}(\Omega))$, respectively. Indeed, $\norm[\big]{R_h^*\ww}_{\widetilde{\bH}^{-1}(\Omega)}\leq\norm[\big]{\ww}_{\bL^2(\Omega)}$ for all $\ww\in\bL^2(\Omega)$ by~\eqref{equ:H1 projection}--\eqref{equ:dual reconstruction}, and $\partial_t(R_h^*\si_{h\tau})=R_h^*(\partial_t\si_{h\tau})$. Extracting subsequences, let $\si$ be the weak limit of $\si_{h\tau}$ in $L^2(0,T;\bH^1(\Omega))$. For every $\zetaa\in L^2(0,T;\bH^1(\Omega))$, we have $R_h\zetaa\to\zetaa$ strongly therein by~\eqref{equ:H1 projection} and dominated convergence. Hence,
\[
    \int_0^T\dual{R_h^*\si_{h\tau}(t)}{\zetaa(t)}\d t
    =\int_0^T\scalarproduct{\si_{h\tau}(t)}{R_h\zetaa(t)}_\Omega\d t
    \to\int_0^T\scalarproduct{\si(t)}{\zetaa(t)}_\Omega\d t.
\]
This identifies the weak limit of the reconstructed sequence, so that
\begin{align}\label{equ:compactness reconstructed s}
    R_h^*\si_{h\tau}\rightharpoonup\si\quad\text{in }H^1(0,T;\widetilde{\bH}^{-1}(\Omega)),
\end{align}
which proves~\eqref{equ:compactness dt s}.

To obtain strong compactness of $\si_{h\tau}$, let $\delta\in (0,T)$. The Cauchy--Schwarz inequality and Fubini's theorem yield
\begin{align*}
    \int_0^{T-\delta}\norm[\big]{R_h^*(\si_{h\tau}(t+\delta)-\si_{h\tau}(t))}_{\widetilde{\bH}^{-1}(\Omega)}^2\d t
    &\leq\delta\int_0^{T-\delta}\int_t^{t+\delta}\norm[\big]{R_h^*(\partial_t\si_{h\tau}(r))}_{\widetilde{\bH}^{-1}(\Omega)}^2\d r\d t\\
    &\leq\delta^2\norm[\big]{R_h^*(\partial_t\si_{h\tau})}_{L^2(0,T;\widetilde{\bH}^{-1}(\Omega))}^2\lesssim\delta^2.
\end{align*}
Using~\eqref{equ:discrete Hminus1 interpolation} for $\si_{h\tau}(t+\delta)-\si_{h\tau}(t)\in\V_h^\Omega$ and applying the uniform $L^2(0,T;\bH^1(\Omega))$ bound in~\eqref{equ:s interpolant bounds}, we infer by the Cauchy--Schwarz inequality in time that
\[
    \int_0^{T-\delta}\norm[\big]{\si_{h\tau}(t+\delta)-\si_{h\tau}(t)}_{\bL^2(\Omega)}^2\d t\lesssim\delta.
\]
By the compact inclusion $\bH^1(\Omega)\Subset\bL^2(\Omega)$ and Simon's time-translation criterion~\cite[Theorem~5]{Simon86}, upon extraction of a further subsequence, $\si_{h\tau}\to\si$ strongly in $\bL^2(\Omega_T)$.} The consistency
	estimate~\eqref{equ:s interpolant L2 consistency} and the triangle inequality
	then give
	\[
		\norm[\big]{\si_{h\tau}^{\pm}-\si_{h\tau}}_{\bL^2(\Omega_T)}
		+
		\norm[\big]{\widehat{\si}_{h\tau}^{+}-\si_{h\tau}}_{\bL^2(\Omega_T)}
		\eqreff*{equ:s interpolant L2 consistency}{\lesssim}
		\tau^{1/2}
		\xrightarrow{(h,\tau)\to(0,0)}
		0,
	\]
	so that $\si_{h\tau}^{\pm}\to\si$ and $\widehat{\si}_{h\tau}^{+}\to\si$ in
	$\bL^2(\Omega_T)$ as well. This proves~\eqref{equ:compactness s strong}.

	It remains to identify the weak and weak-$*$ limits of the piecewise constant
	interpolants. By~\eqref{equ:s interpolant bounds}, the sequences
	$\{\si_{h\tau}\}$, $\{\si_{h\tau}^{\pm}\}$, and $\{\widehat{\si}_{h\tau}^{+}\}$
	are bounded in the reflexive space $L^2(0,T;\bH^1(\Omega))$ as well as in
	$L^\infty(0,T;\bL^2(\Omega))$, the latter being the dual of the separable space
	$L^1(0,T;\bL^2(\Omega))$. Successive application of the Eberlein--Šmulian and of
	the Banach--Alaoglu theorems therefore provides further subsequences converging
	weakly, resp.\ weakly-$*$, in these spaces, with in principle different limits. Since
	both $L^2(0,T;\bH^1(\Omega))$ and $L^\infty(0,T;\bL^2(\Omega))$ are continuously
	embedded into $\bL^2(\Omega_T)$, uniqueness of the limit identifies all of them
	with the limit $\si$. This
	proves~\eqref{equ:compactness s weak star}--\eqref{equ:compactness s weak H1} and
	shows, in particular, that
	$\si\in L^\infty(0,T;\bL^2(\Omega))\cap L^2(0,T;\bH^1(\Omega))\cap
	H^1(0,T;\widetilde{\bH}^{-1}(\Omega))$, as asserted.
	
	\smallskip
	{\bf Step~3 (Proof of the BDF2 derivatives limits (\ref{equ:compactness v}) and (\ref{equ:compactness BDF s})).}
	As in Step~1, we extract a further subsequence such that
	$\vv_{h\tau}^{-}\rightharpoonup\vv$ in $\bL^2(\omega_T)$. To show that
	$\vv=\partial_t\mm$, we use a density argument. From~\eqref{equ:compactness dt m},
	$\partial_t\mm_{h\tau}\rightharpoonup\partial_t\mm$ in $\bL^2(\omega_T)$ and hence
	\[
		\scalarproduct{\vv_{h\tau}^{-}-\partial_t\mm_{h\tau}}{\pphi}_{\omega_T}
		\xrightarrow{(h,\tau)\to(0,0)}
		\scalarproduct{\vv-\partial_t\mm}{\pphi}_{\omega_T}
		\quad\text{for all}\quad
		\pphi\in\bL^2(\omega_T).
	\]
	Equation~\eqref{equ:BDF2 id 3}, together with
	$\vv_{h\tau}^{-}|_{[t_0,t_1)}=\vv_h^0=\dtau\mm_h^1=\partial_t\mm_{h\tau}|_{[t_0,t_1)}$
	and the compact support of $\pphi$ in $\omega_T$, yields
	\begin{equation}\label{equ: v = dtm}
		\int_{0}^{T}\scalarproduct{\vv_{h \tau}^{-}(t) - \partial_t\mm_{h\tau}(t)}{\pphi(t)}_{\omega}\d t
		\xrightarrow{(h,\tau)\to(0,0)}0
		\quad \text{for all} \quad \pphi \in \boldsymbol{C}^{\infty}_c(\omega_T),
	\end{equation}
	and hence, by the density
	of $\boldsymbol{C}^{\infty}_c(\omega_T)$ in $\bL^2(\omega_T)$, for all
	$\pphi\in\bL^2(\omega_T)$; see \cite[Proposition~3.12]{AldFeiPra26} for the
	detailed computation. Therefore, $\vv=\partial_t\mm$,
	which proves~\eqref{equ:compactness v}.

	In a similar fashion, \eqref{equ:s interpolant bounds} guarantees the existence
	of $\boldsymbol{\sigma}\in L^2(0,T;\widetilde{\bH}^{-1}(\Omega))$ such that
	${R_h^*(\Dtau\si_{h\tau}^{+})}\rightharpoonup\boldsymbol{\sigma}$ in
	$L^2(0,T;\widetilde{\bH}^{-1}(\Omega))$. To show that $\boldsymbol{\sigma}=\partial_t\si$, we use
	the same density argument as above, which we now carry out in detail.
	For readability, the duality pairing between $\widetilde{\bH}^{-1}(\Omega)$ and $\bH^1(\Omega)$ is written without subscripts throughout the following argument. From~\eqref{equ:compactness dt s},
	${R_h^*(\partial_t\si_{h\tau})}\rightharpoonup\partial_t\si$ in
	$L^2(0,T;\widetilde{\bH}^{-1}(\Omega))$ and hence
	\[
		\int_{0}^{T}
		\dual{{R_h^*(\Dtau\si_{h \tau}^{+}(t) - \partial_t\si_{h\tau}(t))}}{\zetaa(t)}\d t
		\xrightarrow{(h,\tau)\to(0,0)}
		\int_{0}^{T}
		\dual{\boldsymbol{\sigma}(t) - \partial_t\si(t)}{\zetaa(t)}\d t
	\]
	for all $\zetaa\in L^2(0,T;\bH^1(\Omega))$. It thus remains to show that
	\begin{equation}\label{equ: sigma = dts}
		\int_{0}^{T}
		\dual{{R_h^*(\Dtau\si_{h \tau}^{+}(t) - \partial_t\si_{h\tau}(t))}}{\zetaa(t)}\d t
		\xrightarrow{(h,\tau)\to(0,0)}0
		\quad\text{for all}\quad
		\zetaa\in \boldsymbol{C}_c^\infty\bigl((0,T);\bH^1(\Omega)\bigr).
	\end{equation}
	We consider as a test space
	$\boldsymbol{C}^{\infty}_c\bigl((0,T);\bH^1(\Omega)\bigr)$, which is a dense subset of $L^2(0,T;\bH^1(\Omega))$. Unlike for the
	magnetization case proved above, the set $\boldsymbol{C}^{\infty}_c(\Omega_T)$ would not suffice
	here, since its closure in $L^2(0,T;\bH^1(\Omega))$ is only
	$L^2(0,T;\bH^1_0(\Omega))$, on which the boundary contributions contained in
	$\widetilde{\bH}^{-1}(\Omega)=\bH^1(\Omega)^{*}$ vanish identically.

	Let thus $\zetaa\in\boldsymbol{C}^{\infty}_c((0,T);\bH^1(\Omega))$, extended by
	zero outside $(0,T)$. The convention $\Dtau\si_h^1\coloneqq\dtau\si_h^1$ gives
	$\Dtau\si_{h\tau}^{+}|_{[t_0,t_1)}=\Dtau\si_h^1=\dtau\si_h^1
	=\partial_t\si_{h\tau}|_{[t_0,t_1)}$, whereas, for every $n=1,\ldots,N-1$,
	\[
		\Dtau \si_h^{n+1} -\dtau \si_h^{n+1}
		\eqreff{equ: Dtau}{=}
		\frac{3}{2}\dtau \si_h^{n+1} - \frac{1}{2}\dtau \si_h^n - \dtau \si_h^{n+1}
		= \frac{1}{2}\bigl(\dtau \si_h^{n+1} - \dtau \si_h^n\bigr).
	\]
	Recalling that $\partial_t\si_{h\tau}(t)=\dtau\si_h^{n+1}$, and hence
	$\partial_t\si_{h\tau}(t-\tau)=\dtau\si_h^{n}$, for $t\in[t_n,t_{n+1})$ with
	$n\geq1$, we thus obtain, by linearity of $R_h^*$,
	\begin{align}
		\int_{0}^{T}&\dual{{R_h^*(\Dtau\si_{h \tau}^{+}(t) - \partial_t\si_{h\tau}(t))}}{\zetaa(t)}\d t
		=
		\sum_{n=1}^{N-1}\int_{t_n}^{t_{n+1}}
		\dual{{R_h^*(\Dtau\si_h^{n+1}-\dtau\si_h^{n+1})}}{\zetaa(t)}\d t
		\nonumber
		\\
		&=
		\frac12\sum_{n=1}^{N-1}\int_{t_n}^{t_{n+1}}
		\dual{{R_h^*(\dtau\si_h^{n+1}-\dtau\si_h^{n})}}{\zetaa(t)}\d t
		\nonumber\\
		&=\frac12\sum_{n=1}^{N-1}\int_{t_n}^{t_{n+1}}
		\dual{{R_h^*(\partial_t\si_{h\tau}(t)-\partial_t\si_{h\tau}(t-\tau))}}{\zetaa(t)}\d t
		\nonumber
		\\
		&=
		\frac12\int_{t_1}^{T}\dual{{R_h^*(\partial_t\si_{h\tau}(t))}}{\zetaa(t)}\d t
		-
		\frac12\int_{t_1}^{T}\dual{{R_h^*(\partial_t\si_{h\tau}(t-\tau))}}{\zetaa(t)}\d t ,
		\label{equ: BDF2 s splitting}
	\end{align}

	We consider the two terms in~\eqref{equ: BDF2 s splitting} separately. For the
	first one, the Cauchy--Schwarz inequality and~\eqref{equ:s interpolant bounds}
	give
	\begin{align*}
		\biggl|\int_{0}^{t_1}\dual{{R_h^*(\partial_t\si_{h\tau}(t))}}{\zetaa(t)}\d t\biggr|&
		\leq
		\norm[\big]{{R_h^*(\partial_t\si_{h\tau})}}_{L^2(0,T;\widetilde{\bH}^{-1}(\Omega))}
		\norm[\big]{\zetaa}_{L^2(0,t_1;\bH^1(\Omega))}\\
		&
		\eqreff*{equ:s interpolant bounds}{\lesssim}
		t_1^{1/2}\norm[\big]{\zetaa}_{L^\infty(0,T;\bH^1(\Omega))}
		\lesssim
		\tau^{1/2}\norm[\big]{\zetaa}_{L^\infty(0,T;\bH^1(\Omega))},
	\end{align*}
	so that, by~\eqref{equ:compactness dt s},
	\begin{align}
		\int_{t_1}^{T}\dual{{R_h^*(\partial_t\si_{h\tau}(t))}}{\zetaa(t)}\d t
		&=
		\int_{0}^{T}\dual{{R_h^*(\partial_t\si_{h\tau}(t))}}{\zetaa(t)}\d t
		-
		\int_{0}^{t_1}\dual{{R_h^*(\partial_t\si_{h\tau}(t))}}{\zetaa(t)}\d t
		\nonumber
		\\
		&\xrightarrow{(h,\tau)\to(0,0)}
		\int_{0}^{T}\dual{\partial_t\si(t)}{\zetaa(t)}\d t .
		\label{equ: BDF2 s first term}
	\end{align}
	For the second one,	since $\zetaa(t)=0$ for $t>T$, it holds that
	\begin{align*}
		\int_{t_1}^{T}\dual{{R_h^*(\partial_t\si_{h\tau}(t-\tau))}}{\zetaa(t)}\d t
		&=\int_{0}^{T-\tau}\dual{{R_h^*(\partial_t\si_{h\tau}(r))}}{\zetaa(r+\tau)}\d r\\
		&=\int_{0}^{T}\dual{{R_h^*(\partial_t\si_{h\tau}(t))}}{\zetaa(t+\tau)}\d t .
	\end{align*}
	Since $\zetaa$ is uniformly continuous as a map from $[0,T]$ to $\bH^1(\Omega)$,
	we have $\zetaa(\cdot+\tau)\to\zetaa$ strongly in $L^2(0,T;\bH^1(\Omega))$ as
	$\tau\to0$. Hence, the Cauchy--Schwarz inequality
	and~\eqref{equ:s interpolant bounds} yield
	\begin{align*}
		\biggl|\int_{0}^{T}\dual{{R_h^*(\partial_t\si_{h\tau}(t))}}{\zetaa(t+\tau)-\zetaa(t)}\d t\biggr|
		&\leq
		\norm[\big]{{R_h^*(\partial_t\si_{h\tau})}}_{L^2(0,T;\widetilde{\bH}^{-1}(\Omega))}
		\norm[\big]{\zetaa(\cdot+\tau)-\zetaa}_{L^2(0,T;\bH^1(\Omega))}\\
		&\xrightarrow{(h,\tau)\to(0,0)}0,
	\end{align*}
	so that, using~\eqref{equ:compactness dt s} once more,
	\begin{align}
		\int_{t_1}^{T}\dual{{R_h^*(\partial_t\si_{h\tau}(t-\tau))}}{\zetaa(t)}\d t
		\xrightarrow{(h,\tau)\to(0,0)}
		\int_{0}^{T}\dual{\partial_t\si(t)}{\zetaa(t)}\d t .
		\label{equ: BDF2 s second term}
	\end{align}
	Plugging~\eqref{equ: BDF2 s first term} and~\eqref{equ: BDF2 s second term}
	into~\eqref{equ: BDF2 s splitting}, we conclude~\eqref{equ: sigma = dts} first for
	all $\zetaa\in\boldsymbol{C}^{\infty}_c((0,T);\bH^1(\Omega))$ and then, by
	density, for all $\zetaa\in L^2(0,T;\bH^1(\Omega))$. Therefore,
	$\boldsymbol{\sigma}=\partial_t\si$, which proves~\eqref{equ:compactness BDF s}. This concludes the proof.
\end{proof}


\subsection{Identification of the limit with a weak solution}

In this section we will conclude the proof of Theorem~\ref{the:main convergence} by
showing that the pair $(\mm,\si)$ from Proposition~\ref{prop:compactness} is a weak
solution of the SDLLG system~\eqref{equ:SDLLG} in the sense of
Definition~\ref{def:weak solution}. Throughout, all convergences are understood along
the subsequence extracted in Proposition~\ref{prop:compactness}.

\begin{proof}[{\bfseries Proof of Theorem~\ref{the:main convergence} (Limit $(\mm,\si)$ satisfies Definition~\ref{def:weak solution}\bf(\ref{def: weak solution i})--(\ref{def: weak solution ii}))}]
	Proposition \ref{prop:compactness} shows 
	$\mm\in\bH^1(\omega_T)\cap L^\infty(0,T;\bH^1(\omega))$ with $\abs{\mm}=1$ a.e.\
	in $\omega_T$, as well as
	$\si\in L^2(0,T;\bH^1(\Omega))\cap L^\infty(0,T;\bL^2(\Omega))\cap
	H^1(0,T;\widetilde{\bH}^{-1}(\Omega))$. It thus only remains to verify the
	initial conditions.

	Since $\mm_{h\tau}\rightharpoonup\mm$ in $\bH^1(\omega_T)$
	by~\eqref{equ:compactness m H1}, and since the trace operator $\mathrm{tr}:\bH^1(\omega_T)\to\bH^{1/2}(\omega)$, $\uu\mapsto\uu(0)$, is linear and bounded, hence weakly continuous, it follows that $\mm_{h\tau}(0)\rightharpoonup\mm(0)$ in
	$\bH^{1/2}(\omega)$. By~\eqref{equ:initial convergence m}, it holds that
	$\mm_{h\tau}(0)=\mm_h^0\rightharpoonup\mm_0$ in $\bH^1(\omega)$ and hence in
	$\bH^{1/2}(\omega)$. It follows that $\mm(0)=\mm_0$ in the sense of traces, which
	verifies Definition~\ref{def:weak solution}(\ref{def: weak solution i}).

	As for the spin accumulation, we argue analogously. {By~\eqref{equ:compactness reconstructed s}, $R_h^*\si_{h\tau}\rightharpoonup\si$ in $H^1(0,T;\widetilde{\bH}^{-1}(\Omega))$.} Since
	$H^1(0,T;\widetilde{\bH}^{-1}(\Omega))\subset C([0,T];\widetilde{\bH}^{-1}(\Omega))$
	with continuous inclusion, the evaluation operator
	$H^1(0,T;\widetilde{\bH}^{-1}(\Omega))\to\widetilde{\bH}^{-1}(\Omega)$,
	$\uu\mapsto\uu(0)$, is linear and bounded, hence weakly continuous, so that
	${R_h^*\si_{h\tau}(0)}\rightharpoonup\si(0)$ in $\widetilde{\bH}^{-1}(\Omega)$.
	Furthermore, by~\eqref{equ:initial convergence s} and~\eqref{equ:H1 projection}, for every $\zetaa\in\bH^1(\Omega)$, we have $\dual{R_h^*\si_{h\tau}(0)}{\zetaa}=\scalarproduct{\si_h^0}{R_h\zetaa}_\Omega\to\scalarproduct{\si_0}{\zetaa}_\Omega$. It follows that $\si(0)=\si_0$ in the sense of
	traces, which verifies Definition~\ref{def:weak solution}(\ref{def: weak solution ii}).
\end{proof}

	\begin{proof}[\bfseries Proof of Theorem~\ref{the:main convergence} (Limit $(\mm,\si)$ satisfies Definition~\ref{def:weak solution}\bf(\ref{def: weak solution iii}))]
	The proof is split into four steps, where the core argument is summarized in
	Step~1.

	\smallskip
	{\bf Step~1 (Discrete variational formulation).}
	Let $\phii\in\boldsymbol{C}^\infty({{\overline{\omega_T}}})$ be an arbitrary test
	function and define $\zetaa_h\colon[0,T]\to\V_h^\omega$ by
	\begin{align}\label{equ:test function m}
		\zetaa_h(t)
		\coloneqq
		\III_h\bigl[\widehat{\mm}_{h\tau}^{+}(t)\times\phii(t)\bigr]
		\qquad\text{for all }t\in[0,T].
	\end{align}
	For all $n=0,\ldots,N-1$, all $t\in[t_n,t_{n+1})$, and all
	$\zz\in\mathcal{N}_h^\omega$, we have
	$\zetaa_h(t)(\zz)\cdot\widehat{\mm}_h^{n+1}(\zz)=0$, which implies that
	$\zetaa_h(t)\in\mathcal K_h(\widehat{\mm}_h^{n+1})$. Recalling the convention
	$\widehat{\mm}_h^1\coloneqq\mm_h^0$ and $\widehat{\si}_h^1\coloneqq\si_h^0$, and
	$\boldsymbol{F}_h^{n+1}$ from~\eqref{equ:discrete F}, Algorithm~\ref{alg: full discr} together
	with~\eqref{equ:first step m1} and~\eqref{equ:gen step mn1} yields, for all
	$n=0,\ldots,N-1$ and all $t\in[t_n,t_{n+1})$,
	\begin{align}\label{equ:pointwise discrete m}
		\alpha\scalarproduct{\vv_h^n}{\zetaa_h(t)}_{\omega}
		+
		\scalarproduct{\widehat{\mm}_h^{n+1}\times\vv_h^n}{\zetaa_h(t)}_{\omega}
		+
		\const{C}{ex}\scalarproduct{\nabla\mm_h^{n+1}}{\nabla\zetaa_h(t)}_{\omega}
		=
		\scalarproduct{\boldsymbol{F}_h^{n+1}}{\zetaa_h(t)}_{\omega}.
	\end{align}
	We integrate~\eqref{equ:pointwise discrete m} over $[t_n,t_{n+1})$ and sum these
	integrals over $n=0,\ldots,N-1$. Recalling the definitions of the interpolants,
	this leads to
	\begin{align}\label{equ:discrete var form m}
		&\alpha\int_0^T\scalarproduct{\vv_{h\tau}^{-}(t)}{\zetaa_h(t)}_{\omega}\d t
		+
		\int_0^T\scalarproduct{\widehat{\mm}_{h\tau}^{+}(t)\times\vv_{h\tau}^{-}(t)}{\zetaa_h(t)}_{\omega}\d t
		+
		\const{C}{ex}\int_0^T\scalarproduct{\nabla\mm_{h\tau}^{+}(t)}{\nabla\zetaa_h(t)}_{\omega}\d t
		\nonumber
		\\
		&\qquad
		=
		\int_0^T\scalarproduct{\pp_h\bigl(\widehat{\mm}_{h\tau}^{+}(t)\bigr)}{\zetaa_h(t)}_{\omega}\d t
		+
		\int_0^T\scalarproduct{\ff_{h\tau}^{+}(t)}{\zetaa_h(t)}_{\omega}\d t
		+
		c\int_0^T\scalarproduct{\widehat{\si}_{h\tau}^{+}(t)}{\zetaa_h(t)}_{\omega}\d t .
	\end{align}
	The subsequent steps consider the convergence of each of these six integrals
	separately and we will prove that the limit equation reads
	\begin{align}\label{equ:limit var form m}
		&\alpha\int_0^T\scalarproduct{\partial_t\mm(t)\times\mm(t)}{\phii(t)}_{\omega}\d t
		+
		\int_0^T\scalarproduct{\partial_t\mm(t)}{\phii(t)}_{\omega}\d t\nonumber\\
		&
		+
		\const{C}{ex}\int_0^T\scalarproduct{\nabla\mm(t)\times\mm(t)}{\nabla\phii(t)}_{\omega}\d t
		=
		\int_0^T\scalarproduct{\pp(\mm(t))\times\mm(t)}{\phii(t)}_{\omega}\d t\nonumber\\
		&
		+
		\int_0^T\scalarproduct{\ff(t)\times\mm(t)}{\phii(t)}_{\omega}\d t
		+
		c\int_0^T\scalarproduct{\si(t)\times\mm(t)}{\phii(t)}_{\omega}\d t .
	\end{align}
	Overall, we have then proved that~\eqref{equ:limit var form m} holds for all
	$\phii\in\boldsymbol{C}^\infty({\overline{\omega_T}})$ and hence, by density, for
	all $\phii\in\bH^1(\omega_T)$. Consequently, this verifies
	Definition~\ref{def:weak solution}(\ref{def: weak solution iii}).
	\smallskip
	\noindent

	{\bf Step~2 (Auxiliary convergence results).}
	We recall a few auxiliary results which will be used in the subsequent analysis. For further details, see~\cite{AldFeiPra26}.

	{For all continuous, elementwise $\bH^2$ functions $\ww$ on $\overline\omega$, the nodal interpolant $\III_h\ww\in\V_h^\omega$ satisfies, for each $K\in\mathcal T_h^\omega$,}
	\begin{align}\label{equ:interpolation estimate}
		\norm[\big]{\ww-\III_h\ww}_{\bL^2({K})}
		+
		{h_K}\norm[\big]{\nabla(\ww-\III_h\ww)}_{\bL^2({K})}
		\lesssim
		{h_K^2}\norm[\big]{D^2\ww}_{\bL^2({K})};
	\end{align}
	see, e.g., \cite[Theorem~4.4.4 and Theorem~4.4.20]{BreSco08}. Moreover, for all
	$\pphi_h\in\V_h^\omega$ and all $\phii\in\boldsymbol{C}^\infty(\overline{\omega})$, it holds
	that
	\begin{align}\label{equ:elementwise D2}
		\sum_{K\in\mathcal{T}_h^\omega}
		\norm[\big]{D^2(\pphi_h\times\phii)}_{\bL^2(K)}^2
		\lesssim
		\norm[\big]{\pphi_h}_{\bH^1(\omega)}^2
		\norm[\big]{\phii}_{\mathbf{W}^{2,\infty}(\omega)}^2 .
	\end{align}
	Combining~\eqref{equ:interpolation estimate}--\eqref{equ:elementwise D2}
	elementwise with the uniform bound~\eqref{equ:m interpolant bounds}, we obtain
	\begin{align}\label{equ:interpolation error m}
		\norm[\big]{(\III_h-\mathbf{Id})\bigl[\widehat{\mm}_{h\tau}^{+}\times\phii\bigr]}_{L^2(0,T;\bH^1(\omega))}
		\lesssim
		h\norm[\big]{\widehat{\mm}_{h\tau}^{+}}_{L^2(0,T;\bH^1(\omega))}
		\norm[\big]{\phii}_{\mathbf{W}^{2,\infty}(\omega_T)}
		\eqreff*{equ:m interpolant bounds}{\lesssim}
		h
		\xrightarrow{(h,\tau)\to(0,0)}0 .
	\end{align}
	Next, fix $3/4\leq s<1$. Since $\bH^s(\omega)\subset\bL^4(\omega)$ and
	$\bH^1(\omega)\subset\bL^4(\omega)$ with continuous inclusions, $\mm_{h\tau}\to\mm$ in $L^2(0,T;\bH^s(\omega))$ by
	\eqref{equ:compactness m L2 Hs} and the uniform bounds~\eqref{equ:m interpolant bounds} give
	\begin{align}\label{equ:m L4}
		\mm_{h\tau}^\star\to\mm
		\quad\text{in }L^2(0,T;\bL^4(\omega))
	\end{align}
	for all
	$\mm_{h\tau}^\star\in\{\mm_{h\tau},\mm_{h\tau}^{\pm},\widehat{\mm}_{h\tau}^{+},
	\Pi_h(\widehat{\mm}_{h\tau}^{+})\}$. In particular, together with~\eqref{equ:m interpolant bounds},
	\eqref{equ:interpolation error m} and $\abs{\mm}=1$ a.e.\ in $\omega_T$,
	this implies
	\begin{align}\label{equ:test function m convergence}
		\zetaa_h\to\mm\times\phii
		\quad\text{in }\bL^2(\omega_T),
		\qquad
		\bigl(\widehat{\mm}_{h\tau}^{+}\times\phii\bigr)\times\widehat{\mm}_{h\tau}^{+}
		\to
		(\mm\times\phii)\times\mm
		\quad\text{in }\bL^2(\omega_T),
	\end{align}
	where for the second convergence we used the generalized Hölder inequality with
	the exponents $1/2=1/4+1/4$.

	\smallskip
	{\bf Step~3 (Convergence of the first, second and third integrals
	in~(\ref{equ:discrete var form m})--(\ref{equ:limit var form m})).}
	These three integrals involve the magnetization only, and their treatment is
	identical to the one in~\cite[Theorem~2.2]{AldFeiPra26} or \cite[Theorem~3]{AldFeiPraSoen26_part3}, therefore we only sketch the proof here for completeness, and refer to the cited papers for the
	details. We recall that
	$\scalarproduct{\uu_{h\tau}}{\ww_{h\tau}}\to\scalarproduct{\uu}{\ww}$ provided
	that $\uu_{h\tau}\rightharpoonup\uu$ and $\ww_{h\tau}\to\ww$ in $\bL^2$. Applied
	to the present setting with the convergence
	properties~\eqref{equ:test function m convergence} and the weak convergences from
	Proposition~\ref{prop:compactness}, and using the identity
	$\scalarproduct{\uu}{\ww\times\zetaa}=\scalarproduct{\uu\times\ww}{\zetaa}$, this
	shows
	\begin{align}\label{equ:limit first integral m}
		\alpha\int_0^T\scalarproduct{\vv_{h\tau}^{-}(t)}{\zetaa_h(t)}_{\omega}\d t
		\xrightarrow{(h,\tau)\to(0,0)}
		\alpha\int_0^T\scalarproduct{\partial_t\mm(t)\times\mm(t)}{\phii(t)}_{\omega}\d t .
	\end{align}
	For the second integral, we first use~\eqref{equ:interpolation error m} to replace
	$\zetaa_h$ by $\widehat{\mm}_{h\tau}^{+}\times\phii$ up to an error of order $h$,
	and then use the identity
	$\scalarproduct{\uu\times\ww}{\uu\times\zetaa}
	=\scalarproduct{\ww}{(\uu\times\zetaa)\times\uu}$ together
	with~\eqref{equ:compactness v} and the second convergence
	in~\eqref{equ:test function m convergence} to obtain
	\[
		\int_0^T\scalarproduct{\widehat{\mm}_{h\tau}^{+}(t)\times\vv_{h\tau}^{-}(t)}{\zetaa_h(t)}_{\omega}\d t
		\xrightarrow{(h,\tau)\to(0,0)}
		\int_0^T\scalarproduct{\partial_t\mm(t)}{(\mm(t)\times\phii(t))\times\mm(t)}_{\omega}\d t .
	\]
	Since $\abs{\mm}=1$ a.e.\ in $\omega_T$ and hence
	$\mm\cdot\partial_t\mm=\frac12\partial_t\abs{\mm}^2=0$ a.e.\ in $\omega_T$, it holds that
	\begin{align}\label{equ:limit second integral m}
		\int_0^T\scalarproduct{\widehat{\mm}_{h\tau}^{+}(t)\times\vv_{h\tau}^{-}(t)}{\zetaa_h(t)}_{\omega}\d t
		\xrightarrow{(h,\tau)\to(0,0)}
		\int_0^T\scalarproduct{\partial_t\mm(t)}{\phii(t)}_{\omega}\d t .
	\end{align}
	For the third integral, we first
	use~\eqref{equ:interpolation error m} to
	replace $\zetaa_h$ by $\widehat{\mm}_{h\tau}^{+}\times\phii$ up to an error of
	order $h$, and then split
	$\nabla(\widehat{\mm}_{h\tau}^{+}\times\phii)
	=\nabla\widehat{\mm}_{h\tau}^{+}\times\phii
	+\widehat{\mm}_{h\tau}^{+}\times\nabla\phii$. 	Since $\nabla\widehat{\mm}_{h\tau}^{+}\cdot\bigl(\nabla\widehat{\mm}_{h\tau}^{+}\times\phii\bigr)=0$
	pointwise a.e.\ in $\omega_T$, the consistency estimate~\eqref{equ:m predictor H1 consistency} gives
	\begin{align*}
		\Bigl|
			\int_0^T\scalarproduct{\nabla\mm_{h\tau}^{+}(t)}{\nabla\widehat{\mm}_{h\tau}^{+}(t)\times\phii(t)}_{\omega}\d t
		\Bigr|
		&=
		\Bigl|
			\int_0^T\scalarproduct{\nabla\mm_{h\tau}^{+}(t)-\nabla\widehat{\mm}_{h\tau}^{+}(t)}{\nabla\widehat{\mm}_{h\tau}^{+}(t)\times\phii(t)}_{\omega}\d t
		\Bigr|
		\\
		&\leq
		\norm[\big]{\phii}_{\bL^\infty(\omega_T)}
		\norm[\big]{\mm_{h\tau}^{+}-\widehat{\mm}_{h\tau}^{+}}_{L^2(0,T;\bH^1(\omega))}
		\norm[\big]{\nabla\widehat{\mm}_{h\tau}^{+}}_{\bL^2(\omega_T)}
		\\
		&\eqrefftwo*{equ:m interpolant bounds}{equ:m predictor H1 consistency}{\lesssim}
		\tau^{1/2}
		\xrightarrow{(h,\tau)\to(0,0)}0.
	\end{align*}
	Moreover, thanks to~\eqref{equ:compactness m weak H1}
	and~$\widehat{\mm}_{h\tau}^+\times\nabla \phii\to \mm\times \nabla\phii$ in $\bL^2(\omega_T)$ by \eqref{equ:compactness m strong}, the remaining contribution converges as
	\begin{align*}
			\int_0^T\scalarproduct{\nabla\mm_{h\tau}^{+}(t)}{\widehat{\mm}_{h\tau}^{+}(t)\times\nabla\phii(t)}_{\omega}\d t
			\xrightarrow{(h,\tau)\to (0,0)} \int_0^T \scalarproduct{\nabla \mm(t)\times \mm(t)}{\nabla \phii(t)}_\omega\d t.
	\end{align*}
	Altogether,
	\begin{align}\label{equ:limit third integral m}
		\const{C}{ex}\int_0^T\scalarproduct{\nabla\mm_{h\tau}^{+}(t)}{\nabla\zetaa_h(t)}_{\omega}\d t
		\xrightarrow{(h,\tau)\to(0,0)}
		\const{C}{ex}\int_0^T\scalarproduct{\nabla\mm(t)\times\mm(t)}{\nabla\phii(t)}_{\omega}\d t .
	\end{align}

	\smallskip
	{\bf Step~4 (Convergence of the fourth, fifth and sixth integrals
	in~\eqref{equ:discrete var form m}--\eqref{equ:limit var form m}).}
	By~\eqref{equ:compactness m strong}, we have
	$\widehat{\mm}_{h\tau}^{+}\to\mm$ in $\bL^2(\omega_T)$, so that
	assumption~\eqref{ass:pi consistency} gives
	$\pp_h(\widehat{\mm}_{h\tau}^{+})\rightharpoonup\pp(\mm)$ in $\bL^2(\omega_T)$.
	Together with the strong convergence of $\zetaa_h$
	in~\eqref{equ:test function m convergence}, this yields
	\begin{align}\label{equ:limit fourth integral m}
		\int_0^T\scalarproduct{\pp_h\bigl(\widehat{\mm}_{h\tau}^{+}(t)\bigr)}{\zetaa_h(t)}_{\omega}\d t
		\xrightarrow{(h,\tau)\to(0,0)}
		\int_0^T\scalarproduct{\pp(\mm(t))\times\mm(t)}{\phii(t)}_{\omega}\d t .
	\end{align}
	Analogously, the weak convergence~\eqref{equ:f discrete convergence} of the
	discrete external fields and the strong convergence of $\zetaa_h$ give
	\begin{align}\label{equ:limit fifth integral m}
		\int_0^T\scalarproduct{\ff_{h\tau}^{+}(t)}{\zetaa_h(t)}_{\omega}\d t
		\xrightarrow{(h,\tau)\to(0,0)}
		\int_0^T\scalarproduct{\ff(t)\times\mm(t)}{\phii(t)}_{\omega}\d t .
	\end{align}
	Finally, we consider the coupling term, which constitutes the new contribution
	with respect to~\cite{AldFeiPra26}. Here, both factors converge strongly: indeed,
	the extrapolated spin accumulation satisfies
	$\widehat{\si}_{h\tau}^{+}\to\si$ in $\bL^2(\Omega_T)$
	by~\eqref{equ:compactness s strong} and hence, upon restriction to $\omega_T$, in
	$\bL^2(\omega_T)$. Together
	with~\eqref{equ:test function m convergence}, this gives
	\begin{align}\label{equ:limit sixth integral m}
		c\int_0^T\scalarproduct{\widehat{\si}_{h\tau}^{+}(t)}{\zetaa_h(t)}_{\omega}\d t
		\xrightarrow{(h,\tau)\to(0,0)}
		c\int_0^T\scalarproduct{\si(t)\times\mm(t)}{\phii(t)}_{\omega}\d t .
	\end{align}
	As already noted in Step~1, the
	convergences~\eqref{equ:limit first integral m}--\eqref{equ:limit sixth integral m}
	prove that the discrete variational formulation~\eqref{equ:discrete var form m}
	tends to its continuous counterpart~\eqref{equ:limit var form m}. This proves Definition~\ref{def:weak solution}(\ref{def: weak solution iii}) and concludes
	the proof.
\end{proof}
\begin{proof}[\bfseries Proof of Theorem~\ref{the:main convergence} (Limit $(\mm,\si)$ satisfies Definition~\ref{def:weak solution}\bf(\ref{def: weak solution iv}))]
	Let $\etaa\in C^1([0,T];\boldsymbol{C}^\infty(\overline{\Omega}))$ be an arbitrary test
	function and define $\etaa_h\colon[0,T]\to\V_h^\Omega$ by
	\begin{align}\label{equ:test function s}
		\etaa_h(t)\coloneqq \III_h[\etaa(t)]
		\qquad\text{for all }t\in[0,T].
	\end{align}
	Thanks to the {elementwise $W^{1,\infty}$ approximation properties on shape-regular meshes} of the
	nodal interpolant, we infer that
	\begin{align}\label{equ:etah convergence}
		\etaa_h\to\etaa
		\quad\text{in }L^\infty\bigl(0,T;\mathbf{W}^{1,\infty}(\Omega)\bigr) .
	\end{align}
	Recalling the conventions $\widehat{\mm}_h^1\coloneqq\mm_h^0$ and
	$\Dtau\si_h^1\coloneqq\dtau\si_h^1$, and noting that
	$\Pi_h(\widehat{\mm}_h^1)=\mm_h^0$ since $\mm_h^0\in\mathcal M_h$, the
	equations~\eqref{equ:first step full discr s} and~\eqref{equ:full discr alg s} of
	Algorithm~\ref{alg: full discr} can be written in the unified form
	\begin{align}\label{equ:pointwise discrete s}
		&\scalarproduct{\Dtau\si_h^{n+1}}{\etaa_h(t)}_{\Omega}
		+
		\fra_h\bigl(\Pi_h(\widehat{\mm}_h^{n+1});\si_h^{n+1},\etaa_h(t)\bigr)
		\nonumber
		\\
		&\qquad
		=
		\beta\scalarproduct{\Pi_h(\widehat{\mm}_h^{n+1})\otimes\jj_h^{n+1}}{\nabla\etaa_h(t)}_{\omega}
		-
		\beta\scalarproduct{\jj_h^{n+1}\cdot\boldsymbol{\nu}}{\Pi_h(\widehat{\mm}_h^{n+1})\cdot\etaa_h(t)}_{\Gamma_J}
	\end{align}
	for all $n=0,\ldots,N-1$ and all $t\in[t_n,t_{n+1})$. We
	integrate~\eqref{equ:pointwise discrete s} over $[t_n,t_{n+1})$ and sum these
	integrals over $n=0,\ldots,N-1$. Recalling the definitions of the interpolants,
	this leads to
	\begin{align}\label{equ:discrete var form s}
		&\int_0^T\scalarproduct{\Dtau\si_{h\tau}^{+}(t)}{\etaa_h(t)}_{\Omega}\d t
		+
		\int_0^T\fra_h\bigl(\Pi_h(\widehat{\mm}_{h\tau}^{+}(t));\si_{h\tau}^{+}(t),\etaa_h(t)\bigr)\d t
		\nonumber
		\\
		&\qquad
		=
		\beta\int_0^T\scalarproduct{\Pi_h(\widehat{\mm}_{h\tau}^{+}(t))\otimes\jj_{h\tau}^{+}(t)}{\nabla\etaa_h(t)}_{\omega}\d t
		-
		\beta\int_0^T\scalarproduct{\jj_{h\tau}^{+}(t)\cdot\boldsymbol{\nu}}{\Pi_h(\widehat{\mm}_{h\tau}^{+}(t))\cdot\etaa_h(t)}_{\Gamma_J}\d t,
	\end{align}
	where we recall that 
	\begin{align*}
	\fra_h(\Pi_h(\widehat{\mm}_{h\tau}^{+});\si_{h\tau}^{+},\etaa_h)
	&\coloneqq
	\scalarproduct{D_0\nabla\si_{h\tau}^+}{\nabla\etaa_h}_{\Omega}
	-
	\beta\beta'
	\scalarproduct{D_0\Pi_h(\widehat{\mm}_{h\tau}^{+})\otimes \left(\nabla\si_{h\tau}^+\cdot\Pi_h(\widehat{\mm}_{h\tau}^{+})\right)}{\nabla\etaa_h}_{\omega}
	\nonumber
	\\
	&\quad
	+
	\scalarproduct{D_0\si_{h\tau}^+}{\etaa_h}_{\Omega}
	+
	\scalarproduct{D_0\left(\si_{h\tau}^+\times\Pi_h(\widehat{\mm}_{h\tau}^{+})\right)}{\etaa_h}_{\omega}.
	\end{align*}
	The following proof considers the convergence of each of these terms above separately as $(h,\tau)\to (0,0)$ and we
	will prove that the limit equation reads
	\begin{align}\label{equ:limit var form s}
		&\int_0^T\dual{\partial_t\si(t)}{\etaa(t)}_{\widetilde{\bH}^{-1}(\Omega),\bH^1(\Omega)}\d t
		+
		\int_0^T\fra\bigl(\mm(t);\si(t),\etaa(t)\bigr)\d t
		\nonumber
		\\
		&\qquad
		=
		\beta\int_0^T\scalarproduct{\mm(t)\otimes\jj(t)}{\nabla\etaa(t)}_{\omega}\d t
		-
		\beta\int_0^T\scalarproduct{\jj(t)\cdot\boldsymbol{\nu}}{\mm(t)\cdot\etaa(t)}_{\Gamma_J}\d t .
	\end{align}
	Notice that $\abs{\mm}=1$ a.e.\ in $\omega_T$ by Proposition~\ref{prop:compactness}
	and $\abs{\Pi_h(\widehat{\mm}_{h\tau}^{+})}\leq1$ a.e.\ in $\omega_T$
	by~\eqref{equ:Pi Linfty}, so that both $\fra(\mm(t);\cdot,\cdot)$
	from~\eqref{equ:a} and $\fra_h(\Pi_h(\widehat{\mm}_{h\tau}^{+}(t));\cdot,\cdot)$
	from~\eqref{equ:discrete spin bilinear form} are well-defined (as well as continuous and
	coercive) for almost every $t\in[0,T]$. 
	
	Before passing to the limit, we collect in a preliminary step some
	convergence properties that will be useful in the following.

	\smallskip
	{\bf Step~1 (Auxiliary convergence results).}
	An application of the triangle inequality and~\eqref{equ:Pi Linfty} implies that
	\begin{align*}
		\norm[\big]{\Pi_h(\widehat{\mm}_{h\tau}^+)\otimes \jj_{h\tau}^+ - \mm\otimes\jj}_{\bL^2(\omega_T)}
		&\le \norm[\big]{\Pi_h(\widehat{\mm}_{h\tau}^+)\otimes (\jj_{h\tau}^+ - \jj)}_{\bL^2(\omega_T)} + \norm[\big]{(\Pi_h(\widehat{\mm}_{h\tau}^+)-\mm)\otimes \jj}_{\bL^2(\omega_T)}\\
		&
		\le
		\norm[\big]{\jj_{h\tau}^+-\jj}_{\bL^2(\omega_T)}
		+
		\norm[\big]{\bigl(\Pi_h(\widehat{\mm}_{h\tau}^+)-\mm\bigr)\otimes \jj}_{\bL^2(\omega_T)} .
	\end{align*}
	The first term vanishes in the limit by~\eqref{equ:j discrete convergence}. For the
	second one, $\abs{\Pi_h(\widehat{\mm}_{h\tau}^+)}\leq1$ by~\eqref{equ:Pi Linfty} and $\abs{\mm}=1$ a.e.\ in $\omega_T$ give
	$\abs{\Pi_h(\widehat{\mm}_{h\tau}^+)-\mm}^2\abs{\jj}^2\leq4\abs{\jj}^2\in L^1(\omega_T)$, since $\jj\in L^2(\omega_T)$,
	while $\Pi_h(\widehat{\mm}_{h\tau}^+)\to\mm$ pointwise a.e.\ in $\omega_T$
	by~\eqref{equ:compactness m pointwise}. This implies that 
	$\abs{\Pi_h(\widehat{\mm}_{h\tau}^+)-\mm}^2\abs{\jj}^2\to0$ pointwise a.e.\ in $\omega_T$.
	The dominated convergence theorem thus
	yields
	$\norm[\big]{(\Pi_h(\widehat{\mm}_{h\tau}^+)-\mm)\otimes\jj}_{\bL^2(\omega_T)}\to0$.
	Altogether, it holds that
	\begin{align}
		\Pi_h(\widehat{\mm}_{h\tau}^+)\otimes\jj_{h\tau}^{+}
		&\to
		\mm\otimes\jj
		\quad \text{in }\bL^2(\omega_T).
		\label{equ:M j convergence}
	\end{align}
	Analogously, the triangle inequality gives
	\begin{align*}
		&\norm[\big]{\Pi_h(\widehat{\mm}_{h\tau}^+)\otimes\Pi_h(\widehat{\mm}_{h\tau}^+)-\mm\otimes\mm}_{\bL^2(\omega_T)}\\
		&\le
		\norm[\big]{\bigl(\Pi_h(\widehat{\mm}_{h\tau}^+)-\mm\bigr)\otimes\Pi_h(\widehat{\mm}_{h\tau}^+)}_{\bL^2(\omega_T)}
		+
		\norm[\big]{\mm\otimes\bigl(\Pi_h(\widehat{\mm}_{h\tau}^+)-\mm\bigr)}_{\bL^2(\omega_T)}
	\end{align*}
	and using $\abs{\Pi_h(\widehat{\mm}_{h\tau}^+)}\leq1$ by~\eqref{equ:Pi Linfty} together
	with $\abs{\mm}=1$ a.e.\ in $\omega_T$, we obtain
	\[
		\norm[\big]{\Pi_h(\widehat{\mm}_{h\tau}^+)\otimes\Pi_h(\widehat{\mm}_{h\tau}^+)-\mm\otimes\mm}_{\bL^2(\omega_T)}
		\leq
		2\norm[\big]{\Pi_h(\widehat{\mm}_{h\tau}^+)-\mm}_{\bL^2(\omega_T)}
		\xrightarrow{(h,\tau)\to(0,0)}0
	\]
	by~\eqref{equ:compactness m strong}, that is,
	$\Pi_h(\widehat{\mm}_{h\tau}^+)\otimes\Pi_h(\widehat{\mm}_{h\tau}^+)\to\mm\otimes\mm$
	in $\bL^2(\omega_T)$.
	This, together with the $\bL^\infty$-bound on $\nabla\etaa_h$ coming from~\eqref{equ:etah convergence}, implies that
	\begin{align}
	\bigl(\Pi_h(\widehat{\mm}_{h\tau}^+)\otimes\Pi_h(\widehat{\mm}_{h\tau}^+)\bigr)\nabla\etaa_h
		&\to
		\bigl(\mm\otimes\mm\bigr)\nabla\etaa
		\quad \text{in }\bL^2(\omega_T).
		\label{equ:M M zeta convergence}
	\end{align}
	Moreover, since $\Gamma_J\subseteq\partial\omega$, the trace operator
	$\bH^1(\omega)\to\bL^2(\Gamma_J)$ is bounded, and so is the induced operator
	$L^2(0,T;\bH^1(\omega))\to L^2(0,T;\bL^2(\Gamma_J))$. Being linear and bounded,
	the latter is weakly continuous, so that~\eqref{equ:compactness m weak H1} yields
	\begin{align}\label{equ:M trace convergence}
		\Pi_h(\widehat{\mm}_{h\tau}^+)\rightharpoonup\mm
		\quad\text{in }L^2\bigl(0,T;\bL^2(\Gamma_J)\bigr).
	\end{align}

	Analogously, the trace inequality~\eqref{equ:trace} applied
	to~\eqref{equ:j discrete convergence}, together with the continuous inclusion
	$\bW^{1,\infty}(\Omega)\subset \boldsymbol{C}(\overline{\Omega})$ applied
	to~\eqref{equ:etah convergence}, give
	\begin{align}\label{equ:trace convergences}
		\jj_{h\tau}^{+}\to\jj
		\quad\text{in }L^2\bigl(0,T;\bL^2(\Gamma_J)\bigr),
		\qquad
		\etaa_h\to\etaa
		\quad\text{in }L^\infty\bigl(0,T;\bL^\infty(\Gamma_J)\bigr).
	\end{align}

	\smallskip
	{\bf Step~2 (Convergence of (\ref{equ:discrete var form s}) to (\ref{equ:limit var form s})).}
We now study the limit of~\eqref{equ:discrete var form s} as $(h,\tau)\to(0,0)$.
	For the first term on the left-hand side, since
	${R_h^*(\Dtau\si_{h\tau}^{+})}\rightharpoonup\partial_t\si$ in
	$L^2(0,T;\widetilde{\bH}^{-1}(\Omega))$ by~\eqref{equ:compactness BDF s} and
	$\etaa_h\to\etaa$ in $L^2(0,T;\bH^1(\Omega))$
	by~\eqref{equ:etah convergence}, {noting that $R_h\etaa_h(t)=\etaa_h(t)$,} it holds that
	\begin{align*}
		\int_0^T\scalarproduct{\Dtau\si_{h\tau}^{+}(t)}{\etaa_h(t)}_{\Omega}\d t
		&\eqreff*{equ:dual reconstruction}{=}
		\int_0^T\dual{{R_h^*(\Dtau\si_{h\tau}^{+}(t))}}{\etaa_h(t)}_{\widetilde{\bH}^{-1}(\Omega),\bH^1(\Omega)}\d t\\
		&\xrightarrow{(h,\tau)\to(0,0)}
		\int_0^T\dual{\partial_t\si(t)}{\etaa(t)}_{\widetilde{\bH}^{-1}(\Omega),\bH^1(\Omega)}\d t .
	\end{align*}
	For the first term on the right-hand side of~\eqref{equ:discrete var form s},
	thanks to $\Pi_h(\widehat{\mm}_{h\tau}^+)\otimes\jj_{h\tau}^{+}
		\to
		\mm\otimes\jj$
	in $\bL^2(\omega_T)$ by~\eqref{equ:M j convergence} and $\nabla\etaa_h\to \nabla\etaa$ in $\bL^2(\Omega_T)$ by~\eqref{equ:etah convergence}, we obtain
	that
	\[
		\beta\int_0^T\scalarproduct{\Pi_h(\widehat{\mm}_{h\tau}^+(t))\otimes\jj_{h\tau}^{+}(t)}{\nabla\etaa_h(t)}_{\omega}\d t
		\xrightarrow{(h,\tau)\to(0,0)}
		\beta\int_0^T\scalarproduct{\mm(t)\otimes\jj(t)}{\nabla\etaa(t)}_{\omega}\d t .
	\]
	For the second term on the right-hand side of~\eqref{equ:discrete var form s},
	since $\Pi_h(\widehat{\mm}_{h\tau}^+)\rightharpoonup\mm$ in $L^2(0,T;\bL^2(\Gamma_J))$
	by~\eqref{equ:M trace convergence}, $\jj_{h\tau}^{+}\to\jj$ in $L^2(0,T;\bL^2(\Gamma_J))$
	and $\etaa_h\to \etaa$ in $L^\infty(0,T;\bL^\infty(\Gamma_J))$ by~\eqref{equ:trace convergences}, it holds that
	\[
		\beta\int_0^T\scalarproduct{\jj_{h\tau}^{+}(t)\cdot\boldsymbol{\nu}}{\Pi_h(\widehat{\mm}_{h\tau}^+(t))\cdot\etaa_h(t)}_{\Gamma_J}\d t
		\xrightarrow{(h,\tau)\to(0,0)}
		\beta\int_0^T\scalarproduct{\jj(t)\cdot\boldsymbol{\nu}}{\mm(t)\cdot\etaa(t)}_{\Gamma_J}\d t .
	\]

	We now analyze the second term on the left-hand side
	of~\eqref{equ:discrete var form s} by analyzing the four terms in the definition of the bilinear forms $\fra_h$ and $\fra$. From the convergences
	$\si_{h\tau}^{+}\rightharpoonup\si$ in $L^2(0,T;\bH^1(\Omega))$
	by~\eqref{equ:compactness s weak H1} and $\nabla\etaa_h\to \nabla\etaa$ in $\bL^2(\Omega_T)$ by~\eqref{equ:etah convergence}, we deduce
	that
	\begin{equation}\label{equ:bilinear form 1}
		\int_0^T\scalarproduct{D_0\nabla\si_{h\tau}^{+}(t)}{\nabla\etaa_h(t)}_{\Omega}\d t
		\xrightarrow{(h,\tau)\to(0,0)}
		\int_0^T\scalarproduct{D_0\nabla\si(t)}{\nabla\etaa(t)}_{\Omega}\d t .
	\end{equation}
	Thanks to $\si_{h\tau}^+\rightharpoonup\si$ in $L^2(0,T;\bH^1(\Omega))$ by~\eqref{equ:compactness s weak H1} and~\eqref{equ:M M zeta convergence},
	we have
	\begin{equation}\label{equ:bilinear form 2}
		\begin{split}
			&\beta\beta'\int_0^T\scalarproduct{D_0\Pi_h(\widehat{\mm}_{h\tau}^+)(t)\otimes\bigl(\nabla\si_{h\tau}^{+}(t)\cdot\Pi_h(\widehat{\mm}_{h\tau}^+(t))\bigr)}{\nabla\etaa_h(t)}_{\omega}\d t\\
		&=
		\beta\beta'\int_0^T\scalarproduct[\big]{D_0\nabla\si_{h\tau}^{+}(t)}{\bigl(\Pi_h(\widehat{\mm}_{h\tau}^+(t))\otimes\Pi_h(\widehat{\mm}_{h\tau}^+(t))\bigr)\nabla\etaa_h(t)}_{\omega}\d t\\
		&\xrightarrow{(h,\tau)\to(0,0)}
		\beta\beta'\int_0^T\scalarproduct{D_0\nabla\si(t)}{\bigl(\mm(t)\otimes\mm(t)\bigr)\nabla\etaa(t)}_{\omega}\d t\\
		&\qquad\qquad=
		\beta\beta'\int_0^T\scalarproduct{D_0\mm(t)\otimes\bigl(\nabla\si(t)\cdot\mm(t)\bigr)}{\nabla\etaa(t)}_{\omega}\d t .
		\end{split}		
	\end{equation}
	From $\si_{h\tau}^{+}\to\si$ in $\bL^2(\Omega_T)$
	by~\eqref{equ:compactness s strong} and $\etaa_h\to\etaa$ in $\bL^2(\Omega_T)$ by~\eqref{equ:etah convergence}, we infer
	that
	\begin{equation}\label{equ:bilinear form 3}
		\int_0^T\scalarproduct{D_0\si_{h\tau}^{+}(t)}{\etaa_h(t)}_{\Omega}\d t
		\xrightarrow{(h,\tau)\to(0,0)}
		\int_0^T\scalarproduct{D_0\si(t)}{\etaa(t)}_{\Omega}\d t .
	\end{equation}
	Finally, we have
	\[
		\int_0^T\scalarproduct{D_0\bigl(\si_{h\tau}^{+}(t)\times\Pi_h(\widehat{\mm}^+_{h\tau}(t))\bigr)}{\etaa_h(t)}_{\omega}\d t
		=
		\int_0^T\scalarproduct{\si_{h\tau}^{+}(t)}{\Pi_h(\widehat{\mm}^+_{h\tau}(t))\times\bigl(D_0\etaa_h(t)\bigr)}_{\omega}\d t .
	\]
	Since $\si_{h\tau}^{+}\to\si$ in $\bL^2(\Omega_T)$
	by~\eqref{equ:compactness s strong}, it suffices to show that the second factor
	converges strongly in $\bL^2(\omega_T)$ as well. To this end, the triangle
	inequality, the bound $D_0\leq D^*$ a.e.\ in $\Omega$, and $\abs{\mm}=1$ a.e.\ in
	$\omega_T$ give
	\begin{align*}
		&\norm[\big]{\Pi_h(\widehat{\mm}^+_{h\tau})\times(D_0\etaa_h)-\mm\times(D_0\etaa)}_{\bL^2(\omega_T)}
		\\
		&\qquad
		\leq
		\norm[\big]{\bigl(\Pi_h(\widehat{\mm}^+_{h\tau})-\mm\bigr)\times(D_0\etaa_h)}_{\bL^2(\omega_T)}
		+
		\norm[\big]{\mm\times\bigl(D_0(\etaa_h-\etaa)\bigr)}_{\bL^2(\omega_T)}
		\\
		&\qquad
		\leq
		D^*\sup_h\norm[\big]{\etaa_h}_{\bL^\infty(\Omega_T)}
		\norm[\big]{\Pi_h(\widehat{\mm}^+_{h\tau})-\mm}_{\bL^2(\omega_T)}
		+
		D^*\norm[\big]{\etaa_h-\etaa}_{\bL^2(\Omega_T)}
		\xrightarrow{(h,\tau)\to(0,0)}0 ,
	\end{align*}
	where we used~\eqref{equ:compactness m strong} for the first term, and the
	convergence $\etaa_h\to \etaa$ in $\bL^2(\Omega_T)$, together with the uniform $\bL^\infty(\Omega_T)$-bound provided
	by~\eqref{equ:etah convergence}. Therefore, we conclude that
	\begin{align}\label{equ:bilinear form 4}
		\int_0^T\scalarproduct{D_0\bigl(\si_{h\tau}^{+}(t)\times\Pi_h(\widehat{\mm}^+_{h\tau}(t))\bigr)}{\etaa_h(t)}_{\omega}\d t
		\xrightarrow{(h,\tau)\to(0,0)}
		\int_0^T\scalarproduct{D_0\bigl(\si(t)\times\mm(t)\bigr)}{\etaa(t)}_{\omega}\d t .
	\end{align}
	Altogether, recalling the definition~\eqref{equ:a} of $\fra$ and that
	$\abs{\mm}=1$ a.e.\ in $\omega_T$, the limits \eqref{equ:bilinear form 1}--\eqref{equ:bilinear form 4} show that 
	\begin{align}\label{equ:limit bilinear form}
		\int_0^T\fra_h(\Pi_h(\widehat{\mm}_{h\tau}^{+}(t));\si_{h\tau}^{+}(t),\etaa_h(t))\d t
		\xrightarrow{(h,\tau)\to(0,0)}
		\int_0^T\fra\bigl(\mm(t);\si(t),\etaa(t)\bigr)\d t .
	\end{align}
	Therefore, the limit of~\eqref{equ:discrete var form s} reads
	\begin{align}
		&\int_0^T\dual{\partial_t\si(t)}{\etaa(t)}_{\widetilde{\bH}^{-1}(\Omega),\bH^1(\Omega)}\d t
		+
		\int_0^T\fra\bigl(\mm(t);\si(t),\etaa(t)\bigr)\d t
		\nonumber
		\\
		&\qquad
		=
		\beta\int_0^T\scalarproduct{\mm(t)\otimes\jj(t)}{\nabla\etaa(t)}_{\omega}\d t
		-
		\beta\int_0^T\scalarproduct{\jj(t)\cdot\boldsymbol{\nu}}{\mm(t)\cdot\etaa(t)}_{\Gamma_J}\d t
	\end{align}
	for all $\etaa\in C^1([0,T];\boldsymbol{C}^\infty(\overline{\Omega}))$ and
	hence, by density, for all $\etaa\in L^2(0,T;\bH^1(\Omega))$. In particular, this shows that
	for almost every $t\in[0,T]$ and $\zetaa\in \bH^1(\Omega)$ it holds that 
	\begin{align*}
			& \dual{\partial_t \si(t)}{\zetaa}_{\widetilde{\bH}^{-1}(\Omega),\bH^1(\Omega)}
			+ \fra(\mm(t);\si(t),\zetaa)
			=
			\beta \scalarproduct{\mm(t)\otimes \jj(t)}{\nabla \zetaa}_{\omega}
			- \beta \scalarproduct{\jj(t) \cdot \boldsymbol{\nu}}{\mm(t)\cdot \zetaa}_{\Gamma_J},
		\end{align*}
	which is exactly Definition~\ref{def:weak solution}(\ref{def: weak solution iv}). 
	\end{proof}
	
Combining the three proofs above with Proposition~\ref{prop:compactness}, the pair
$(\mm,\si)$ is a weak solution of the SDLLG system~\eqref{equ:SDLLG} in the sense of
Definition~\ref{def:weak solution}, which concludes the proof of
Theorem~\ref{the:main convergence}.\qed

\begin{proof}[Proof of Corollary~\ref{cor:full sequence}]
	Suppose that the weak solution $(\mm,\si)$ of the SDLLG system~\eqref{equ:SDLLG}
	is unique, and assume that one of the convergences of
	Theorem~\ref{the:main convergence} fails for the full sequence. Then there exist
	$\varepsilon>0$ and a subsequence which stays at distance at least $\varepsilon$
	from the corresponding limit in the corresponding topology. Since the assumptions
	of Theorem~\ref{the:main convergence} are satisfied for this subsequence as well,
	Proposition~\ref{prop:compactness} and the identification above provide a further
	subsequence which converges to a weak solution of the SDLLG
	system~\eqref{equ:SDLLG}. By uniqueness, this weak solution coincides with
	$(\mm,\si)$, which contradicts the choice of $\varepsilon$.
\end{proof}


\section{Energy estimate}\label{sec:energy}

In this section, we prove that the weak solution $(\mm,\si)$ constructed in
Theorem~\ref{the:main convergence} satisfies a weak counterpart of the energy
law~\eqref{equ:energy law}.

Throughout this section, we impose the following additional assumptions.
\begin{enumerate}[label=\rm (H\arabic*), ref=\rm H\arabic*, start=4]
	\item\label{ass:assumption pi linear self-adjoint} 
	There exists a functional
	$\mathcal E_{\rm low}\colon \bL^2(\omega)\to\R$ that is Lipschitz continuous on bounded subsets of
	$\bL^2(\omega)$ and satisfies
	\begin{equation*}
		\frac{\mathrm d}{\mathrm dt}
		\mathcal E_{\rm low}(\mm(t))
		=
		-\scalarproduct{\pp(\mm(t))}{\partial_t\mm(t)}_\omega
	\end{equation*}
	for almost every $t\in(0,T)$ and every $\mm\in H^1(0,T;\bL^2(\omega))$ satisfying $|\mm|=1$ almost everywhere in $\omega_T$.
	\item\label{ass:assumption pi strong}
	For any sequence $\ww_{h\tau}\to\ww$ strongly in $\bL^2(\omega_T)$,
	it holds that $\pp_h(\ww_{h\tau})\to\pp(\ww)$ strongly in
	$\bL^2(\omega_T)$, which strengthens~\eqref{ass:pi consistency}.
	\item\label{ass:assumption f smooth}
	It holds that
	$\ff\in C^1([0,T];\boldsymbol{C}(\overline{\omega}))$,
	and the discretization
	$\ff_{h\tau}^{+}\in L^2(0,T;\V_h^\omega)$
	is given by
	\[
		\ff_{h\tau}^{+}(t,\zz)
		\coloneqq
		\ff_h^{j+1}(\zz)
		\coloneqq
		\ff(t_{j+1},\zz)
	\]
	for all $j=0,\ldots,N-1$, $t\in[t_j,t_{j+1})$, and $\zz\in\mathcal{N}_h^\omega$.
	Thus, $\ff_{h\tau}^{+}$ is piecewise constant in time and piecewise affine in space.
	In particular,
	\[
		\ff_{h\tau}^{+}\to\ff
		\qquad\text{strongly in }\bL^2(\omega_T)
	\]
	as $(h,\tau)\to(0,0)$, which strengthens
	\eqref{equ:f discrete convergence}.
	\item\label{ass:assumption m0 strong}
	The discrete initial data satisfy $\mm_h^0\to\mm_0$ strongly in
	$\bH^1(\omega)$, which strengthens~\eqref{equ:initial convergence m}.
\end{enumerate}

Assumptions~\eqref{ass:assumption pi linear self-adjoint}--\eqref{ass:assumption pi strong} cover standard uniaxial anisotropy and, with a suitable cutoff extension, nonlinear anisotropies such as cubic anisotropy; see~\cite{AldFeiPraSoen26_part3}. These also accommodate standard $L^2$-stable and strongly consistent hybrid FEM--BEM discretizations of the stray field; see~\cite[Remarks~10 and~23]{AbeHrk_etal14} and~\cite{BruSue_etal14}, as well as~\cite{AldFeiPraSoen26_part3}.
Unlike the energy analysis in~\cite{AbeHrk_etal14}, our argument does not require $\pp$ to be linear and self-adjoint, but instead assumes the existence of an associated energy functional satisfying the chain-rule property. Accordingly, the micromagnetic energy is given by
\begin{align}\label{equ:energy low}
	\mathcal{E}(\mm,\ff)
	\coloneqq
	\frac{\const{C}{ex}}{2}\norm[\big]{\nabla\mm}_{\bL^2(\omega)}^2
	+\mathcal E_{\rm low}(\mm)
	-\scalarproduct{\ff}{\mm}_{\omega}.
\end{align}
In particular, if $\pp$ is linear and self-adjoint,  then~\eqref{ass:assumption pi linear self-adjoint} holds with $\mathcal E_{\rm low}(\mm)=-\frac12\scalarproduct{\pp(\mm)}{\mm}_{\omega}$, and~\eqref{equ:energy low} is the energy~\eqref{equ:energy}.

We first rewrite the identity derived in the proof of
Lemma~\ref{lem:stability magnetization}, retaining the term $\mathcal{G}_h^n$
from~\eqref{equ:G definition}. Since certain parts of the proof are analogous
to~\cite[Lemma~3.5]{AldFeiPra26} and
\cite[Lemma~7]{AldFeiPraSoen26_part3}, we only sketch the main steps for brevity and
refer to these references for the details.

\begin{lemma}[Discrete energy identity]\label{lem:discrete energy}
	Recall $\mathcal{G}_h^n$ from~\eqref{equ:G definition} and $\boldsymbol{F}_h^{n+1}$, $n=0,\ldots,N-1$, from~\eqref{equ:discrete F}. Moreover, define
	\begin{align}\label{equ:eta}
		\eta_0\coloneqq
		\norm[\big]{\nabla\mm_h^1}_{\bL^2(\omega)}^2
		-\norm[\big]{\nabla\mm_h^0}_{\bL^2(\omega)}^2.
	\end{align}
	Then, for all $1\leq n\leq N$, Algorithm~\ref{alg: full discr} guarantees
	that
	\begin{align}\label{equ:discrete energy identity}
		&\const{C}{ex}\mathcal{G}_h^n
		+\alpha\tau\sum_{j=0}^{n-1}
		\norm[\big]{\vv_h^j}_{\bL^2(\omega)}^2
		+\frac{\const{C}{ex}\tau^4}{4}\sum_{j=1}^{n-1}
		\norm[\big]{\nabla\dtau^2\mm_h^{j+1}}_{\bL^2(\omega)}^2
		\nonumber\\
		&\qquad=
		\frac{\const{C}{ex}}2
		\norm[\big]{\nabla\mm_h^0}_{\bL^2(\omega)}^2
		+\tau\sum_{j=0}^{n-1}
		\scalarproduct{\boldsymbol{F}_h^{j+1}}{\vv_h^j}_{\omega}
		+\frac{\const{C}{ex}}4\eta_0.
	\end{align}
\end{lemma}

\begin{proof}
	Testing~\eqref{equ:first step full discr m} with
	$\phii_h=\tau\vv_h^0$ and using~\eqref{equ:first step m1} gives
	\[
		\alpha\tau\norm[\big]{\vv_h^0}_{\bL^2(\omega)}^2
		+\frac{\const{C}{ex}}2\norm[\big]{\nabla\mm_h^1}_{\bL^2(\omega)}^2
		+\frac{\const{C}{ex}\tau^2}{2}
		\norm[\big]{\nabla\vv_h^0}_{\bL^2(\omega)}^2
		=
		\frac{\const{C}{ex}}2\norm[\big]{\nabla\mm_h^0}_{\bL^2(\omega)}^2
		+\tau\scalarproduct{\boldsymbol{F}_h^1}{\vv_h^0}_{\omega}.
	\]
	Since~\eqref{equ:G definition} yields
	\[
		\mathcal{G}_h^1
		=\frac12\norm[\big]{\nabla\mm_h^1}_{\bL^2(\omega)}^2
		+\frac{\tau^2}{2}\norm[\big]{\nabla\vv_h^0}_{\bL^2(\omega)}^2
		+\frac14\eta_0,
	\]
	this proves~\eqref{equ:discrete energy identity} for $n=1$.
	For $n\geq2$, we sum the identity
	\[
		\const{C}{ex}(\mathcal{G}_h^{j+1}-\mathcal{G}_h^j)
		+\alpha\tau\norm[\big]{\vv_h^j}_{\bL^2(\omega)}^2
		+\frac{\const{C}{ex}\tau^4}{4}
		\norm[\big]{\nabla\dtau^2\mm_h^{j+1}}_{\bL^2(\omega)}^2
		=\tau\scalarproduct{\boldsymbol{F}_h^{j+1}}{\vv_h^j}_{\omega},
	\]
	derived in the proof of Lemma~\ref{lem:stability magnetization}, over
	$j=1,\ldots,n-1$ and add the identity for $n=1$.
	This concludes the proof.
\end{proof}

\begin{theorem}[Energy inequality]\label{the:energy}
	Under the assumptions of Theorem~\ref{the:main convergence}
	and~\eqref{ass:assumption pi linear self-adjoint}--\eqref{ass:assumption m0 strong},
	the weak solution $(\mm,\si)$ of the SDLLG system~\eqref{equ:SDLLG}
	obtained as the limit in Theorem~\ref{the:main convergence} satisfies,
	for almost every $t'\in(0,T)$, the energy inequality with $\mathcal E$ from~\eqref{equ:energy low}
	\begin{align}\label{equ:energy inequality}
		\mathcal{E}\bigl(\mm(t'),\ff(t')\bigr)
		&+\alpha\int_0^{t'}
		\norm[\big]{\partial_t\mm(t)}_{\bL^2(\omega)}^2\d t
		+\int_0^{t'}
		\scalarproduct{\partial_t\ff(t)}{\mm(t)}_{\omega}\d t
		\nonumber\\
		&\qquad\leq
		\mathcal{E}\bigl(\mm_0,\ff(0)\bigr)
		+c\int_0^{t'}
		\scalarproduct{\si(t)}{\partial_t\mm(t)}_{\omega}\d t.
	\end{align}
\end{theorem}

\begin{proof}
	The proof is split into four steps and follows
	\cite[Theorem~2.2]{AldFeiPra26},
	\cite[Theorem~3]{AldFeiPraSoen26_part3}, and
	\cite[Theorem~24]{AbeHrk_etal14}. We sketch the arguments for the
	lower-order and applied fields, and detail the spin contribution and the
	passage to the limit in $\mathcal{G}_h^n$.

	\smallskip
	{\bf Step~1 (Discrete energy inequality).}
	In analogy with~\eqref{equ:time reconstructions m}, we define
	$\ff_\tau$ and $\ff_\tau^{\pm}$ as piecewise linear and piecewise constant
	interpolations of $\ff^n\coloneqq\ff(t_n)$.
	Since
	$\ff\in C^1([0,T];\bL^2(\omega))$ by~\eqref{ass:assumption f smooth}, it holds that
	\begin{align}\label{equ:energy f reconstructions}
	\ff_\tau,\ff_\tau^{\pm}\to\ff
	\quad\text{ and } \quad \partial_t\ff_\tau\to\partial_t\ff\quad \text{ in } \quad\bL^2(\omega_T) \quad\text{ as } \quad\tau\to0.
	\end{align}

	Now, for $n=1,\ldots,N$, define the following discrete modified energy, where $\const{C}{ex}\mathcal{G}_h^n$ replaces the exchange contribution:
	\begin{align}\label{equ:energy discrete modified energy}
		\mathcal{E}_h^n
		\coloneqq\const{C}{ex}\mathcal{G}_h^n
		+\mathcal E_{\rm low}(\mm_h^n)
		-\scalarproduct{\ff^n}{\mm_h^n}_{\omega}.
	\end{align}
	Recalling $\boldsymbol{F}_h^{n+1}=\pp_h(\widehat{\mm}_h^{n+1})+\ff_h^{n+1}+c\,\widehat{\si}_h^{n+1}$
	and the definition~\eqref{equ:energy low} of the energy,~\eqref{equ:discrete energy identity} from
	Lemma~\ref{lem:discrete energy} yields, after dropping the non-negative second-difference term on its left-hand side,
	\begin{align}
		\mathcal{E}_h^n-\mathcal{E}(\mm_h^0,\ff^0)
		&=
		\const{C}{ex}\mathcal{G}_h^n
		-\frac{\const{C}{ex}}{2}
		\norm[\big]{\nabla\mm_h^0}_{\bL^2(\omega)}^2
		+\mathcal E_{\rm low}(\mm_h^n)
		-\mathcal E_{\rm low}(\mm_h^0)
		\nonumber\\
		&\qquad
		-\scalarproduct{\ff^n}{\mm_h^n}_{\omega}
		+\scalarproduct{\ff^0}{\mm_h^0}_{\omega}
		\nonumber\\
		&\leq
		-\alpha\tau\sum_{j=0}^{n-1}
		\norm[\big]{\vv_h^j}_{\bL^2(\omega)}^2
		+T_{\pp}^n+T_{\ff}^n
		+c\tau\sum_{j=0}^{n-1}
		\scalarproduct{\widehat{\si}_h^{j+1}}{\vv_h^j}_{\omega}
		+\frac{\const{C}{ex}}4\eta_0,
		\label{eq:ineq En min E0}
	\end{align}
	where
	\begin{align*}
		T_{\pp}^n
		&\coloneqq\sum_{j=0}^{n-1}\Big[
		\tau\scalarproduct{\pp_h(\widehat{\mm}_h^{j+1})}{\vv_h^j}_{\omega}
		+\mathcal E_{\rm low}(\mm_h^{j+1})
		-\mathcal E_{\rm low}(\mm_h^j)
		\Big],
		\\
		T_{\ff}^n
		&\coloneqq\sum_{j=0}^{n-1}\Big[
		\tau\scalarproduct{\ff_h^{j+1}}{\vv_h^j}_{\omega}
		-\scalarproduct{\ff^{j+1}}{\mm_h^{j+1}}_{\omega}
		+\scalarproduct{\ff^j}{\mm_h^j}_{\omega}
		\Big].
	\end{align*}
	For every $t'\in[0,T)$, let $n=n(t')$ be the unique index satisfying
	$t'\in[t_{n-1},t_n)$. Since $t'\leq t_n$, the preceding inequality \eqref{eq:ineq En min E0} can be written in terms of the time reconstructions as
	\begin{align}\label{equ:discrete energy inequality}
		&\mathcal{E}_h^{n(t')}-\mathcal{E}(\mm_{h\tau}^-(0),\ff(0))
		+\alpha\int_0^{t'}
		\norm[\big]{\vv_{h\tau}^{-}(t)}_{\bL^2(\omega)}^2\d t
		\nonumber\\
		&\qquad\leq
		T_{\pp}^{n(t')}+T_{\ff}^{n(t')}
		+c\int_0^{t_{n(t')}}
		\scalarproduct{\widehat{\si}_{h\tau}^{+}(t)}{\vv_{h\tau}^{-}(t)}_{\omega}\d t
		+\frac{\const{C}{ex}}4\eta_0.
	\end{align}

	\smallskip
	{\bf Step~2 (Dissipation and work terms).}
	By~\eqref{equ:compactness v}, we have
	$\vv_{h\tau}^{-}\rightharpoonup\partial_t\mm$
	in $\bL^2(\omega_T)$. Weak lower semicontinuity 
	therefore gives
	\begin{align*}
		&\int_0^{t'}
		\norm[\big]{\partial_t\mm(t)}_{\bL^2(\omega)}^2\d t
		\leq
		\liminf_{(h,\tau)\to(0,0)}
		\int_0^{t'}
		\norm[\big]{\vv_{h\tau}^{-}(t)}_{\bL^2(\omega)}^2\d t.
	\end{align*}
	
	The treatment of $T_{\pp}^n$ and $T_{\ff}^n$ follows
	Step~2 of the proof of~\cite[Theorem~3]{AldFeiPraSoen26_part3}.
	{Note that the relevant arguments therein do not require quasi-uniformity of the triangulations, but rely only on the uniform bounds and convergence properties already established above, together with the chain rule for $\mathcal E_{\rm low}$ and integration by parts in time.}
	Thanks to the available convergences in \eqref{equ:compactness m strong} and~\eqref{equ:compactness v}, together with assumptions~\eqref{ass:assumption pi linear self-adjoint}--\eqref{ass:assumption f smooth}, and~\eqref{equ:energy f reconstructions}, we obtain,
	for every measurable $I\subseteq(0,T)$,
	\begin{align}\label{equ:energy field remainders}
		\int_I T_{\pp}^{n(t')}\d t'
		\to 0,
		\qquad
		\int_I T_{\ff}^{n(t')}\d t'
		\to
		-\int_I\int_0^{t'}
		\scalarproduct{\partial_t\ff(t)}{\mm(t)}_{\omega}
		\d t\,\d t'.
	\end{align}

	For the spin coupling term, we proceed as follows.
	The extrapolated spin accumulation satisfies $\widehat{\si}_{h\tau}^{+}\to\si$ strongly in $\bL^2(\Omega_T)$ by~\eqref{equ:compactness s strong} and hence, upon restriction, in
	$\bL^2(\omega_T)$. Moreover, $\vv_{h\tau}^{-}\rightharpoonup\partial_t\mm$ in $\bL^2(\omega_T)$
	by~\eqref{equ:compactness v}. Thus, their scalar product integrated over
	$[0,t']$ converges to the corresponding limit for every fixed $t'\in[0,T]$.
	The remaining interval is controlled by
	\begin{align*}
		\left|\int_{t'}^{t_{n(t')}}
		\scalarproduct{\widehat{\si}_{h\tau}^{+}(t)}{\vv_{h\tau}^{-}(t)}_{\omega}\d t
		\right|
		&\leq\tau^{1/2}
		\norm[\big]{\widehat{\si}_{h\tau}^{+}}_{L^\infty(0,T;\bL^2(\Omega))}
		\norm[\big]{\vv_{h\tau}^{-}}_{\bL^2(\omega_T)}
		\eqrefftwo*{equ:s interpolant bounds}{equ:m interpolant bounds}{\lesssim} \tau^{1/2}\to0.
	\end{align*}
	Consequently,
	\begin{align}\label{equ:Ts limit}
		c\int_0^{t_{n(t')}}
		\scalarproduct{\widehat{\si}_{h\tau}^{+}(t)}{\vv_{h\tau}^{-}(t)}_{\omega}\d t
		\xrightarrow{(h,\tau)\to(0,0)}
		c\int_0^{t'}
		\scalarproduct{\si(t)}{\partial_t\mm(t)}_{\omega}\d t.
	\end{align}
	Since the left-hand side of~\eqref{equ:Ts limit} is bounded uniformly in $t'$ by~\eqref{equ:s interpolant bounds} and~\eqref{equ:m interpolant bounds}, the dominated convergence theorem shows that~\eqref{equ:Ts limit} also holds after integration over~$I$.
	
	\smallskip
	{\bf Step~3 ({Energy terms and initial-step correction}).}
	First, $\mm_h^0\to\mm_0$ strongly in $\bH^1(\omega)$
	by~\eqref{ass:assumption m0 strong}, and hence
	\[
		\mathcal{E}(\mm_h^0,\ff(0))
		\to\mathcal{E}(\mm_0,\ff(0)).
	\]
	Moreover, the initial-step estimate~\eqref{equ:initial m stability} gives
	$\const{C}{ex}\eta_0\leq\tau\norm[\big]{\boldsymbol{F}_h^1}_{\bL^2(\omega)}^2/\alpha$.
	The triangle inequality, together with~\eqref{ass:pi bounded},
	\eqref{equ:initial convergence m},
	\eqref{equ:initial convergence s},
	and~\eqref{ass:assumption f smooth}, gives
	\begin{align*}
		\tau\norm[\big]{\boldsymbol{F}_h^1}_{\bL^2(\omega)}^2
		&\leq
		3\tau\const{C}{$\pi$}^2
		\norm[\big]{\mm_h^0}_{\bL^2(\omega)}^2
		+3\tau\norm[\big]{\ff_h^1}_{\bL^2(\omega)}^2
		+3c^2\tau\norm[\big]{\si_h^0}_{\bL^2(\Omega)}^2
		\lesssim \tau\to0.
	\end{align*}
	In particular,
	$
		\limsup_{(h,\tau)\to(0,0)}\eta_0\leq0.
	$
	Let $I\subseteq[0,T]$ be a measurable set. The definition~\eqref{equ:G definition}
	and the time reconstructions give
	\begin{align*}
		\mathcal{G}_h^{n(t')}
		=\frac14\Bigl[
		\norm[\big]{\nabla\mm_{h\tau}^{+}(t')}_{\bL^2(\omega)}^2
		+\norm[\big]{2\nabla\mm_{h\tau}^{+}(t')-
		\nabla\mm_{h\tau}^{-}(t')}_{\bL^2(\omega)}^2
		\Bigr].
	\end{align*}
	By~\eqref{equ:compactness m weak H1}, both
	$\nabla\mm_{h\tau}^{+}$ and
	$2\nabla\mm_{h\tau}^{+}-\nabla\mm_{h\tau}^{-}$ converge weakly to
	$\nabla\mm$ in $\bL^2(\omega_T)$. Weak lower semicontinuity therefore yields
	\begin{align}\label{equ:energy G liminf}
		\frac12\int_I\norm[\big]{\nabla\mm(t')}_{\bL^2(\omega)}^2\d t'
		\leq\liminf_{(h,\tau)\to(0,0)}
		\int_I\mathcal{G}_h^{n(t')}\d t'.
	\end{align}
	{The contribution of $\mathcal E_{\rm low}$ to~\eqref{equ:energy discrete modified energy} converges thanks to the Lipschitz continuity of $\mathcal E_{\rm low}$ on bounded subsets of $\bL^2(\omega)$ from~\eqref{ass:assumption pi linear self-adjoint}, together with the uniform bounds~\eqref{equ:m interpolant bounds} and the strong convergence~\eqref{equ:compactness m strong}, as
	\begin{align*}
		\int_I
		\mathcal E_{\rm low}\bigl(\mm_{h\tau}^{+}(t')\bigr)\d t'
		\xrightarrow{(h,\tau)\to(0,0)}
		\int_I
		\mathcal E_{\rm low}\bigl(\mm(t')\bigr)\d t'.
	\end{align*}}
	Similarly, the contribution of $\ff$ converges thanks to $\ff_\tau^{+}\to\ff$ in $\bL^2(\omega_T)$
	by~\eqref{equ:energy f reconstructions} and
	$\mm_{h\tau}^{+}\to\mm$ in $\bL^2(\omega_T)$ as
	\begin{align*}
		&\int_I
		\scalarproduct{\ff_\tau^{+}(t')}
		{\mm_{h\tau}^{+}(t')}_{\omega}\d t'
		\xrightarrow{(h,\tau)\to(0,0)}
		\int_I
		\scalarproduct{\ff(t')}{\mm(t')}_{\omega}\d t'.
	\end{align*}
	These convergences,
	together with~\eqref{equ:energy G liminf}, yield
	\begin{align}\label{equ:energy liminf}
		\int_I
		\mathcal{E}\bigl(\mm(t'),\ff(t')\bigr)\d t'
		\leq\liminf_{(h,\tau)\to(0,0)}
		\int_I\mathcal{E}_h^{n(t')}\d t'.
	\end{align}

	\smallskip
	{\bf Step~4 (Continuous energy inequality).}
	Integrating~\eqref{equ:discrete energy inequality} over $I$,
	applying Fatou's lemma to the dissipation term, and using
	Steps~2--3, we obtain
		\begin{align*}
		&\int_I\Bigl[
		\mathcal{E}\bigl(\mm(t'),\ff(t')\bigr)
		+\alpha\int_0^{t'}
		\norm[\big]{\partial_t\mm(t)}_{\bL^2(\omega)}^2\d t
		\Bigr]\d t'
		\\
		&\quad\leq\liminf_{(h,\tau)\to(0,0)}
		\int_I\Bigl[
		\mathcal{E}_h^{n(t')}
		+\alpha\int_0^{t'}
		\norm[\big]{\vv_{h\tau}^{-}(t)}_{\bL^2(\omega)}^2\d t
		\Bigr]\d t'
		\\
		&\quad\eqreff*{equ:discrete energy inequality}{\leq}\limsup_{(h,\tau)\to(0,0)}
		\int_I\Bigl[
		\mathcal{E}\bigl(\mm_{h\tau}^-(0),\ff(0)\bigr)
		+T_{\pp}^{n(t')}+T_{\ff}^{n(t')}
		+c\int_0^{t_{n(t')}}
		\scalarproduct{\widehat{\si}_{h\tau}^{+}(t)}
		{\vv_{h\tau}^{-}(t)}_{\omega}\d t
		+\frac{\const{C}{ex}}4\eta_0
		\Bigr]\d t'
		\\
		&\quad\eqrefftwo*{equ:energy field remainders}{equ:Ts limit}{=}
		\quad\int_I\Bigl[
		\mathcal{E}\bigl(\mm_0,\ff(0)\bigr)
		-\int_0^{t'}
		\scalarproduct{\partial_t\ff(t)}{\mm(t)}_{\omega}\d t
		+c\int_0^{t'}
		\scalarproduct{\si(t)}{\partial_t\mm(t)}_{\omega}\d t
		\Bigr]\d t'
		\\[0.8ex]
		&\qquad\;
		+\frac{\const{C}{ex}}4|I|
		\limsup_{(h,\tau)\to(0,0)}\eta_0
		\\
		&\quad\leq
		\int_I\Bigl[
		\mathcal{E}\bigl(\mm_0,\ff(0)\bigr)
		-\int_0^{t'}
		\scalarproduct{\partial_t\ff(t)}{\mm(t)}_{\omega}\d t
		+c\int_0^{t'}
		\scalarproduct{\si(t)}{\partial_t\mm(t)}_{\omega}\d t
		\Bigr]\d t'.
	\end{align*}
	Since $I$ is arbitrary, the integrands satisfy the same inequality
	for almost every $t'\in(0,T)$, which is
	exactly~\eqref{equ:energy inequality} and concludes the proof.
\end{proof}

\begin{remark}[Time regularity of the applied field]
    The conclusions of Theorem~\ref{the:energy} remain valid
    if~\eqref{ass:assumption f smooth} is replaced by
    \[
        \ff\in H^1(0,T;\bL^2(\omega)),
        \qquad
        \ff_h^n=P_h^\omega\ff(t_n)
        \quad\text{for }n=1,\ldots,N,
    \]
    where $P_h^\omega$ is the $\bL^2$-orthogonal projection
    from~\eqref{equ:L2 projection};
    see also~\cite[Section~5]{AbeHrk_etal14}
    and~\cite[Section~4.1]{Rug16} for the weaker time regularity.
    Thus, the nodal interpolation in~\eqref{ass:assumption f smooth}
    is replaced by an approximation which is well-defined for fields which are $\bL^2$ in space. 
    The stronger $C^1$-assumption in time is adopted for simplicity, as in~\cite[Section~4.3]{Rug16}.

    In the proof, the time regularity of $\ff$ enters through the convergence of its time reconstructions in Step~1 and the estimates for $T_{\ff}^n$ in Step~2. The uniform $L^\infty(0,T;\bL^2(\omega))$ bound of $\ff_{h\tau}^+$ follows also in this case from $\bL^2$-contractivity of $P_h^\omega$ and $\ff\in C([0,T];\bL^2(\omega))$.
    Under the weaker assumption, the continuous embedding
    \[
        H^1(0,T;\bL^2(\omega))
        \hookrightarrow C^0([0,T];\bL^2(\omega))
    \]
    provides well-defined values $\ff(t_n)$ and $\ff(0)$.
    Moreover, the piecewise linear interpolant of these values satisfies
    \[
        \partial_t\ff_\tau(t)
        =\frac1\tau\int_{t_j}^{t_{j+1}}
        \partial_t\ff(r)\d r
        \qquad\text{for a.e. }t\in[t_j,t_{j+1}),
    \]
    and therefore
    \[
        \norm[\big]{\partial_t\ff_\tau}_{\bL^2(\omega_T)}
        \leq\norm[\big]{\partial_t\ff}_{\bL^2(\omega_T)},
        \qquad
        \partial_t\ff_\tau\to\partial_t\ff
        \quad\text{in }\bL^2(\omega_T).
    \]
    Together with the continuity of $\ff$, this proves that
    the convergences in Step~1 remain valid.
    For the initial-step argument in Step~3, the stability
    of $P_h^\omega$ gives
    \[
        \norm[\big]{\ff_h^1}_{\bL^2(\omega)}
        \leq\norm[\big]{\ff(t_1)}_{\bL^2(\omega)}
        \leq\norm[\big]{\ff}_{C^0([0,T];\bL^2(\omega))}.
    \]
    Hence, the estimate
    $\tau\norm[\big]{\boldsymbol{F}_h^1}_{\bL^2(\omega)}^2
    \lesssim \tau\to0$
    remains valid.

    Finally, $\ff_{h\tau}^{+}=P_h^\omega\ff_\tau^{+}$,
    and hence the stability of $P_h^\omega$ gives
    \begin{align*}
        \norm[\big]{\ff_{h\tau}^{+}-\ff}_{\bL^2(\omega_T)}
        &\leq
		\norm[\big]{P_h^\omega\big(\ff_\tau^{+}-\ff\big)}_{\bL^2(\omega_T)}
        +\norm[\big]{P_h^\omega\ff-\ff}_{\bL^2(\omega_T)}\\
		&
		\leq
        \norm[\big]{\ff_\tau^{+}-\ff}_{\bL^2(\omega_T)}
        +\norm[\big]{P_h^\omega\ff-\ff}_{\bL^2(\omega_T)}
        \to0.
	\end{align*}
    Thus, the strong convergence of the discrete applied field
    is also preserved.
\end{remark}


\section{Numerical experiments}\label{sec:num-exp}

In this section, we present several numerical experiments to support our theoretical results using the open-source computing platform \textsc{FEniCS}~\cite{AlnaesEtal15}.

\subsection{Manufactured solution convergence test}
\label{subsec:manuf}

This subsection investigates the convergence properties of
Algorithm~\ref{alg: full discr} by means of a manufactured-solution test.
We set $\pp_h\equiv0$ and augment the spin-diffusion equation by an
additional source term $\boldsymbol{g}$. Throughout the experiment, we consider
\[
\Omega=(0,1)^3
\qquad\text{and}\qquad
\omega=(0,1)^2\times \Bigl(\frac13 , \frac23\Bigr),
\]
so that the ferromagnetic region $\omega$ is a slab strictly contained
in the conductor $\Omega$.

We use the parameters
\[
\alpha=0.5,
\qquad
\const{C}{ex}=0.1,
\qquad
c=0.2,
\qquad
D_0=0.05,
\qquad
\beta=0.4,
\qquad
\beta'=0.5,
\]
and take the final time $T=0.01$.

For $\chi\in\{0,1\}$, define
\begin{align}
	\vartheta_{\chi}(\xx,t)
	&\coloneqq
	0.35
	\left(
	1+0.30\chi
	\cos(\pi x_1)\cos(\pi x_2)
	\right)
	\left(
	1+0.20\sin t
	\right).
	\label{equ:MMS theta}
\end{align}
The manufactured magnetization is
\begin{align}
	\mm(\xx,t)
	&\coloneqq
	\begin{pmatrix}
		\cos\vartheta_{\chi}(\xx,t)\\
		\sin\vartheta_{\chi}(\xx,t)\\
		0
	\end{pmatrix},
	\label{equ:MMS m}
\end{align}
so that $|\mm|=1$ pointwise. The manufactured spin accumulation is
\begin{align}
	\si(\xx,t)
	\coloneqq
	\begin{pmatrix}
		0.20
		\left(
		1+0.25\chi\cos(\pi x_1)\cos(\pi x_2)
		\right)
		\cos(2t)
		\\[1mm]
		0.15
		\left(
		1+0.20\chi\cos(2\pi x_1)\cos(\pi x_2)
		\right)
		\sin(1.5t)
		\\[1mm]
		0.10
		\left(
		1+0.15\chi\cos(\pi x_1)\cos(2\pi x_2)
		\right)
		\cos t
	\end{pmatrix}.
	\label{equ:MMS s}
\end{align}
Both fields are independent of $x_3$. {The parameter $\chi\in\{0,1\}$ switches the spatial variation of the manufactured solution on ($\chi=1$) or off ($\chi=0$): for $\chi=0$, the fields $\mm$ and $\si$ are spatially constant and hence represented exactly by $\P_1$ finite elements, so that the error is dominated by the temporal discretization; see Section~\ref{subsec:numerical temporal convergence}.} Moreover, the trigonometric
factors are chosen such that the corresponding normal derivatives
vanish on the relevant planar boundaries.

The prescribed current density is
\begin{align}
	\jj(\xx,t)
	&\coloneqq
	\bigl(0.70+0.10\cos t\bigr)
	\begin{pmatrix}
		\sin(\pi x_1)\cos(\pi x_2)
		\\
		-\cos(\pi x_1)\sin(\pi x_2)
		\\
		0
	\end{pmatrix}.
	\nonumber
\end{align}

The external field $\ff$ in the LLG equation is chosen such that
\eqref{equ:MMS m}--\eqref{equ:MMS s} satisfy \eqref{equ:SDLLG m}.
In~\eqref{equ:SDLLG s}, we add a source term
$\boldsymbol g$ to the right-hand side to ensure \eqref{equ:MMS m}--\eqref{equ:MMS s} also satisfy the augmented spin-diffusion equation. {Explicitly, $\boldsymbol g\coloneqq\partial_t\si+\nabla\cdot\js+D_0\si+D_0(\si\times\mm)$ is computed from~\eqref{equ:MMS m}--\eqref{equ:MMS s}, and the corresponding term $\scalarproduct{\boldsymbol g(t_{n+1})}{\zetaa_h}_{\Omega}$ is added to the right-hand sides of~\eqref{equ:first step full discr s} and~\eqref{equ:full discr alg s}.} The initial data are given by the exact solution at $t=0$,
\[
\mm_0=\mm(\cdot,0),
\qquad
\si_0=\si(\cdot,0),
\]
and the corresponding discrete initial data are their nodal interpolants.

We discretize $\Omega$ by a uniform tetrahedral mesh obtained from an
${M\times M\times M}$ subdivision of the unit cube and construct the
mesh of $\omega$ as the corresponding fitted submesh. Since all values
of {$M$} used below are divisible by $3$, the interfaces
$\{x_3=\frac13\}$ and $\{x_3=\frac23\}$ coincide exactly with mesh planes. For these
meshes, $h=\dfrac{\sqrt{3}}{{M}}$.

For the errors, we write
\begin{align*}
	e_{\mm,0}
	&\coloneqq
	\norm{\mm(T)-\mm_h^N}_{\bL^2(\omega)},
	&
	e_{\mm,1}
	&\coloneqq
	\norm{\nabla(\mm(T)-\mm_h^N)}_{\bL^2(\omega)},
	\\
	e_{\si,0}
	&\coloneqq
	\norm{\si(T)-\si_h^N}_{\bL^2(\Omega)},
	&
	e_{\si,1}
	&\coloneqq
	\norm{\nabla(\si(T)-\si_h^N)}_{\bL^2(\Omega)}.
\end{align*}
The experimental order of convergence between two consecutive
discretizations is calculated in the standard manner.

\subsubsection{Spatial convergence}
\label{subsec:numerical spatial convergence}

For the spatial convergence test, we set $\chi=1$ in~\eqref{equ:MMS theta}--\eqref{equ:MMS s} 
and fix $\tau=10^{-3}$, while varying ${M}=3,6,9,12,15$.
Table~\ref{tab:spatial convergence MMS} reports the resulting errors.
The observed rates approach second order in the $\bL^2$-norm and first
order in the $\bH^1$-seminorm for both $\mm_h$ and $\si_h$, in agreement
with the expected spatial accuracy of piecewise affine finite
elements.

\begin{table}[!ht]
	\centering
	\small
	\setlength{\tabcolsep}{5pt}
		\begin{tabular}{
				c
				@{\qquad}
				cc 
				@{\qquad}
				cc 
				@{\qquad}
				cc 
				@{\qquad}
				cc 
			}
			\toprule
			$M$
			& $e_{\mm,0}$ & rate
			& $e_{\mm,1}$ & rate
			& $e_{\si,0}$ & rate
			& $e_{\si,1}$ & rate
			\\
			\midrule
			3
			& $6.236{\times}10^{-3}$ & --
			& $6.602{\times}10^{-2}$ & --
			& $6.139{\times}10^{-3}$ & --
			& $6.513{\times}10^{-2}$ & --
			\\
			6
			& $1.670{\times}10^{-3}$ & $1.90$
			& $3.460{\times}10^{-2}$ & $0.93$
			& $1.682{\times}10^{-3}$ & $1.87$
			& $3.460{\times}10^{-2}$ & $0.91$
			\\
			9
			& $7.529{\times}10^{-4}$ & $1.97$
			& $2.330{\times}10^{-2}$ & $0.97$
			& $7.606{\times}10^{-4}$ & $1.96$
			& $2.338{\times}10^{-2}$ & $0.97$
			\\
			12
			& $4.258{\times}10^{-4}$ & $1.98$
			& $1.755{\times}10^{-2}$ & $0.99$
			& $4.306{\times}10^{-4}$ & $1.98$
			& $1.762{\times}10^{-2}$ & $0.98$
			\\
			15
			& $2.732{\times}10^{-4}$ & $1.99$
			& $1.407{\times}10^{-2}$ & $0.99$
			& $2.764{\times}10^{-4}$ & $1.99$
			& $1.414{\times}10^{-2}$ & $0.99$
			\\
			\bottomrule
	\end{tabular}
	\caption{
		Spatial convergence for $\tau=10^{-3}$ and $T=0.01$.
	}
	\label{tab:spatial convergence MMS}
\end{table}

\subsubsection{Temporal convergence and nodal unit-length defect}
\label{subsec:numerical temporal convergence}

To isolate the temporal discretization error, we set $\chi=0$ in~\eqref{equ:MMS theta}--\eqref{equ:MMS s}.
Consequently, $\mm$ and $\si$ remain time dependent but are spatially
constant, and thus are represented exactly in space by the $\P_1$
finite element. We fix $M=12$ and consider $N=2,4,8,16,24$, where $N= T/\tau$.
Table~\ref{tab:temporal convergence MMS} shows second-order convergence
for both $\mm_h$ and $\si_h$ in the $\bL^2$-norm.

\begin{figure}[!htb]
	\centering
	\includegraphics[width=\textwidth]
	{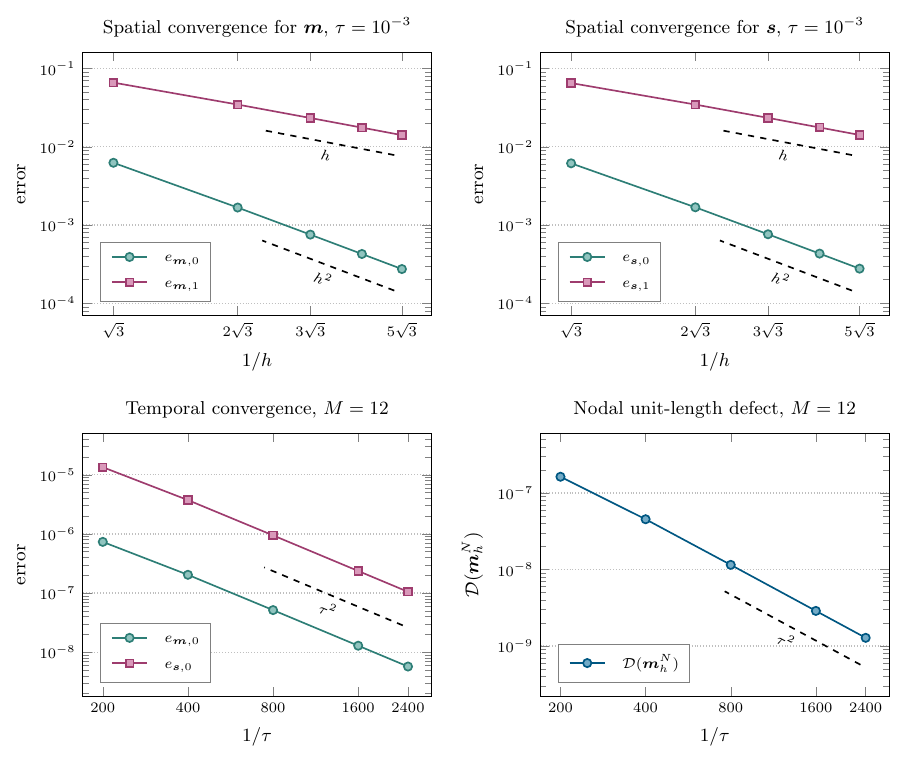}
	\caption{
		Spatial and temporal convergence of
		Algorithm~\ref{alg: full discr}, together
		with the final-time nodal unit-length defect~$\mathcal D(\mm_h^N)$ for the
		manufactured-solution experiment.
	}
	\label{fig:SDLLG convergence MMS}
\end{figure}

Finally, we record the {nodal unit-length defect
\begin{equation}\label{equ:nodal defect}
	\mathcal D(\mm_h^n)
	\coloneqq
	\max_{\zz\in\mathcal N_h^\omega}
	\bigl|\abs{\mm_h^n(\zz)}^2-1\bigr|,
	\qquad n=0,\ldots,N,
\end{equation}
at the final time, i.e., $\mathcal D(\mm_h^N)$.}
As expected for the projection-free tangent-plane update, this
quantity is nonzero. The numerical results indicate
$\mathcal D(\mm_h^N)=O(\tau^2)$, see {the last column of Table~\ref{tab:temporal convergence MMS} and} Figure~\ref{fig:SDLLG convergence MMS}.

\begin{table}[!ht]
	\centering
	\small
	\begin{tabular}{c @{\qquad} cc @{\qquad} cc @{\qquad} cc}
		\toprule
		$N$
		& $e_{\mm,0}$ & rate
		& $e_{\si,0}$ & rate
		& {$\mathcal D(\mm_h^N)$} & {rate}
		\\
		\midrule
		2
		& $7.330{\times}10^{-7}$ & --
		& $1.343{\times}10^{-5}$ & --
		& $1.629{\times}10^{-7}$ & {--}
		\\
		4
		& $2.034{\times}10^{-7}$ & $1.85$
		& $3.731{\times}10^{-6}$ & $1.85$
		& $4.532{\times}10^{-8}$ & {$1.85$}
		\\
		8
		& $5.147{\times}10^{-8}$ & $1.98$
		& $9.443{\times}10^{-7}$ & $1.98$
		& $1.148{\times}10^{-8}$ & {$1.98$}
		\\
		16
		& $1.287{\times}10^{-8}$ & $2.00$
		& $2.361{\times}10^{-7}$ & $2.00$
		& $2.870{\times}10^{-9}$ & {$2.00$}
		\\
		24
		& $5.719{\times}10^{-9}$ & $2.00$
		& $1.049{\times}10^{-7}$ & $2.00$
		& $1.276{\times}10^{-9}$ & {$2.00$}
		\\
		\bottomrule
	\end{tabular}
	\caption{
		Temporal convergence and final-time nodal unit-length defect for
		$M=12$ and $T=0.01$.
	}
	\label{tab:temporal convergence MMS}
\end{table}

Figure~\ref{fig:SDLLG convergence MMS} summarizes the convergence
behavior on logarithmic scales. In particular, the $\bL^2$ and $\bH^1$ spatial errors follow the expected $O(h^2)$ and $O(h)$ reference slopes, respectively. Both temporal $\bL^2$-errors and the final-time nodal unit-length defect $\mathcal D(\mm_h^N)$ exhibit
second-order decay with respect to $\tau$.

\subsection{Current-driven propagation of a N\'eel wall}
\label{subsec:Neel-wall}

We next consider a physically motivated domain-wall motion experiment. The purpose of this test is to illustrate the current-driven dynamics of a domain wall, in particular highlighting the energy and nodal unit-length defect evolution for our BDF2 approximation.

We consider the conducting domain and ferromagnetic region
\[
\Omega=(0,6)\times(0,1)\times(0,1),
\qquad
\omega=(0,6)\times(0,1)\times \Bigl(\frac13 , \frac23\Bigr).
\]
The geometry represents a long ferromagnetic strip embedded in a larger conducting region. This experiment is motivated by numerical studies of current-driven N\'eel-wall dynamics based on coupled spin-diffusion--LLG models; see, e.g.,~\cite{AbeRug_etal15}.

We set the external field $\ff=0$ and include an easy-axis anisotropy along the longitudinal direction $\ee_1=(1,0,0)^\top$ by setting $\pp(\mm)= K(\mm\cdot\ee_1)\ee_1$.
The parameters are chosen as
\[
\alpha=0.5,
\quad
\const{C}{ex}=0.1,
\quad
K=0.5,
\quad
c=0.2,
\quad
D_0=0.05,
\quad
\beta=0.4,
\quad
\beta'=0.5.
\]
The parameters above are dimensionless and are selected to produce clearly observable current-driven dynamics.

The initial magnetization is a N\'eel wall centered at
$x_1=3$,
\begin{equation*}
	\mm_0(\xx)
	=
	\begin{pmatrix}
		\tanh\bigl((x_1-3)/\delta\bigr)
		\\[1mm]
		\sech\bigl((x_1-3)/\delta\bigr)
		\\[1mm]
		0
	\end{pmatrix},
	\qquad
	\text{where }\;
	\delta
	=
	\sqrt{\frac{\const{C}{ex}}{K}}
	=
	\sqrt{0.2}.
\end{equation*}
In particular, $|\mm_0|=1$ pointwise. The initial spin accumulation
is taken to be $\si_0=0$.
We impose the longitudinal current
\begin{equation*}
	\jj(t)
	=
	J(t)\ee_1,
	\qquad \text{where }\;
	J(t)
	=
	\begin{cases}
		20,
		& 0\leq t\leq 1.5,
		\\
		0,
		& 1.5 < t\leq 2.5.
	\end{cases}
\end{equation*}
Note that the current is switched off at $t=1.5$, while the computation is
continued until $T=2.5$ in order to observe the subsequent relaxation.

We use the fitted tetrahedral mesh obtained from a
$54\times 9 \times 9$ subdivision of $\Omega$, corresponding to the mesh
parameter $M=9$. Since $M$ is divisible by $3$, the interfaces
$\{x_3=\frac13\}$ and $\{x_3=\frac23\}$ coincide exactly with mesh planes.
The time-step size is $\tau=2.5\times 10^{-3}$.

We also monitor the micromagnetic energy
\begin{equation}
	\label{equ:Neel-wall-energy}
	\mathcal E(\mm_h^n)
	=
	\frac{\const{C}{ex}}{2}
	\norm{\nabla\mm_h^n}_{\bL^2(\omega)}^2
	+
	\frac{K}{2}
	\int_\omega
	\bigl(
	1-(\mm_{h}^n\cdot \ee_1)^2
	\bigr)
	\,\mathrm{d}x
\end{equation}
as well as the nodal unit-length defect $\mathcal D(\mm_h^n)$ from~\eqref{equ:nodal defect}.
{Note that~\eqref{equ:Neel-wall-energy} is the energy~\eqref{equ:energy low} with $\ff=\boldsymbol{0}$ and $\mathcal E_{\rm low}(\mm)\coloneqq\frac K2\int_\omega\bigl(1-(\mm\cdot\ee_1)^2\bigr)\d x$, which satisfies~\eqref{ass:assumption pi linear self-adjoint} for $\pp(\mm)=K(\mm\cdot\ee_1)\ee_1$.}
Since the system is driven by a nonzero current for $0\leq t\leq 1.5$,
the micromagnetic energy is not expected to be monotone during this time.

\begin{figure}[!htb]
	\centering
	\includegraphics[width=\textwidth]
	{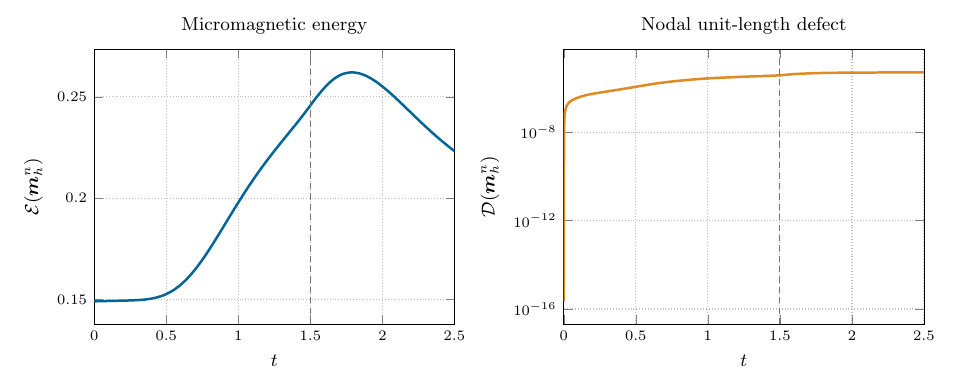}
	\caption{
		Evolution of the micromagnetic energy~\eqref{equ:Neel-wall-energy}
		and the nodal unit-length defect~\eqref{equ:nodal defect}.
		The applied current is switched off at $t=1.5$.
	}
	\label{fig:Neel-wall-diagnostics}
\end{figure}

Figure~\ref{fig:Neel-wall-diagnostics} displays the evolution of
$\mathcal E(\mm_h^n)$ and $\mathcal D(\mm_h^n)$. The current is active until $t= 1.5$ and is switched
off thereafter. In particular, the computation allows us to monitor
the dynamics of the system during the current-driven regime and the subsequent relaxation.

\begin{figure}[!ht]
	\centering
	\setlength{\tabcolsep}{3pt}
	\renewcommand{\arraystretch}{1.05}
	
	\begin{tabular}{@{}c@{\hspace{0.5em}}c@{\hspace{0.8em}}c@{}}
		&
		{\small\textbf{Magnetization $\mm_h^n$}}
		&
		{\small\textbf{Spin accumulation $\si_h^n$}}
		\\[1ex]
		
		{\small $t=0.5$}
		&
		\raisebox{-0.4\height}{%
			\includegraphics[
			width=0.45\textwidth
			]{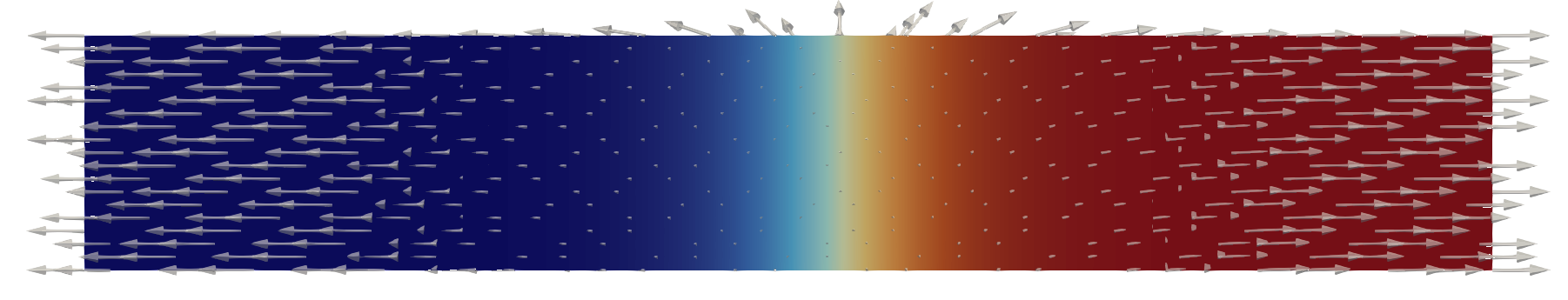}%
		}
		&
		\raisebox{-0.4\height}{%
			\includegraphics[
			width=0.45\textwidth
			]{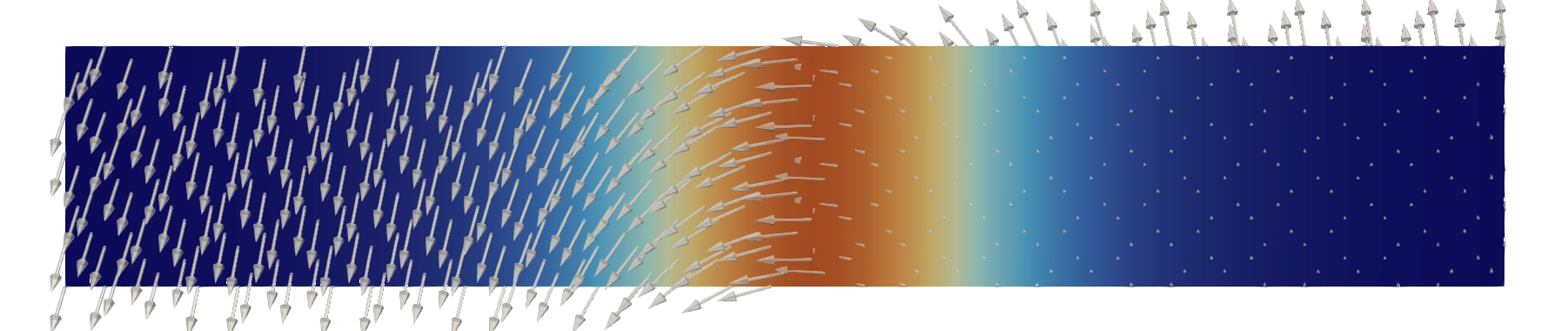}%
		}
		\\[1ex]
		
		{\small $t=1.0$}
		&
		\raisebox{-0.4\height}{%
			\includegraphics[
			width=0.45\textwidth
			]{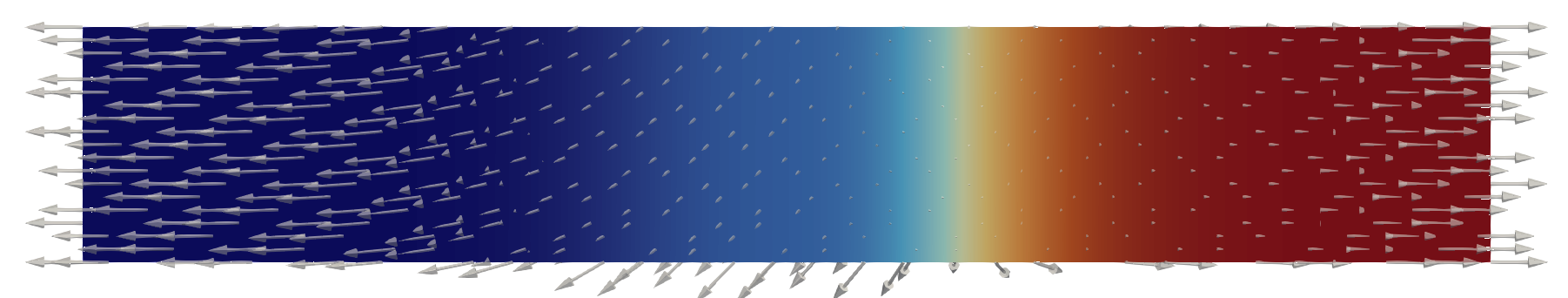}%
		}
		&
		\raisebox{-0.4\height}{%
			\includegraphics[
			width=0.45\textwidth
			]{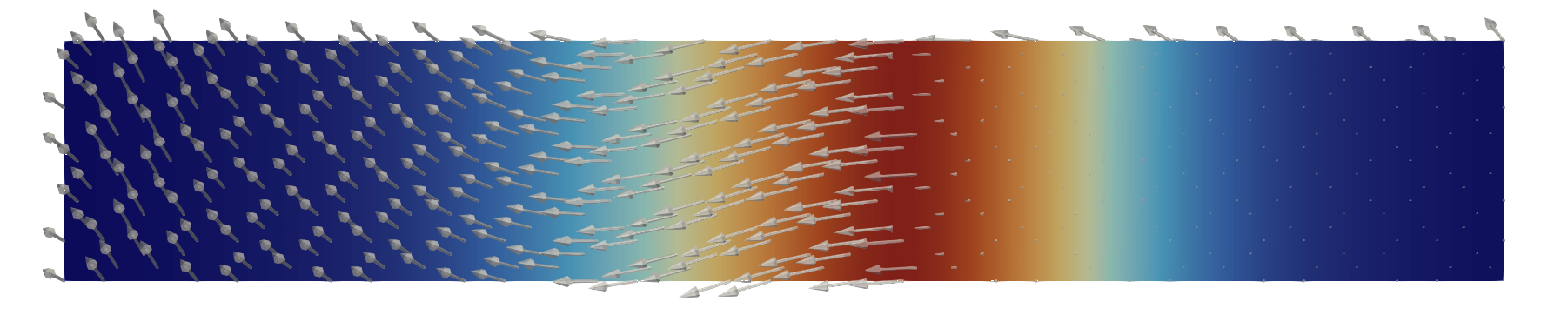}%
		}
		\\[1ex]
		
		{\small $t=1.5$}
		&
		\raisebox{-0.4\height}{%
			\includegraphics[
			width=0.45\textwidth
			]{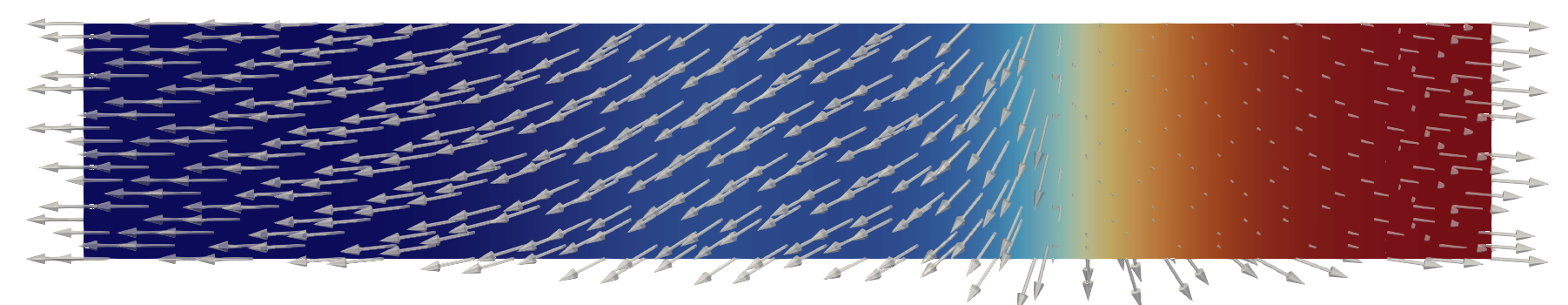}%
		}
		&
		\raisebox{-0.4\height}{%
			\includegraphics[
			width=0.45\textwidth
			]{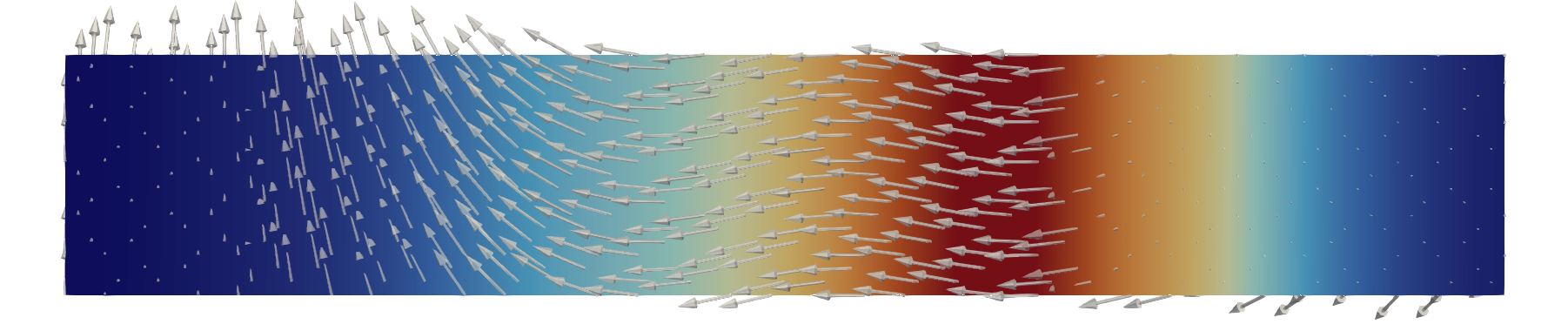}%
		}
		\\[1ex]
		
		{\small $t=2.5$}
		&
		\raisebox{-0.4\height}{%
			\includegraphics[
			width=0.45\textwidth
			]{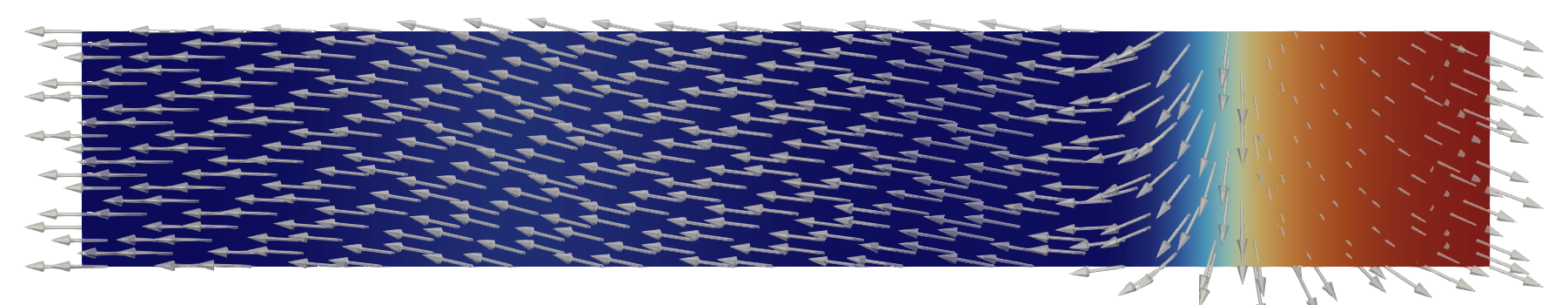}%
		}
		&
		\raisebox{-0.4\height}{%
			\includegraphics[
			width=0.45\textwidth
			]{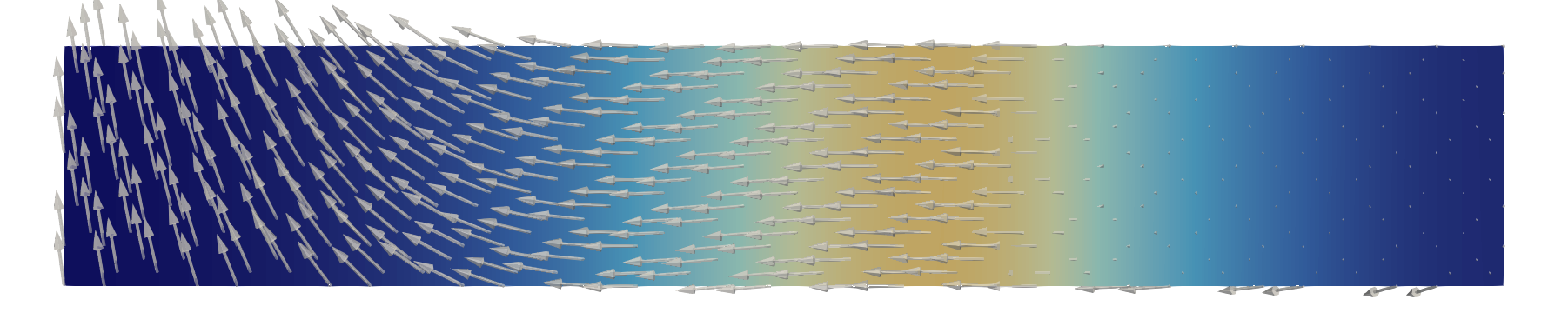}%
		}
		\\[1ex]
		
		&
		\includegraphics[
		width=0.37\textwidth
		]{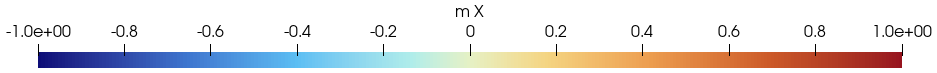}
		&
		\includegraphics[
		width=0.37\textwidth
		]{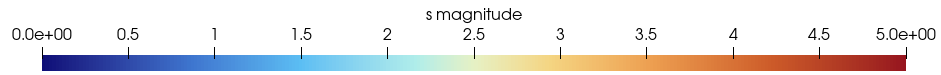}
	\end{tabular}
	
	\caption{
		Snapshots of the current-driven N\'eel-wall experiment on the
		midplane $x_3= \frac12$.
		In the left column, the background color represents
		$m_{h,1}$ and the arrows indicate the direction of $\mm_h$.
		In the right column, the background color represents
		$|\si_h|$ and the arrows indicate the direction of $\si_h$.
		The current is active for $0\leq t\leq 1.5$ and is
		switched off thereafter.
	}
	\label{fig:Neel-wall-snapshots}
\end{figure}

Figure~\ref{fig:Neel-wall-snapshots} shows snapshots
of the coupled dynamics on the midplane $x_3= \frac12$. For the
magnetization, the background color represents the longitudinal
component $m_{h,1}$, while the arrows indicate the direction of
$\mm_h$. For the spin accumulation, the background color represents
$|\si_h|$, while the arrows indicate the direction of $\si_h$.

The applied current drives the N\'eel wall towards the end of the
ferromagnetic strip. At $t=0.5$ and $t=1$, the wall has undergone a
clear longitudinal displacement, while the spin accumulation develops
a broader profile concentrated around and lagging behind the moving
magnetic texture. At $t=1.5$, the current is switched off, and the remaining magnetization and spin-accumulation undergo relaxation.

\subsection{Current-induced switching in a multilayer structure}
\label{subsec:multilayer-switching}

As a final experiment, we consider a simplified current-perpendicular-to-plane
multilayer switching problem motivated by the Permalloy spin-valve
configuration of~\cite{AbeRug_etal15}. This experiment illustrates the current-induced transfer of spin angular momentum between two ferromagnetic layers separated by a non-magnetic spacer, which is the main principle underlying spin-transfer torque magnetoresistive RAM (STT-MRAM)~\cite{FioEnd_etal21}.

\begin{figure}[!htb]
	\centering
	\begin{minipage}[c]{0.54\textwidth}
		\centering
		\small
		\renewcommand{\arraystretch}{1.4}
		
		\begin{tabular}{l|c}
			\textbf{Layer} & \textbf{$x_3$-interval} \\
			\hline
			Bottom contact & $(0,0.6)$ \\
			Fixed ferromagnetic layer $\omega_{\rm fix}$
			& $(0.6,2.6)$ \\
			Non-magnetic spacer & $(2.6,3.2)$ \\
			Free ferromagnetic layer $\omega_{\rm free}$
			& $(3.2,3.6)$ \\
			Top contact & $(3.6,4.2)$
		\end{tabular}
	\end{minipage}
	\hfill
	\begin{minipage}[c]{0.44\textwidth}
		\centering
		\includegraphics[width=\linewidth]
		{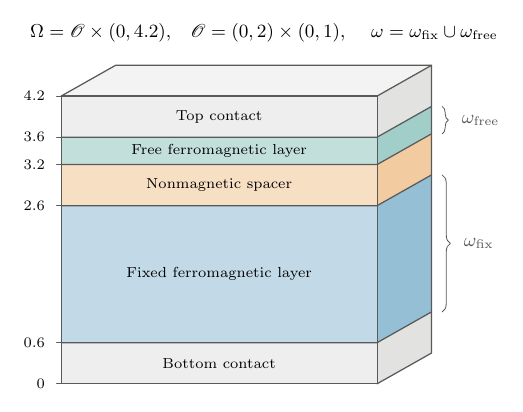}
	\end{minipage}
	\caption{
		Geometry of the multilayer spin-valve configuration.
		The conducting region consists of two ferromagnetic layers
		separated by a non-magnetic spacer and enclosed between
		two conducting contacts.
	}
	\label{fig:multilayer-geometry}
\end{figure}

The multilayer is described by Figure~\ref{fig:multilayer-geometry}. Fix the conducting region $\Omega= \mathscr{O} \times (0, 4.2)$, where $\mathscr{O} \coloneqq (0,2) \times (0,1)$.
The ferromagnetic region is the disconnected set $\omega= \omega_{\rm fix}\cup\omega_{\rm free}$, where~$\omega_{\rm fix} \coloneqq \mathscr{O} \times(0.6,2.6)$ and $\omega_{\rm free} \coloneqq \mathscr{O} \times(3.2,3.6)$.
While the layer thicknesses are proportional to those in the experiment of~\cite{AbeRug_etal15}, the present computation is dimensionless and is
intended only as a qualitative illustration.
The domain is discretized by a uniform $10\times 5\times 21$ subdivision, so that the mesh spacing in each coordinate direction is $0.2$ and all material interfaces coincide exactly with mesh
planes. The final time is set to be $T=4$.

We neglect the nonlocal magnetostatic field and instead represent its dominant shape-anisotropy effects by the local energy
\begin{align}
	\mathcal E_{\rm ani}(\mm)
	&=
	\frac12
	\int_{\omega}
	\left(
	K_2(\xx)m_2^2
	+
	K_3(\xx)m_3^2
	\right)
	\d x
	\label{equ:multilayer-anisotropy-energy}
\end{align}
corresponding to the effective field $\pp(\mm) = -K_2m_2\ee_2 -K_3m_3\ee_3$, where $m_j=\mm\cdot \ee_j$, $j=1,2,3$.
The coefficients are chosen piecewise as
\[
(K_2,K_3)
=
\begin{cases}
	(20,50),
	&
	\text{in }\omega_{\rm fix},
	\\
	(0.15,2),
	&
	\text{in }\omega_{\rm free}.
\end{cases}
\]
Hence, $\pm\ee_1$ are the preferred magnetic orientations, while the
out-of-plane direction $\ee_3$ is strongly penalized. A substantially larger anisotropy in $\omega_{\rm fix}$ models the pinning of the reference layer.

The remaining parameters are
\[
\alpha=0.2,
\qquad
\const{C}{ex}=0.1,
\qquad
c=1,
\qquad
\beta=0.9,
\qquad
\beta'=0.8,
\]
and the spin-diffusion coefficient is
\[
D_0(\xx)
=
\begin{cases}
	0.05,
	&
	\xx\in\omega,
	\\
	0.25,
	&
	\xx\in\Omega\setminus\omega.
\end{cases}
\]
In particular, the spin-diffusion coefficient in the non-magnetic material
is five times that in the ferromagnetic layers.

Initially, the fixed layer is magnetized along the easy axis,
\(
\mm_0=\ee_1
\)
in
\(
\omega_{\rm fix},
\)
while the free layer is initially almost parallel to the fixed layer,
but with a small in-plane perturbation,
\[
\mm_0
=
\begin{pmatrix}
	\cos\theta_0\\
	\sin\theta_0\\
	0
\end{pmatrix},
\qquad
\theta_0=12^\circ
\qquad
\text{in }\omega_{\rm free}.
\]
The perturbation breaks the symmetry of the exactly parallel
configuration. The initial spin accumulation is set as $\si_0 =0$.

A current is applied perpendicular to the layers, namely
\begin{align}
	\jj(t)
	&=
	J(t)\ee_3,
	\qquad \text{where } \;
	J(t)
	=
	\begin{cases}
		-10,
		&
		0\le t\le3,
		\\
		0,
		&
		3<t\le4.
	\end{cases}
	\label{equ:multilayer-current}
\end{align}
Physically, the fixed ferromagnet acts as a reference or polarizing layer: the
current-driven spin accumulation generated through the magnetic stack diffuses across the nonmagnetic spacer and produces a spin-transfer torque on the free layer.

\subsubsection{Switching dynamics}

\begin{figure}[!htb]
	\centering
	\includegraphics[width=\textwidth]
	{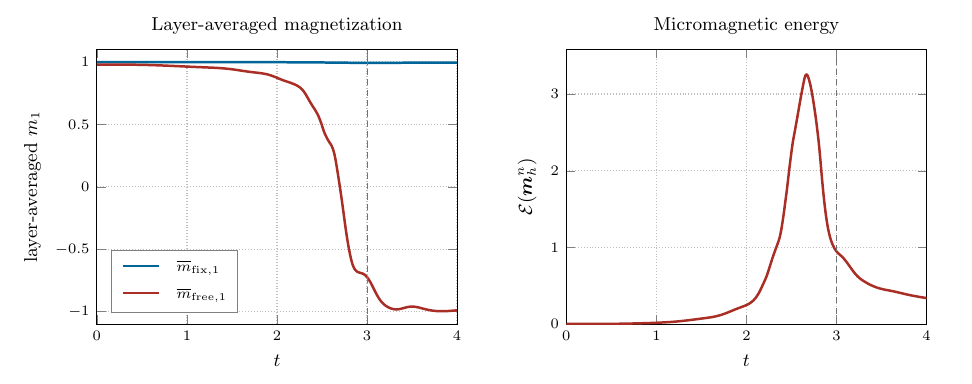}
	\caption{
		Current-induced switching experiment. The current $J=-10$ is applied for $0\le t\le 3$ and is switched off thereafter. Left: the free layer reverses from an initially almost parallel state to an antiparallel state, whereas the more strongly pinned fixed layer remains close to its initial orientation. Right: micromagnetic energy evolution.
	}
	\label{fig:multilayer-switching-diagnostics}
\end{figure}

To quantify the switching process, we consider the layer-averaged
magnetizations
\begin{align}
	\overline{\mm}_{\rm fix}(t)
	&\coloneqq
	\frac{1}{|\omega_{\rm fix}|}
	\int_{\omega_{\rm fix}}
	\mm_h(t)\d \xx,
	&
	\overline{\mm}_{\rm free}(t)
	&\coloneqq
	\frac{1}{|\omega_{\rm free}|}
	\int_{\omega_{\rm free}}
	\mm_h(t)\d \xx.
	\label{equ:layer-averaged-magnetisations}
\end{align}
Figure~\ref{fig:multilayer-switching-diagnostics} shows the longitudinal
components of these averages together with the magnetic energy
\begin{align}\label{equ:multilayer-discrete-energy}
	\mathcal E(\mm_h^n)
	&=
	\frac{\const{C}{ex}}{2}
	\norm{\nabla\mm_h^n}_{\bL^2(\omega)}^2
	+
	\mathcal E_{\rm ani}(\mm_h^n),
\end{align}
where $\mathcal E_{\rm ani}$ is given by \eqref{equ:multilayer-anisotropy-energy}. {Thus,~\eqref{equ:multilayer-discrete-energy} is the energy~\eqref{equ:energy low} with $\ff=\boldsymbol{0}$ and $\mathcal E_{\rm low}=\mathcal E_{\rm ani}$, which satisfies~\eqref{ass:assumption pi linear self-adjoint}.}

In this experiment, let $\tau= 2\times 10^{-3}$. For the sign convention in
\eqref{equ:multilayer-current}, the applied current drives the system
from the initially almost parallel configuration towards an
anti-parallel configuration. At the current switch-off time $t=3$, the free layer relaxes further towards the anti-parallel state after the current is removed. In contrast, the fixed layer remains close to its initial orientation throughout the computation. This confirms that the stronger anisotropy in the fixed layer provides an effective magnetic reference while the softer free layer undergoes reversal. 

The magnetic energy increases strongly during the switching. It subsequently decreases as the system approaches the reversed state and, after the current is switched off at $t=3$, continues to relax under Gilbert damping.
The maximum nodal unit-length defect satisfies
\[
\max_{0\le n\le N} \mathcal D(\mm_h^n)
=
\max_{0\le n\le N}
\max_{\zz\in\mathcal N_h^\omega}
\left|
|\mm_h^n(\zz)|^2-1
\right|=\mathcal D(\mm_h^N)
\simeq 6\times10^{-2}.
\]
{Here, the last equality follows from~\eqref{equ:explicit nodal BDF2 constraint}, which shows that $\abs{\mm_h^n(\zz)}^2-1\geq0$ is nondecreasing in~$n$ for every node $\zz\in\mathcal N_h^\omega$, since all coefficients on its right-hand side are nonnegative and nondecreasing in~$n$. Thus, for this time-step size, the projection-free approximation remains close to the unit sphere throughout the switching, whereas the defect becomes larger for larger time steps; see Figure~\ref{fig:multilayer-nodal}.}

\subsubsection{Nodal unit-length defect for various $\tau$}

\begin{figure}[!htb]
	\centering
	\includegraphics[width=\textwidth]
	{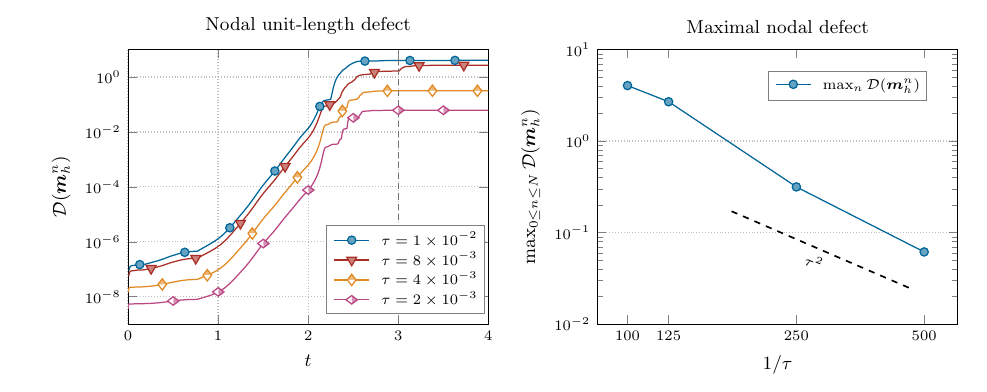}
	\caption{
		Nodal unit-length defect in the current-induced multilayer switching
		experiment. Left: evolution of
		$\mathcal D(\mm_h^n)=\max_{\zz\in\mathcal N_h^\omega}
		\bigl||\mm_h^n(\zz)|^2-1\bigr|$ 
		for different time-step sizes.
		Right: the maximum defect $\max_{0\le n\le N}\mathcal D(\mm_h^n)$ as a function of
		$1/\tau$. The dashed reference line indicates second-order decay
		with respect to $\tau$.
		}
	\label{fig:multilayer-nodal}
\end{figure}

We now let $\tau$ vary.
Figure~\ref{fig:multilayer-nodal} presents the nodal unit-length violation of the projection-free scheme. For all tested time-step sizes, the defect remains small initially and increases most noticeably during the rapid switching stage. Refining the time step systematically reduces the defect over the whole simulation. {The right panel shows that $\max_{0\le n\le N}\mathcal D(\mm_h^n)=\mathcal D(\mm_h^N)$ decays roughly like $\tau^2$ over the considered range of time-step sizes, in agreement with Remark~\ref{rem:quadratic defect}.}

\appendix

\section{Identities and inequalities}\label{sec:appendix}

We collect some algebraic identities and inequalities which are used throughout the analysis.


\begin{lemma}[{\cite[Lemma~3.1]{AldFeiPra26}}]\label{lem:BDF2 identities}
	The discrete solutions generated by Algorithm~\ref{alg: full discr} satisfy
	\begin{align}
		\vv_h^0
		&=
		\dtau\mm_h^1,
		\label{equ:BDF2 id 1}
		\\
		\intertext{{and, for $n=1,\ldots,N-1$,}}
		\vv_h^n
		&=
		\Dtau\mm_h^{n+1},
		\label{equ:BDF2 id 2}
		\\
		2\vv_h^n
		&=
		3\dtau\mm_h^{n+1}
		- \dtau\mm_h^n,
		\label{equ:BDF2 id 3}
		\\
		2\left(
		\vv_h^n-\dtau\mm_h^{n+1}
		\right)
		&=
		\tau\dtau^2\mm_h^{n+1}.
		\label{equ:BDF2 id 4}
		\end{align}
\end{lemma}


\begin{lemma}[BDF2 norm equivalence] \label{lem:BDF2 norm equivalence}
	Let $X$ be a Hilbert space and let $\{\ww^n\}_{n=0}^N\subset X$.
	Then, for every $n=1,\ldots,N$,
	\begin{align}
		\label{equ:BDF2 norm equivalence}
		C_1
		\bigg(
		\tau\norm[\big]{\dtau\ww^1}_{X}^2
		+
		\tau\sum_{j=1}^{n-1}
		\norm[\big]{\Dtau\ww^{j+1}}_{X}^2
		\bigg)^{1/2}
		&\leq
		\bigg(
		\tau\sum_{j=0}^{n-1}
		\norm[\big]{\dtau\ww^{j+1}}_{X}^2
		\bigg)^{1/2}
		\nonumber
		\\
		&\mkern-50mu\leq
		C_2
		\bigg(
		\tau\norm[\big]{\dtau\ww^1}_{X}^2
		+
		\tau\sum_{j=1}^{n-1}
		\norm[\big]{\Dtau\ww^{j+1}}_{X}^2
		\bigg)^{1/2},
	\end{align}
	where $C_1={\frac12}$ and $C_2={\frac{\sqrt5}{2}}$ {are sharp}.
\end{lemma}

\begin{proof}
	Set $\qq_0:=\dtau\ww^1$ and $\qq_j:=\Dtau\ww^{j+1}$ for all $j\geq 1$.
	Then, by the definition of $\Dtau$,
	\[
	\qq_j
	=
	\frac32\dtau\ww^{j+1}
	-
	\frac12\dtau\ww^j,
	\qquad j\geq1.
	\]
	Hence, {since $\norm[\big]{\qq_j}_X^2\leq 3\norm[\big]{\dtau\ww^{j+1}}_X^2+\norm[\big]{\dtau\ww^j}_X^2$ for $j\geq1$ by the Cauchy--Schwarz inequality,}
	\begin{equation}\label{equ:BDF2 norm up}
	\tau\sum_{j=0}^{n-1}\norm[\big]{\qq_j}_{X}^2
	\leq
	{4}\tau\sum_{j=0}^{n-1}
	\norm[\big]{\dtau\ww^{j+1}}_{X}^2.
	\end{equation}
	Conversely, we have $\dtau\ww^1=\qq_0$ and, for $j\geq1$,
	\[
	\dtau\ww^{j+1}
	=
	\frac23\qq_j+\frac13\dtau\ww^j,
	\]
	and thus, {by Young's inequality with parameter $5/4$,}
	\[
	\tau\sum_{j=0}^{n-1}
	\norm[\big]{\dtau\ww^{j+1}}_{X}^2
	\leq
	\tau\sum_{j=0}^{n-1}\norm[\big]{\qq_j}_{X}^2
	+
	{\frac15} \tau \sum_{j=0}^{n-1} \norm[\big]{\dtau\ww^{j+1}}_{X}^2.
	\]
	Absorbing the last term and combining the result with \eqref{equ:BDF2 norm up} proves~\eqref{equ:BDF2 norm equivalence}.
\end{proof}


\begin{lemma}[\texorpdfstring{\cite[Proposition~3.8]{AldFeiPra26}}{}]\label{lem:explicit nodal BDF2 constraint}
	For every $n\geq 1$ and $\vec{z}\in \mathcal{N}_h^\omega$,
	\begin{align}
		\label{equ:explicit nodal BDF2 constraint}
		\abs{\mm_h^n(\vec{z})}^2
		=
		1 + \frac32
		\left( 1-3^{-n} \right) \tau^2 \abs{\vv_h^0(\vec{z})}^2
		+
		\frac32\tau^4
		\sum_{j=2}^{n}
		\left(1-3^{-(n+1-j)}\right) \abs{\dtau^2\mm_h^j(\vec{z})}^2.
	\end{align}
\end{lemma}

\begin{proof}
We refer to \cite[Proposition~3.8]{AldFeiPra26}. Notice that the only difference is that \cite[Proposition~3.8]{AldFeiPra26} states and proves this result for $n\geq 2$. However, if the sum is interpreted as an empty sum when $n=1$, it holds that $|\mm_h^1(\vec{z})|^2 = |\mm_h^0(\vec{z})|^2+\tau^2|\vv_h^0(\vec{z})|^2$, and~\eqref{equ:explicit nodal BDF2 constraint} coincides with~\eqref{equ:first nodal constraint}. This proves the lemma.
\end{proof}

\begin{lemma}[Discrete Gronwall inequality \texorpdfstring{\cite[Lemma~10.5]{thomee2013galerkin}}{}]\label{lem: gronwall}
Let $\left\{a_n\right\}_{n \geq 0}$, $\left\{b_n\right\}_{n \geq 0}$, and $\left\{w_n\right\}_{n \geq 0}$ be sequences of real numbers such that
	\[
0 \leq a_n \leq a_{n+1}, \quad b_n \geq 0, \quad \text { and } \quad w_n \leq a_n+\sum_{j=0}^{n-1} b_j w_j \quad \text { for all } n \geq 0 .
\]
Then, it holds that
\begin{equation*}
	w_n \leq a_n \exp \sum_{j=0}^{n-1} b_j \quad \text { for all } n \geq 0.\tag*{\qed}
\end{equation*}
\end{lemma}

\printbibliography

@Article {AbeHrk_etal14,
	AUTHOR = {Abert, Claas and Hrkac, Gino and Page, Marcus and Praetorius,
	Dirk and Ruggeri, Michele and Suess, Dieter},
	TITLE = {Spin-polarized transport in ferromagnetic multilayers: an
	unconditionally convergent {FEM} integrator},
	JOURNAL = {Comput. Math. Appl.},
	FJOURNAL = {Computers \& Mathematics with Applications. An International
	Journal},
	VOLUME = {68},
	YEAR = {2014},
	NUMBER = {6},
	PAGES = {639--654},
	ISSN = {0898-1221,1873-7668},
	MRCLASS = {65M60 (35B45 35D30 35K51 35Q60 65M12 82D40)},
	MRNUMBER = {3253344},
	DOI = {10.1016/j.camwa.2014.07.010},
}

@article{AbeRug_etal16,
	title={A self-consistent spin-diffusion model for micromagnetics},
	author={Abert, Claas and Ruggeri, Michele and Bruckner, Florian and Vogler, Christoph and Manchon, Aurelien and Praetorius, Dirk and Suess, Dieter},
	JOURNAL = {Scientific Reports},
	volume={6},
	number={1},
	pages={16},
	year={2016},
	DOI = {10.1038/s41598-016-0019-y},
}

@article{AbeRug_etal15,
	title={A three-dimensional spin-diffusion model for micromagnetics},
	author={Abert, Claas and Ruggeri, Michele and Bruckner, Florian and Vogler, Christoph and Hrkac, Gino and Praetorius, Dirk and Suess, Dieter},
	JOURNAL = {Scientific Reports},
	volume={5},
	number={1},
	pages={14855},
	year={2015},
	DOI = {10.1038/srep14855},
}

@Article {GarWan07,
	AUTHOR = {Garc\'ia-Cervera, Carlos J. and Wang, Xiao-Ping},
	TITLE = {Spin-polarized currents in ferromagnetic multilayers},
	JOURNAL = {J. Comput. Phys.},
	FJOURNAL = {Journal of Computational Physics},
	VOLUME = {224},
	YEAR = {2007},
	NUMBER = {2},
	PAGES = {699--711},
	ISSN = {0021-9991,1090-2716},
	MRCLASS = {82D40},
	MRNUMBER = {2330291},
	DOI = {10.1016/j.jcp.2006.10.029},
}

@article{AlnaesEtal15,
	title     = {The {FEniCS} Project Version 1.5},
	author    = {Aln{\ae}s, Martin S. and Blechta, Jan and Hake, Johan and Johansson, August and Kehlet, Benjamin and Logg, Anders and Richardson, Chris N. and Ring, Johannes and Rognes, Marie E. and Wells, Garth N.},
	journal   = {Archive of Numerical Software},
	year      = {2015},
	volume    = {3},
	number    = {100},
	pages     = {9--23},
	DOI       = {10.11588/ans.2015.100.20553},
}

@book{thomee2013galerkin,
  author    = {Vidar Thom{\'e}e},
  title     = {Galerkin finite element methods for parabolic problems},
  edition   = {2},
  series    = {Springer Series in Computational Mathematics},
  volume    = {25},
  publisher = {Springer},
  address   = {Berlin, Heidelberg},
  year      = {2006},
  isbn      = {978-3-540-33122-3},
  doi       = {10.1007/3-540-33122-0}
}

@article {GarCerWan_dcds_07,
    AUTHOR = {Garc\'ia-Cervera, Carlos J. and Wang, Xiao-Ping},
     TITLE = {Spin-polarized transport: existence of weak solutions},
   JOURNAL = {Discrete Contin. Dyn. Syst. Ser. B},
  FJOURNAL = {Discrete and Continuous Dynamical Systems. Series B. A Journal
              Bridging Mathematics and Sciences},
    VOLUME = {7},
      YEAR = {2007},
    NUMBER = {1},
     PAGES = {87--100},
      ISSN = {1531-3492,1553-524X},
   MRCLASS = {82D40 (35D05 35Q60)},
  MRNUMBER = {2257453},
       DOI = {10.3934/dcdsb.2007.7.87},
       URL = {https://doi.org/10.3934/dcdsb.2007.7.87},
}

@book {BreSco08,
    AUTHOR = {Brenner, Susanne C. and Scott, L. Ridgway},
     TITLE = {The mathematical theory of finite element methods},
    SERIES = {Texts in Applied Mathematics},
    VOLUME = {15},
   EDITION = {Third},
 PUBLISHER = {Springer, New York},
      YEAR = {2008},
     PAGES = {xviii+397},
      ISBN = {978-0-387-75933-3},
   MRCLASS = {65-01 (65-02)},
  MRNUMBER = {2373954},
       DOI = {10.1007/978-0-387-75934-0},
}

@article{Simon86,
  author  = {Simon, Jacques},
  title   = {Compact sets in the space $L^p(0,T;B)$},
  journal = {Ann. Mat. Pura Appl. (4)},
  volume  = {146}, pages = {65--96}, year = {1986},
  doi     = {10.1007/BF01762360},
}

@misc{AldFeiPraSoen26_part3,
    AUTHOR = {Ald{\'e}, Michele and Praetorius, Dirk and Soenjaya, Agus L.},
     TITLE = {{BDF2-type} integrator for {Landau--Lifshitz--Gilbert} equation in micromagnetics: {IMEX} time-stepping and full effective field},
   	year 		  = {2026},	
   archiveprefix = {arXiv},
   eprint 		  = {2610.10291},
   url = {https://arxiv.org/abs/2610.10291}
}

@article {Soe26,
	AUTHOR = {Soenjaya, Agus L.},
	TITLE = {The {L}andau-{L}ifshitz-{B}loch equation with spin diffusion:
	global strong solution and finite element approximation},
	JOURNAL = {Numer. Methods Partial Differential Equations},
	FJOURNAL = {Numerical Methods for Partial Differential Equations. An
	International Journal},
	VOLUME = {42},
	YEAR = {2026},
	NUMBER = {1},
	PAGES = {Paper No. e70070, 24},
	ISSN = {0749-159X,1098-2426},
	MRCLASS = {65M60 (35A01 35Q60)},
	MRNUMBER = {5018588},
	DOI = {10.1002/num.70070},
	URL = {https://doi.org/10.1002/num.70070},
}

@thesis{Rug16,
   author      = {Ruggeri, Michele},
   title       = {Coupling and numerical integration of the
                  {Landau--Lifshitz--Gilbert} equation},
   type        = {PhD thesis},
   institution = {TU Wien},
   year        = {2016},
}

@article{FioEnd_etal21,
	title = {Coupled spin and charge drift-diffusion approach applied to magnetic tunnel junctions},
	journal = {Solid-State Electronics},
	volume = {186},
	pages = {108103},
	year = {2021},
	issn = {0038-1101},
	doi = {10.1016/j.sse.2021.108103},
	url = {https://www.sciencedirect.com/science/article/pii/S0038110121001489},
	author = {S. Fiorentini and J. Ender and S. Selberherr and R.L. {de Orio} and W. Goes and V. Sverdlov},
}

@article {BruSue_etal14,
	AUTHOR = {Bruckner, F. and Suess, D. and Feischl, M. and F\"uhrer, T.
	and Goldenits, P. and Page, M. and Praetorius, D. and Ruggeri,
	M.},
	TITLE = {Multiscale modeling in micromagnetics: existence of solutions
	and numerical integration},
	JOURNAL = {Math. Models Methods Appl. Sci.},
	FJOURNAL = {Mathematical Models and Methods in Applied Sciences},
	VOLUME = {24},
	YEAR = {2014},
	NUMBER = {13},
	PAGES = {2627--2662},
	ISSN = {0218-2025,1793-6314},
	MRCLASS = {78A25 (35Q60 65N30)},
	MRNUMBER = {3260281},
	DOI = {10.1142/S0218202514500328},
}

@article{WorledgeHu2024,
	author  = {Worledge, Daniel C. and Hu, Guohan},
	title   = {Spin-transfer torque magnetoresistive random access memory technology status and future directions},
	journal = {Nature Reviews Electrical Engineering},
	volume  = {1},
	pages   = {730--747},
	year    = {2024},
	doi     = {10.1038/s44287-024-00111-z}
}

@article{FiorentiniEtAl2022,
	author  = {Fiorentini, Simone and Bendra, Mario and Ender, Johannes and de Orio, Roberto L. and Goes, Wolfgang and Selberherr, Siegfried and Sverdlov, Viktor},
	title   = {Spin and charge drift-diffusion in ultra-scaled {MRAM} cells},
	journal = {Scientific Reports},
	volume  = {12},
	pages   = {20958},
	year    = {2022},
	doi     = {10.1038/s41598-022-25586-4}
}

@article{BendraEtAl2026,
	author  = {Bendra, Mario and Goes, Wolfgang and Selberherr, Siegfried and Sverdlov, Viktor},
	title   = {Simulation of {SAF}-enhanced multilayered {STT-MRAM} structures},
	journal = {Microelectronic Engineering},
	volume  = {302},
	pages   = {112426},
	year    = {2026},
	doi     = {10.1016/j.mee.2025.112426}
}

@article{ValFer93,
	title = {Theory of the perpendicular magnetoresistance in magnetic multilayers},
	author = {Valet, T. and Fert, A.},
	journal = {Phys. Rev. B},
	volume = {48},
	issue = {10},
	pages = {7099--7113},
	year = {1993},
	month = sep,
	publisher = {American Physical Society},
	doi = {10.1103/PhysRevB.48.7099},
	url = {https://link.aps.org/doi/10.1103/PhysRevB.48.7099}
}

@article{ZhaLi04,
	title = {Roles of Nonequilibrium Conduction Electrons on the Magnetization Dynamics of Ferromagnets},
	author = {Zhang, S. and Li, Z.},
	journal = {Phys. Rev. Lett.},
	volume = {93},
	issue = {12},
	pages = {127204},
	numpages = {4},
	year = {2004},
	month = sep,
	publisher = {American Physical Society},
	doi = {10.1103/PhysRevLett.93.127204},
	url = {https://link.aps.org/doi/10.1103/PhysRevLett.93.127204}
}

@article{Slo96,
	title = {Current-driven excitation of magnetic multilayers},
	journal = {J. Magn. Magn. Mater.},
	volume = {159},
	number = {1},
	pages = {L1--L7},
	year = {1996},
	issn = {0304-8853},
	doi = {10.1016/0304-8853(96)00062-5},
	url = {https://www.sciencedirect.com/science/article/pii/0304885396000625},
	author = {J.C. Slonczewski},
}

@article{ZhaLevFer02,
	title = {Mechanisms of Spin-Polarized Current-Driven Magnetization Switching},
	author = {Zhang, S. and Levy, P. M. and Fert, A.},
	journal = {Phys. Rev. Lett.},
	volume = {88},
	issue = {23},
	pages = {236601},
	numpages = {4},
	year = {2002},
	month = may,
	publisher = {American Physical Society},
	doi = {10.1103/PhysRevLett.88.236601},
	url = {https://link.aps.org/doi/10.1103/PhysRevLett.88.236601}
}

@article {PuGuo10,
	AUTHOR = {Pu, Xueke and Guo, Boling},
	TITLE = {Global smooth solutions for the one-dimensional spin-polarized
	transport equation},
	JOURNAL = {Nonlinear Anal.},
	FJOURNAL = {Nonlinear Analysis. Theory, Methods \& Applications. An
	International Multidisciplinary Journal},
	VOLUME = {72},
	YEAR = {2010},
	NUMBER = {3-4},
	PAGES = {1481--1487},
	ISSN = {0362-546X,1873-5215},
	MRCLASS = {35Q60},
	MRNUMBER = {2577549},
	DOI = {10.1016/j.na.2009.08.032},
	URL = {https://doi.org/10.1016/j.na.2009.08.032},
}

@article {GuoPu08,
	AUTHOR = {Guo, Boling and Pu, Xueke},
	TITLE = {Global smooth solutions of the spin polarized transport equation},
	JOURNAL = {Electron. J. Differential Equations},
	FJOURNAL = {Electronic Journal of Differential Equations},
	VOLUME = {2008},
	YEAR = {2008},
	NUMBER = {63},
	PAGES = {1--15},
	URL = {https://ejde.math.txstate.edu/Volumes/2008/63/abstr.html},
}

@article {AloJai06,
	AUTHOR = {Alouges, Fran\c{c}ois and Jaisson, Pascal},
	TITLE = {Convergence of a finite element discretization for the
	{L}andau-{L}ifshitz equations in micromagnetism},
	JOURNAL = {Math. Models Methods Appl. Sci.},
	FJOURNAL = {Mathematical Models and Methods in Applied Sciences},
	VOLUME = {16},
	YEAR = {2006},
	NUMBER = {2},
	PAGES = {299--316},
	ISSN = {0218-2025,1793-6314},
	MRCLASS = {65M60 (35K55 35Q60 82D40)},
	MRNUMBER = {2210092},
	MRREVIEWER = {B\'eatrice\ M.\ Rivi\`ere},
	DOI = {10.1142/S0218202506001169},
}

@article {Alo08,
	AUTHOR = {Alouges, Fran\c{c}ois},
	TITLE = {A new finite element scheme for {L}andau-{L}ifchitz equations},
	JOURNAL = {Discrete Contin. Dyn. Syst. Ser. S},
	FJOURNAL = {Discrete and Continuous Dynamical Systems. Series S},
	VOLUME = {1},
	YEAR = {2008},
	NUMBER = {2},
	PAGES = {187--196},
	ISSN = {1937-1632,1937-1179},
	MRCLASS = {65M60 (35K55 82B80 82D40)},
	MRNUMBER = {2379897},
	MRREVIEWER = {Etienne\ Emmrich},
	DOI = {10.3934/dcdss.2008.1.187},
	URL = {https://doi.org/10.3934/dcdss.2008.1.187},
}

@article {Ron23,
	AUTHOR = {Rong, Rong},
	TITLE = {Partial regular solution to the {L}andau-{L}ifshitz-{M}axwell
	system with spin accumulation in two dimensions},
	JOURNAL = {Math. Methods Appl. Sci.},
	FJOURNAL = {Mathematical Methods in the Applied Sciences},
	VOLUME = {46},
	YEAR = {2023},
	NUMBER = {4},
	PAGES = {4303--4330},
	ISSN = {0170-4214,1099-1476},
	MRCLASS = {35Q35},
	MRNUMBER = {4551367},
	DOI = {10.1002/mma.8757},
}

@article {CheWan23,
	AUTHOR = {Chen, Bo and Wang, Youde},
	TITLE = {Very regular solution to {L}andau-{L}ifshitz-{G}ilbert system
	with spin-polarized transport},
	JOURNAL = {Front. Math.},
	FJOURNAL = {Frontiers of Mathematics},
	VOLUME = {18},
	YEAR = {2023},
	NUMBER = {4},
	PAGES = {751--795},
	ISSN = {2731-8648,2731-8656},
	MRCLASS = {35K59 (35Q60 78A25)},
	MRNUMBER = {4619012},
	DOI = {10.1007/s11464-021-0209-1},
}

@article{FraJunPraSla23,
	author  = {Di Fratta, Giovanni and J{\"u}ngel, Ansgar and Praetorius, Dirk and Slastikov, Valeriy},
	title   = {Spin-diffusion model for micromagnetics in the limit of long times},
	journal = {Journal of Differential Equations},
	volume  = {343},
	pages   = {467--494},
	year    = {2023},
	doi     = {10.1016/j.jde.2022.10.012}
}

@article {BarPro06,
	AUTHOR = {Bartels, S\"oren and Prohl, Andreas},
	TITLE = {Convergence of an implicit finite element method for the
	{L}andau-{L}ifshitz-{G}ilbert equation},
	JOURNAL = {SIAM J. Numer. Anal.},
	FJOURNAL = {SIAM Journal on Numerical Analysis},
	VOLUME = {44},
	YEAR = {2006},
	NUMBER = {4},
	PAGES = {1405--1419},
	ISSN = {0036-1429,1095-7170},
	MRCLASS = {65M60 (82D40)},
	MRNUMBER = {2257110},
	MRREVIEWER = {Anne\ Nouri},
	DOI = {10.1137/050631070},
}

@article{RugAbeHrk_etal16,
	title = {Coupling of dynamical micromagnetism and a stationary spin drift-diffusion equation: A step towards a fully self-consistent spintronics framework},
	journal = {Physica B: Condensed Matter},
	volume = {486},
	pages = {88--91},
	year = {2016},
	issn = {0921-4526},
	doi = {10.1016/j.physb.2015.09.003},
	author = {Michele Ruggeri and Claas Abert and Gino Hrkac and Dieter Suess and Dirk Praetorius},
}

@article {AloKriSteTou14,
	AUTHOR = {Alouges, Fran\c{c}ois and Kritsikis, Evaggelos and Steiner,
	Jutta and Toussaint, Jean-Christophe},
	TITLE = {A convergent and precise finite element scheme for
	{L}andau-{L}ifschitz-{G}ilbert equation},
	JOURNAL = {Numer. Math.},
	FJOURNAL = {Numerische Mathematik},
	VOLUME = {128},
	YEAR = {2014},
	NUMBER = {3},
	PAGES = {407--430},
	ISSN = {0029-599X,0945-3245},
	MRCLASS = {65M60 (35K55 35Q60 65M12)},
	MRNUMBER = {3268842},
	MRREVIEWER = {Nicolae\ Pop},
	DOI = {10.1007/s00211-014-0615-3},
	URL = {https://doi.org/10.1007/s00211-014-0615-3},
}

@article {AkrFeiKovLub21,
	AUTHOR = {Akrivis, Georgios and Feischl, Michael and Kov\'{a}cs,
	Bal\'{a}zs and Lubich, Christian},
	TITLE = {Higher-order linearly implicit full discretization of the
	{L}andau-{L}ifshitz-{G}ilbert equation},
	JOURNAL = {Math. Comp.},
	FJOURNAL = {Mathematics of Computation},
	VOLUME = {90},
	YEAR = {2021},
	NUMBER = {329},
	PAGES = {995--1038},
	ISSN = {0025-5718,1088-6842},
	MRCLASS = {65M60 (35Q60 65L06 65M12 65M15)},
	MRNUMBER = {4232216},
	MRREVIEWER = {Hamdullah\ Y\"{u}cel},
	DOI = {10.1090/mcom/3597},
	URL = {https://doi.org/10.1090/mcom/3597},
}

@article{AldFeiPra26,
	author  = {Ald{\'e}, Michele and Feischl, Michael and Praetorius, Dirk},
	title   = {{BDF2}-type integrator for {Landau--Lifshitz--Gilbert}
	equation in micromagnetics: Unconditional weak convergence
	to weak solutions},
	journal = {Math. Comp.},
	year    = {2026},
	doi     = {10.1090/mcom/4210},
	url     = {https://doi.org/10.1090/mcom/4210}
}

@article {Dah78,
	AUTHOR = {Dahlquist, Germund},
	TITLE = {{$G$}-stability is equivalent to {$A$}-stability},
	JOURNAL = {BIT},
	FJOURNAL = {BIT. Nordisk Tidskrift for Informationsbehandling (BIT)},
	VOLUME = {18},
	YEAR = {1978},
	NUMBER = {4},
	PAGES = {384--401},
	ISSN = {0006-3835},
	MRCLASS = {65L05 (65J10)},
	MRNUMBER = {520750},
	MRREVIEWER = {Ramon\ E.\ Moore},
	DOI = {10.1007/BF01932018},
	URL = {https://doi.org/10.1007/BF01932018},
}

@article {NevOde81,
	AUTHOR = {Nevanlinna, Olavi and Odeh, F.},
	TITLE = {Multiplier techniques for linear multistep methods},
	JOURNAL = {Numer. Funct. Anal. Optim.},
	FJOURNAL = {Numerical Functional Analysis and Optimization. An
	International Journal},
	VOLUME = {3},
	YEAR = {1981},
	NUMBER = {4},
	PAGES = {377--423},
	ISSN = {0163-0563,1532-2467},
	MRCLASS = {65L05},
	MRNUMBER = {636736},
	MRREVIEWER = {Peter\ Alfeld},
	DOI = {10.1080/01630568108816097},
	URL = {https://doi.org/10.1080/01630568108816097},
}

@article {Bar05,
	AUTHOR = {Bartels, S\"oren},
	TITLE = {Stability and convergence of finite-element approximation
	schemes for harmonic maps},
	JOURNAL = {SIAM J. Numer. Anal.},
	FJOURNAL = {SIAM Journal on Numerical Analysis},
	VOLUME = {43},
	YEAR = {2005},
	NUMBER = {1},
	PAGES = {220--238},
	ISSN = {0036-1429,1095-7170},
	MRCLASS = {65N30 (58E20 76A15 82D30)},
	MRNUMBER = {2177142},
	MRREVIEWER = {R.\ Kodn\'ar},
	DOI = {10.1137/040606594},
	URL = {https://doi.org/10.1137/040606594},
}

@book {Bar15,
	AUTHOR = {Bartels, S\"oren},
	TITLE = {Numerical methods for nonlinear partial differential
	equations},
	SERIES = {Springer Series in Computational Mathematics},
	VOLUME = {47},
	PUBLISHER = {Springer, Cham},
	YEAR = {2015},
	PAGES = {x+393},
	ISBN = {978-3-319-13796-4; 978-3-319-13797-1},
	MRCLASS = {65-01 (35A15 35A35 65Mxx 65Nxx)},
	MRNUMBER = {3309171},
	MRREVIEWER = {Karsten\ Urban},
	DOI = {10.1007/978-3-319-13797-1},
	URL = {https://doi.org/10.1007/978-3-319-13797-1},
}

@article {PuWan20,
	AUTHOR = {Pu, Xueke and Wang, Wendong},
	TITLE = {Partial regularity to the {L}andau-{L}ifshitz equation with
	spin accumulation},
	JOURNAL = {J. Differential Equations},
	FJOURNAL = {Journal of Differential Equations},
	VOLUME = {268},
	YEAR = {2020},
	NUMBER = {2},
	PAGES = {707--737},
	ISSN = {0022-0396,1090-2732},
	MRCLASS = {35Q60 (35B65 58J35)},
	MRNUMBER = {4021901},
	DOI = {10.1016/j.jde.2019.08.034},
	URL = {https://doi.org/10.1016/j.jde.2019.08.034},
}

@article {Ron23D,
	AUTHOR = {Rong, Rong},
	TITLE = {Partial regular solution to the coupled spin polarized system
	in three dimensions},
	JOURNAL = {Discrete Contin. Dyn. Syst. Ser. B},
	FJOURNAL = {Discrete and Continuous Dynamical Systems. Series B. A Journal
	Bridging Mathematics and Sciences},
	VOLUME = {28},
	YEAR = {2023},
	NUMBER = {8},
	PAGES = {4293--4310},
	ISSN = {1531-3492,1553-524X},
	MRCLASS = {35Q35},
	MRNUMBER = {4561507},
	MRREVIEWER = {Taylan\ Sengul},
	DOI = {10.3934/dcdsb.2023013},
	URL = {https://doi.org/10.3934/dcdsb.2023013},
}

@article {AnMajProTra22,
	AUTHOR = {An, Xin and Majee, Ananta K. and Prohl, Andreas and Tran,
	Thanh},
	TITLE = {Optimal control for a coupled spin-polarized current and
	magnetization system},
	JOURNAL = {Adv. Comput. Math.},
	FJOURNAL = {Advances in Computational Mathematics},
	VOLUME = {48},
	YEAR = {2022},
	NUMBER = {3},
	PAGES = {Paper No. 28, 40},
	ISSN = {1019-7168,1572-9044},
	MRCLASS = {49J20 (49K20 49M25 49M41 78M30)},
	MRNUMBER = {4419351},
	MRREVIEWER = {Wei\ Gong},
	DOI = {10.1007/s10444-022-09947-1},
	URL = {https://doi.org/10.1007/s10444-022-09947-1},
}

@article{ShpLevZha03,
	title = {Self-consistent treatment of nonequilibrium spin torques in magnetic multilayers},
	author = {Shpiro, Asya and Levy, Peter M. and Zhang, Shufeng},
	journal = {Phys. Rev. B},
	volume = {67},
	issue = {10},
	pages = {104430},
	numpages = {17},
	year = {2003},
	month = {Mar},
	publisher = {American Physical Society},
	doi = {10.1103/PhysRevB.67.104430},
	url = {https://link.aps.org/doi/10.1103/PhysRevB.67.104430}
}

\end{document}